\documentclass{amsart}

\usepackage{amsmath}
\usepackage{amsfonts}
\usepackage{amssymb}
\usepackage{a4wide}
\usepackage{amsthm}
\usepackage{mathrsfs}
\usepackage{epic}
\usepackage{hyperref}
\usepackage{array,booktabs}
\usepackage{tikz}

\usepackage{mathtools}

\usepackage{pgfplots}
\usepackage{anyfontsize}
\pgfplotsset{compat=newest}
\usetikzlibrary{calc,intersections,datavisualization.formats.functions, patterns}
\usepgfplotslibrary{fillbetween, decorations.softclip}

\newtheorem{Theorem}{Theorem}

\theoremstyle{definition}

\def\komma{\raise.5ex\hbox{,}}
\def\punt{\raise.5ex\hbox{.}}
\newcommand{\R}{\mathbb R}
\newcommand{\Z}{\mathbb Z}
\newcommand{\Q}{\mathbb Q}
\newcommand{\N}{\mathbb N}

\begin{document}

\date{\today}
\title[Natural extensions for Nakada's $\alpha$-expansions: descending from $1$ to $g^2$]{Natural extensions for Nakada's $\alpha$-expansions: \\ descending from $1$ to $g^2$}
\author{Jaap de Jonge}
\address{Jaap de Jonge, Delft University of Technology, department of Electrical Engineering, Mathematics and Computer Science, Mekelweg 4, 2628 CD  Delft and University of Amsterdam, Korteweg - de Vries Institute for Mathematics, Science park 105 - 107, 1098 XG Amsterdam, The Netherlands}
\email{j.dejonge@felisenum.nl, c.j.dejonge@uva.nl, c.j.dejonge@tudelft.nl}
\author{Cor Kraaikamp}
\address{Cor Kraaikamp, Delft University of Technology, department of Electrical Engineering, Mathematics and Computer Science, Mekelweg 4, 2628 CD  Delft, The Netherlands}
\email{C.Kraaikamp@tudelft.nl}
\keywords{Continued fractions, Metric theory}
\subjclass[2000]{Primary 11J70; Secondary 11K50}
\date{\today}

\begin{abstract}
By means of singularisations and insertions in Nakada's $\alpha$-expansions, which involves the removal of partial quotients $1$ while introducing partial quotients with a minus sign, the natural extension of Nakada's continued fraction map $T_\alpha$ is given for $(\sqrt{10}-2)/3\leq\alpha<1$. From our construction it follows that $\Omega_\alpha$, the domain of the natural extension of  $T_\alpha$, is metrically isomorphic to $\Omega_g$ for $\alpha \in [g^2,g)$, where $g$ is the small golden mean. Finally, although $\Omega_\alpha$ proves to be very intricate and unmanageable for $\alpha \in [g^2, (\sqrt{10}-2)/3)$, the $\alpha$-Legendre constant $L(\alpha)$ on this interval is explicitly given. 
\end{abstract}

\maketitle

\section{Introduction}
\label{Introduction}
In 1981, Hitoshi Nakada introduced in \cite{Na} a family of continued fraction maps $T_{\alpha}: [\alpha -1,\alpha ]\to [\alpha -1,\alpha )$, with $\alpha \in [0,1]$, and obtained and studied their natural extensions for $\alpha \in [1/2,1]$; see \cite{Na}, \cite{DK} and \cite{K}. Nakada defined these maps by 
\begin{equation}
\label{alphaCF} T_{\alpha}(x) := \frac{\varepsilon (x)}{x} -a(x) ,\quad x\neq 0; \quad T_{\alpha}(0):=0,
\end{equation}
where $\varepsilon (x)=\text{sign} (x)$ and $a: (-1,1) \to \N \cup \{\infty\}$ is defined by 
$$
a(x):=\left \lfloor \left |\frac1x \right |+1-\alpha\right \rfloor, \quad x\neq 0; \quad a(0):=\infty.
$$
 In case $\alpha =1$, we have the \emph{regular} (or: \emph{simple})
\emph{continued fraction expansion} (rcf), while the case $\alpha
=\frac{1}{2}$ is the \emph{nearest integer continued fraction
expansion} (nicf). In the case $\alpha=0$, we have the {\it{by-access continued fraction expansion}}; see for instance \cite{Z}.
\medskip\

For $\alpha \in [0,1]$, if $T_{\alpha}^{n-1}(x)\neq 0$ for $n \in \N$, we define
$$
\varepsilon_n=\varepsilon_n(x):=\varepsilon \left(
T_{\alpha}^{n-1}(x)\right) ,\quad \text{and}\quad
a_n=a_n(x):=\left\lfloor \left| \frac{1}{T_{\alpha}^{n-1}(x)}\right|
+1-\alpha \right\rfloor.
$$
Applying~(\ref{alphaCF}), we obtain
\begin{equation}\label{expansion}
x=T_\alpha^0(x)=\cfrac{\varepsilon_1}{a_1+T_{\alpha}(x)}=
\cfrac{\varepsilon_1}{a_1+\cfrac{\displaystyle
\varepsilon_2}{\displaystyle a_2+T_{\alpha}^2(x)}}=\dots=\cfrac{\varepsilon_1}{a_1+\cfrac{\displaystyle
\varepsilon_2}{\displaystyle a_2+\cfrac{\varepsilon_3}{a_3 + \ddots }}}\quad,
\end{equation}
which we will throughout this paper write as $x=[0;\varepsilon_1a_1,\varepsilon_2a_2,\varepsilon_3a_3,\dots ]$; we will call the numbers $a_i,\, i \geq 1$, the {\it{partial quotients}} of the continued fraction expansion of $x$.\\
\\
Let $\Omega_{\alpha} \in \R^2$ be the domain for the natural extension of $T_{\alpha}$. We define ${\mathcal T}_{\alpha} :\Omega_{\alpha}\to \Omega_{\alpha}$ by
$$
{\mathcal T}_{\alpha}(x,y):=\left(T_{\alpha}(x),\frac{1}{a(x)+\varepsilon(x) y}\right) .
$$
Omitting the suffix `$(x)$', we define the sequence of `futures' $(t_n)_{n\geq 0}$ and
the sequence of `pasts' $(v_n)_{n\geq 0}$ of $x$ by $v_0=0$, and
$$
(t_n,v_n):={\mathcal T}^n_{\alpha}(x,0),\quad n\in\N \cup \{ 0\}.
$$
We speak of 'futures' and 'pasts' because of the well-known relations (see for instance \cite{DK})
\begin{equation}\label{t en v}
t_n=[0;\varepsilon_{n+1}a_{n+1},\varepsilon_{n+2}a_{n+2},\dots] \quad {\text{and}} \quad v_n=[0;a_n,\varepsilon_na_{n-1},\dots,\varepsilon_2a_1].
\end{equation}
The sequence $(t_n,v_n)_{n\geq 0}$ plays a central role in the theory of approximation by continued fractions and is strongly connected with the sequence of convergents $p_n/q_n$ of $x$. This sequence is given
by $p_n/q_n=[0;\varepsilon_1a_1,\varepsilon_2a_2,\dots,\varepsilon_na_n]$, $n\geq 1$. In \cite{DK}, for instance,  a lot of well-known properties concerning the sequence $(p_n,q_n)_{n\geq -1}$ can be found, such as $\gcd(p_n,q_n)=1$, $n\geq -1$, and 
\begin{equation}\label{pn en qn}
\begin{aligned}
p_{-1}:=1,\quad p_0:=a_0, \qquad p_n=a_np_{n-1}+\varepsilon_np_{n-2}, \, n \geq 1;\hspace{2cm}\\
q_{-1}:=0,\quad q_0:=1, \qquad \, \, \, \, q_n=a_nq_{n-1}+\varepsilon_nq_{n-2}, \, n \geq 1.\hspace{2cm}
\end{aligned}
\end{equation}
The equation $v_n=q_{n-1}/q_n,\, n\geq -1$, serves as an example for the connection between the sequences $(t_n,v_n)_{n\geq 0}$ and $p_n/q_n$; the {\it{approximation coefficients}} 
$$
\theta_n(x):=q_n^2|x-p_n/q_n|,
$$ 
that are used to indicate the quality of the approximation of $x$ by $p_n/q_n$, is another important example. Again from \cite{DK},                                                                                                                                                                                                                                                                                                                                                                                                                                                                                                                                                                             we have, for instance (again omitting the suffix `$(x)$'):
$$
\theta_{n-1}=\frac{v_n}{1+t_nv_n} \quad {\text{and}} \quad \theta_n=\frac{\varepsilon_{n+1}t_n}{1+t_nv_n},
$$
the first of which we will use in the section on the $\alpha$-Legendre constant at the end of this paper.\\
Without giving $\Omega_\alpha$ yet (the determination of which is the main topic of this paper), we note that -- apart from a set of (Lebesgue) measure $0$ -- the map $\mathcal T_\alpha$ is bijective almost everywhere on $\Omega_\alpha$, with inverse map
$$
{\mathcal T}^{-1}_{\alpha}(x,y):=\left(\frac{\varepsilon(x)}{x+a(x)},\frac{1-a(x)y}{\varepsilon(x)y}\right)=\left(\frac{\varepsilon(x)}{x+a(x)},\frac{\varepsilon(x)}y-\varepsilon(x)a(x)\right);
$$
note how the coordinates giving information about `pasts' and `futures' of the continued fraction expansion traded places with those of $\mathcal T_\alpha$.\smallskip

A very helpful tool in the arithmetic of $t_n$ and $v_n$ is the well-known \emph{M\"obius transformation}: \\
Let $M=\left (\begin{array}{cc}
a & b \\
c & d \end{array} \right)$, with $a,b,c,d \in \Z, ad-bc=\pm 1$. The M\"obius transformation $M:\R \cup \{\infty\} \to \R \cup \{\infty\}$ associated with $M$\footnote{Throughout this paper we will omit the distinction between a matrix and its associated M\"obius transformation.} is defined by
$$
M(x):=\frac{ax+b}{cx+d}, \quad x\neq\frac{-d}c,\infty \quad \text{and} \quad M\left (\frac{-d}c \right ) = \infty; \quad M(\infty)=\frac ac.
$$
So, using (\ref{t en v}) and applying M\"obius transformations, we can rephrase (\ref{expansion}) as
$$x=\left (\begin{array}{cc} 0 & \varepsilon_1\\
1 & a_1 \end{array}\right )(t_1)=\left (\begin{array}{cc} 0 & \varepsilon_1\\
1 & a_1 \end{array}\right )\left (\begin{array}{cc} 0 & \varepsilon_2\\
1 & a_2 \end{array}\right )(t_2)=\left (\begin{array}{cc} 0 & \varepsilon_1\\
1 & a_1 \end{array}\right )\left (\begin{array}{cc} 0 & \varepsilon_2\\
1 & a_2 \end{array}\right )\left (\begin{array}{cc} 0 & \varepsilon_3\\
1 & a_3 \end{array}\right )\dots .
$$
Writing $M_n:=\left (\begin{array}{cc} 0 & \varepsilon_n\\
1 & a_n \end{array}\right )$, we also see that 
$$
t_n=M_{n+1} (t_{n+1}) \quad ({\text{so}} \,\,\, t_{n+1}=M_{n+1}^{-1} (t_n))
$$
and generally,
$$
t_n=M_{n+1} \cdots M_{n+k} (t_{n+k}) \quad ({\text{and}} \,\,\, t_{n+k}=M_{n+k}^{-1} \cdots M_{n+1}^{-1}(t_n)).
$$

Likewise,
$$
v_{n+1}=M_{n+1}^T (v_n) 
$$
and generally
$$
v_{n+k}=M_{n+k}^T \cdots M_{n+1}^T (v_n),
$$
from which we derive the very useful fact that 
\begin{equation}\label{relation v_{n+k} and t_{n+k}}
t_{n+k}=M(t_n)\quad {\text{if and only if}} \quad
v_{n+k}=(M^T)^{-1}(v_n),
\end{equation}
with $M$ a M\"obius transformation.\smallskip

The main topic of this paper is the construction of all domains $\Omega_\alpha$ with $\alpha \in [g^2,1)$, where $g:=(\sqrt{5}-1)/2=0.6180\dots$ (while $G:=(\sqrt{5}+1)/2=1.6180\dots$; $g$ and $G$ are for the well-known {\it{golden ratios}}). We will show how to obtain these by transforming $\Omega_{\alpha}$ into $\Omega_{\alpha'}$, with $\alpha' <\alpha$, starting with $\alpha=1$. We will also show that for $\alpha \leq (\sqrt{10}-2)/3$, $\Omega_\alpha$ becomes rapidly extremely intricate, finally making it senseless if not impossible to give the associated set of pairs $(t,v)$ explicitly when $\alpha$ approaches $g^2$. Our approach is based on \emph{singularisations} and \emph{insertions} in the continued fraction expansion of $\alpha$; see \cite{B} and also \cite{K}, where singularisations are used to define a suitable induced transformation. For $A, B\in\N$, $B \geq 2$, $\xi\in \R$ and $\varepsilon \in \{-1,1\}$, the first of these operations is based on the equation
\begin{equation}\label{singularisation}
A+\frac{\displaystyle 1}{\displaystyle 1+\frac{\displaystyle
\varepsilon}{\displaystyle B+\xi}}=A+1-\frac{\varepsilon}{B+\varepsilon+\xi} \quad{\text{(partial quotient 1 between $A$ and $B$ is {\it{singularised)}}}};
\end{equation}
and the second one on  
\begin{equation}\label{insertion}
A+\frac{\varepsilon}{B+\xi}=A+\varepsilon-\frac {\varepsilon}{1+\frac{\displaystyle
1}{\displaystyle B-1+\xi}}\quad{\text{(partial quotient 1 is {\it{inserted}} between $A$ and $B$)}}. 
\end{equation}
In terms of M\"obius transformations, (\ref{singularisation}) is similar to
$$
\left (\begin{array}{cc} 0 & \varepsilon_n\\
1 & a_n \end{array}\right )\left (\begin{array}{cc} 0 & 1\\
1 & 1 \end{array}\right )\left (\begin{array}{cc} 0 & \varepsilon_{n+2}\\
1 & a_{n+2} \end{array}\right )=\left (\begin{array}{cc} 0 & \varepsilon_n\\
1 & a_n+1 \end{array}\right )\left (\begin{array}{cc} 0 & -\varepsilon_{n+2}\\
1 & a_{n+2}+\varepsilon_{n+2} \end{array}\right )
$$
and (\ref{insertion}) to
$$
\left (\begin{array}{cc} 0 & \varepsilon_n\\
1 & a_n \end{array}\right )\left (\begin{array}{cc} 0 & \varepsilon_{n+1}\\
1 & a_{n+1} \end{array}\right )=\left (\begin{array}{cc} 0 & \varepsilon_n\\
1 & a_n +\varepsilon_{n+1} \end {array}\right )\left (\begin{array}{cc} 0 & -\varepsilon_{n+1}\\
1 & 1 \end{array}\right )\left (\begin{array}{cc} 0 & \varepsilon_{n+1}\\
1 & a_{n+1} - \varepsilon_{n+1} \end{array}\right ).
$$
Singularisations {\it{shorten}} the (regular) continued fraction expansion of numbers $x \in \R$. Considering (\ref{expansion}), this implies a {\it{loss}} of convergents (as will be illustrated in the next section). Insertions obviously lengthen a continued fraction expansion, yielding so-called {\it{mediants}}, that are usually of lesser approximative quality compared to convergents. It be remarked, however, that inserting partial quotients $1$ in a given continued fraction expansion involves the possibility of creating points $(t,v)$ lying in a part of $\Omega_1$ that had been removed in a previous stage of constructing. As long as we can compensate each insertion with a singularisation, we can avoid this. Applying this {\it{compensated insertion}} is actually the base of our construction of $\Omega_{\alpha}$, from $\alpha=g$ downwards, that we will gradually unfold in the next sections. \\
It is easy to see that compensation is only possible as long as an insertion involves the introduction of an extra 1, that is, if $B-1$ in (\ref{insertion}) equals $1$, implying $B=2$. In this case we deal with
$$
A+\frac{\displaystyle \varepsilon_1}{\displaystyle 2+\frac{\displaystyle
\varepsilon_2}{\displaystyle C+\xi}}=A+\varepsilon_1-\frac{\displaystyle \varepsilon_1}{\displaystyle 1+\frac{\displaystyle
1}{\displaystyle 1+\frac{\displaystyle \varepsilon_2}{C+\xi}}}=A+\varepsilon_1-\frac{\displaystyle \varepsilon_1}{\displaystyle 2-\frac{\displaystyle
\varepsilon_2}{\displaystyle C+\varepsilon_2+\xi}}.
$$
Of course, similar equations could be given in terms of M\"obius transformations, but at this point we want to stress that either notation has its pros and cons as to readibility and illuminating power. Depending on the context, for our computations we will use either continued fractions or M\"obius transformations.\smallskip

In the next sections we will step-by-step transform the domain $\Omega_1:=[0,1)\times[0,1]$ for the natural extension of the regular continued fraction to the domain $\Omega_\alpha$ of the natural extension of any $\alpha$-expansion for $\alpha \in [g^2,1)$. Note that $t_n=1$ only occurs in the case $x=1$ and $n=0$, associated with the single point $(1,0)$, that is sent to $(0,1)$ under $T_1$. So, in view of (\ref{alphaCF}), one should actually define $\Omega_1:=[0,1)\times[0,1] \cup \{(1,0)\}$. We will suppress subtleties like these in similar cases throughout this paper.\smallskip

We will often use some well-known equations involving $g$ and $G$, such as $g^2=1-g$ and $\tfrac1g=G=g+1$. Although we let $\alpha$ decrease continuously from $1$ downwards, we speak of 'step-by-step', as the generic form of  $\Omega_\alpha$ varies over several intervals, according to feasability for singularisation (possible in Section \ref{The case alpha in (g,1]} only) or compensated insertion. In Section \ref{intervals} we will explain how to determine these intervals. We remark that the results for $\alpha\geq\sqrt{2}-1$ are not new (see for instance \cite{K} or \cite{Na}). However, our approach is best introduced by starting at $\alpha=1$ and has some illuminative qualities. For one thing, we will show how our construction of $\Omega_\alpha$ calls for the term {\it{quilting}}, as introduced in \cite{KSS}. More importantly, some interesting results on $\alpha$-expansions come with the determination of $\Omega_\alpha$. Since the construction of $\Omega_\alpha$ will prove to be quite intricate, we will present these results mostly separately at the end of this paper. In Section \ref{ergodic systems} we will go into ergodic properties of $\mathcal T_\alpha$. Representing the collection of subsets of $\Omega_\alpha$ by $\mathcal B$ and defining $\mu_\alpha$ as the probability measure with density 
$$
\frac1{N_\alpha}\,\cdot\, \frac1{(1+tv)^2}
$$ 
on $(\Omega_\alpha,\mathcal B)$ (with $N_\alpha$ a normalising constant), we will also show that $(\Omega_\alpha, \mu_\alpha, \mathcal T_\alpha)$ is an ergodic dynamical system.\smallskip

Finally, in Section \ref{legendre} we will extend a result of Rie Natsui (\cite{Nat}), based on the well-known theorem of Legendre in the theory of regular continued fractions; see for instance \cite{BJ}:
\begin{Theorem}\label{Legendre Theorem}
Let $x \in \R\setminus \Q$ and $p/q \in \Q$, with $\gcd(p,q)=1$ and $q>0$, such that $q^2|x-p/q|<1/2$. Then there exists a non-negative integer $n$ such that $p/q=p_n/q_n$, in words: $p/q$ is a regular convergent  of $x$. This constant $1/2$ is best possible in the sense  that for every $\varepsilon >0 $ an irrational $x$ and a rational number $p/q$, with $\gcd(p,q)=1$ and $q>0$, exist such that $p/q$ is not a regular convergent and $q^2|x-p/q|<1/2+\varepsilon$.
\end{Theorem}
In \cite{Nat}, Natsui proves the existence of the {\it{$\alpha$-Legendre constant}} for $\alpha$-continued fractions, defined by
$$
L(\alpha):=\sup\{c>0: q^2\left |x-\tfrac pq \right |<c, \gcd(p,q)=1 \Rightarrow \tfrac pq = \tfrac{p_n}{q_n}, n \geq 0\}
$$
for $0<\alpha \leq 1$; here $p_n/q_n$ is an $\alpha$-convergent of $x$. She expands the result given by Ito (\cite{I2}; see also \cite{K}) for values $1/2\leq\alpha\leq1$ to giving explicit values of $L(\alpha)$ for $\sqrt{2}-1 \leq \alpha <1/2$. In Section \ref{legendre} we will augment this result with the interval $[g^2,\sqrt{2}-1)$ and the numbers $\alpha=1/n$, for $n \in \N$, $n \geq 3$. \smallskip

We will first concentrate on the construction of $\Omega_\alpha$ for various values of $\alpha$.

\section{Determining the intervals according to type of singularisation}
\label{intervals}
One way to look at the domains $\Omega_\alpha$, $\alpha < 1$, is that in $\Omega_\alpha$ no points $(t,v)$ exist for which $\alpha \leq t \leq 1$. In combination with the definition of $(t_n,v_n)$ as $\mathcal T^n_{\alpha}(x,0),\, n \in\N \cup \{ 0\}$, for any $x \in [\alpha-1,\alpha]$, this forms the basis of our approach for constructing the $\Omega_\alpha$ with $\alpha \in [g^2,1)$. We start with $\Omega_1$ and fix $\alpha$ such that each time $\alpha \leq t_n \leq 1$, we replace the point $(t_n,v_n)$ by a point $(t_n^*,v_n^*)$ belonging to the continued fraction expansion of $x$ after singularising. When fixing $\alpha$, we use the continued fraction map $T_1$ of the associated ergodic system. Then we follow the same procedure, this time starting with $\Omega_\alpha$, fixing an $\alpha' < \alpha$ and applying $T_\alpha$ in order to construct $\Omega_{\alpha'}$. From Section \ref{The case alpha in (frac12,g]} on, we will show how to make use of compensated insertion when $\alpha \leq g$.\smallskip

We remark that under $T_\alpha$ the continued fraction expansion of numbers can vary considerably, according to the value of $\alpha$. Denoting with $x_\alpha$ the $\alpha$-expansion of $x \in (0,1)$, we have, for instance, $g_1=[0;\overline{1}]$ and $g_g=[0;2,\overline{-3}]$, the bar indicating infinite repetition. Constructing $\Omega_{\alpha'}$ from $\Omega_\alpha$, the way points are replaced depends on the way $x_{\alpha}$ changes under singularisation or (compensated) insertion. As a start, in Section \ref{The case alpha in (g,1]} we will remove all points $(t,v)$ for which $t=[0;1,\dots]$, to which end mere singularisation suffices. In Sections \ref{The case alpha in (frac12,g]} and \ref{The case alpha in (sqrt{2}-1,tfrac12)} we will remove all points $(t,v)$ for which $t=[0;2,\dots]$ by means of a single compensated insertion. In Section \ref{The case alpha in (tfrac{sqrt{10}-2}3,sqrt{2}-1]} we will show that removing points $(t,v)$ for which $t=[0;3,\dots]$ (the ultimate one being $t=g^2$) is possible only in special cases, applying compensated insertion twice. 

Not only the continued fraction expansion of a number $x$ depends on which continued fraction map we use, the expansion of the associated numbers $t_n$ does also. In order to  determine the interval $[\alpha', \alpha]$, $\alpha' <\alpha \leq 1$, such that $t_n \in (\alpha', \alpha)$ implies $t_n=[0;1,\dots]$, we observe that on $(0,1]$ both $\tfrac1\alpha$ and $1-\alpha$ are decreasing functions of $\alpha$. So, for $\alpha \in (\alpha',1]$,
\begin{equation}\label{boundary between 1 and 2}
1 = \left \lfloor \frac11+1-1\right \rfloor \leq  \left \lfloor \frac1\alpha +1-\alpha \right \rfloor \leq a_1(\alpha'_\alpha)=\left \lfloor \frac1{\alpha'}+1-\alpha \right \rfloor \leq \left \lfloor \frac1{\alpha'}+1-\alpha' \right \rfloor, 
\end{equation}
from which we derive $a_1(\alpha'_\alpha)=1$ if and only if $1/\alpha'+1-\alpha' <2$, which is the case if and only if $g<\alpha'<1$. This is why the next section is about the case $\alpha \in (g,1]$. There (and more so in the sections following it) we will show that although replacing points by other points is in fact only part of the construction of the $\Omega_\alpha$, the heart of the construction is still the procedure sketched above, which is actually a way to determine {\it{$\alpha$-fundamental intervals}} $\Delta_{n,\alpha}=\Delta_\alpha(i_1,i_2,\dots,i_n)$ of rank $n$, that we define as follows:
$$
\Delta_\alpha(i_1,i_2,\dots,i_n):=\{\alpha \in [0,1]: a_1(\alpha_\alpha)=i_1, a_2(\alpha_\alpha)=i_2, \dots,a_n(\alpha_\alpha)=i_n\},
$$
where $ i_j \in \Z \setminus \{0\}$ for each $1\leq j\leq n$. In the previous paragraph, for instance, we found $\Delta_\alpha(1)=\{\alpha \in [0,1]: a_1(\alpha_\alpha)=1\}=(g,1]$.

\section{The case $\alpha \in (g,1]$}
\label{The case alpha in (g,1]}
We start with the natural extension domain for the regular continued fraction, the square $\Omega_1$,  and the classic {\it{Gauss map}} $T_1$, given by Nakada in \cite{Na}. Let $x \in [0,1]$ and $n \geq 0$ be the smallest integer for which $t_n \in [\alpha,1)$, i.e., for which $(t_n,v_n) \in R_\alpha:=[\alpha,1)\times[0,1]$. We then have that the rcf of $x$ satisfies
$$
x=[0;a_1,\dots,a_n,1,a_{n+2},a_{n+3},\dots];
$$
 observe that in $\Omega_1$ we have $\varepsilon_k=1$, $k\geq1$. Singularising $a_{n+1}=1$, we obtain 
$$
x=[0;a_1,\dots,(a_n+1),-(a_{n+2}+1),a_{n+3},\dots].
$$
Let the sequence $(t_n^*,v_n^*)_{n\geq 0}$ denote the futures and pasts of the singularised continued fraction expansion of $x$. Then 
$$
(t_n^*,v_n^*)=([0;-(a_{n+2}+1),a_{n+3},\dots],[0;a_n+1,a_{n-1},\dots,a_1]),
$$

while 
$$
(t_n,v_n)=([0;1,a_{n+2},a_{n+3},\dots],[0;a_n,a_{n-1},\dots,a_1]).
$$

Note that $v_n=1/(a_n+v_{n-1})$, which is equivalent to $a_n+v_{n-1}=1/v_n$, and so
$$
v_n^*=\frac1{a_n+1+v_{n-1}}=\frac1{1+\frac1{v_n}}=\frac{v_n}{1+v_n}.
$$
Applying (\ref{insertion}) with $A=0$, $\varepsilon=-1$, $B=a_n+1$ and $\xi=t_{n+2}$, we get
$$
t_n^*=\frac{-1}{a_n+1+t_{n+2}}=-1+\frac 1{1+\frac{\displaystyle
1}{\displaystyle a_n+t_{n+2}}}=t_n-1.
$$
Since $(t_n,v_n) \in R_\alpha$, we find that $(t^*,v^*) \in A_\alpha:=[\alpha-1,0)\times[0,1/2]$.\smallskip

Note that the map $\mathcal M: R_{\alpha} \to A_{\alpha}$, defined by $\mathcal M(t,v):=(t-1,v/(1+v))$ is a bijection. In Section \ref{ergodic systems} we will use this to obtain ergodicity of the dynamical system $(\Omega_\alpha,\mathcal B, \mu_\alpha,\mathcal T_\alpha)$, $g \leq \alpha \leq 1$. 

It is obvious that $t_{n+1}^*=t_{n+2}$, and since 
$$
v^*_{n+1}=\frac1{a_{n+2}+1-v_n^*}=\frac1{a_{n+2}+1-\tfrac{v_n}{1+v_n}}=\frac1{a_{n+2}+\tfrac1{1+v_n}}=\frac1{a_{n+2}+v_{n+1}}=v_{n+2},
$$
we conclude that 
\begin{equation}\label{(t,v)^* en (t,v)1}
(t_{n+1}^*, v_{n+1}^*)=(t_{n+2},v_{n+2}). 
\end{equation}

Having singularised the first partial quotient $1$ in the continued fraction expansion of $x$ and so obtained an alternative expansion, we repeat the procedure ad infinitum. When doing so, we actually not only remove $R_a$ from $\Omega_1$, but $\mathcal T_1(R_\alpha):=(0,(1-\alpha)/\alpha]\times [1/2,1]$ as well, since  $(t_{n+1},v_{n+1})=((1-t_n)/t_n,1/(1+v_n))$. Equation (\ref{(t,v)^* en (t,v)1}) implies $\mathcal T_1(A_\alpha)=\mathcal T_1^2(R_\alpha)$, so the loss of points is confined to $R_\alpha$ and $\mathcal T_1(R_\alpha)$. Removing from $\Omega_1$ the strips $R_\alpha$ and $\mathcal T_1(R_\alpha)$ and adding $A_\alpha$ yields $\Omega_\alpha$; see Figures \ref{fig: omega1} and \ref{fig: omega0.8}. Here we see a first glimpse of the similarity with the process of quilting as described in \cite{KSS}. Expanding this metaphor, our construction in the next sections may be described as a form of recursive quilting.\smallskip

We remark that $\Omega_\alpha$ can be constructed from $\Omega_1$ immediately, using $\mathcal T_1$, or from some $\Omega_{\alpha'}$, $\alpha<\alpha'<1$, using $\mathcal T_{\alpha'}$. Indeed, in these cases we have something similar to (\ref{boundary between 1 and 2}):
$$
\left \lfloor \frac1\alpha +1-\alpha' \right \rfloor < \left \lfloor \frac1g+1-g\right \rfloor = 2,
$$
implying that changing ``$t_n \in [\alpha,1)$" in the beginning of this section into ``$t_n \in [\alpha,\alpha')$" would leave the construction of $\Omega_\alpha$ unaffected, save for the initial domain. \smallskip

Rewriting $x=[0;a_1,\dots,a_n,1,a_{n+2},\dots]$ as $x=[0;a_1,\dots,(a_n+1),-(a_{n+2}+1),\dots]$ does {\it{not}} leave the sequence of convergents unaffected: since $p_k^*=p_k$ and $q_k^*=q_k$, $k<n$, we find, applying (\ref{pn en qn}), 
\begin{align*}
& p_n^*=(a_n+1)p_{n-1}^*+p_{n-2}^*=a_np_{n-1}^*+p_{n-1}^*+p_{n-2}^*=a_np_{n-1}+p_{n-2}+p_{n-1}=p_n+p_{n-1}=p_{n+1};\\
& p_{n+1}^*=2p_n^*-p_{n-1}^*=a_{n+2}p_{n+1}+p_{n+1}-p_{n-1}=a_{n+2}p_{n+1}+p_n=p_{n+2}.
\end{align*}
Similarly, we find $q_n^*=q_{n+1}$ and $q_{n+1}^*=q_{n+2}$.  We conclude that in the transformation from $\Omega_1$ to $\Omega_\alpha$, $\alpha \in (g,1)$, a convergent $p_n/q_n$ is lost every time that $(t_n,v_n) \in R_\alpha$. \smallskip

Note that $\lim_{\alpha \downarrow g} (1-\alpha)/\alpha = (1-g)/g=g=\lim_{\alpha \downarrow g} \alpha$, which is mirrored in the fact that for $\alpha=g$ the left boundary of $R_\alpha$ coincides with the right boundary of $\mathcal T_1(R_\alpha)$. We now have obtained the following result:
\begin{Theorem}
Let $\alpha \in (g,1)$. Then $\Omega_\alpha=[\alpha-1,0)\times[0,\tfrac12]\cup[0,\tfrac{1-\alpha}{\alpha}]\times[0,\tfrac12)\cup(\tfrac{1-\alpha}{\alpha},\alpha)\times[0,1].$
\end{Theorem}
This theorem was already obtained by Nakada in \cite{Na}; see also \cite{K}. Figure \ref{fig: omega0.8} shows an example of the generic form of $\Omega_\alpha$, $\alpha \in (g,1)$. In Figure  \ref{fig: omegag} we have written $-g^2$ for the equivalent $g-1$.

\begin{figure}[!htb]
\minipage{0.32\textwidth}
$$
\begin{tikzpicture}[scale =4] 
  \draw[draw=white,fill=gray!20!white] 
 plot[smooth,samples=100,domain=.8:1] (\x,0) --
 plot[smooth,samples=100,domain=1:.8] (\x,1);
 \draw[draw=white,fill=gray!20!white] 
 plot[smooth,samples=100,domain=0:.25] (\x,.5) --
 plot[smooth,samples=100,domain=.25:0] (\x,1);
 \node at (.125,.75){$\mathcal T_\alpha(R_\alpha)$};
  \draw [dashed] (.25,0) -- (.25,1);
   \draw [dashed] (0,.5) -- (.25,.5);
 \draw (0,0) -- (1,0);
 \draw (0,0) -- (0,1);
  \draw  [dashed] (1,0) -- (1,1);
 \draw (0,1) -- (1,1);
 \draw (.8,0) -- (.8,1);
 \node at (.9,.5){$R_\alpha$};
 \node at (0,-.07){$0$};
  \node at (.8,-.07){$\alpha$};
 \node at (1,-.07){$1$};
  \node at (-.04,1){$1$};
  \node at (-.04,.5) {$\tfrac12$};
       \node at (.25,-.07){$\tfrac{1-\alpha}{\alpha}$};
 
 \end{tikzpicture}
 $$
 \caption[$\Omega_1, \, \alpha=0.8.$] {\label{fig: omega1}
$\Omega_1, \, \alpha=0.8.$}
\endminipage\hfill
\minipage{0.32\textwidth}
$$
\begin{tikzpicture}[scale =4] 
  \draw (-.2,0) -- (.8,0);
 \draw (-.2,0) -- (-.2,.5);
  \draw (-.2,.5) -- (0,.5);
  \draw [dashed] (0,.5) -- (.25,.5);
   \draw [dashed] (.8,0) -- (.8,1);
 \draw (0,0) -- (0,.5);
 \draw [dashed] (.25,0) -- (.25,1);
 \draw (.25,1) -- (.8,1); 
   \node at (-.2,-.07){$\alpha-1$};
  \node at (0,-.07){$0$};
 \node at (.8,-.07){$\alpha$};
   \node at (-.24,.5){$\tfrac12$};
   \node at (-.1,.25){$A_\alpha$};
   \node at (.21,1) {$1$};
    \node at (.25,-.07){$\tfrac{1-\alpha}{\alpha}$};
      \node at (0,1.05) {};
 \end{tikzpicture}
 $$
 \caption[$\Omega_\alpha$ with $\alpha=0.8$.] {\label{fig: omega0.8}
$\Omega_\alpha$, $\alpha=0.8$.}
\endminipage\hfill
\minipage{0.32\textwidth}
$$
\begin{tikzpicture}[scale =4] 
 \draw (-.382,0) -- (.618,0);
 \draw (-.382,0) -- (-.382,.5);
  \draw (-.382,.5) -- (0,.5);
  \draw [dashed] (0,.5) -- (.618,.5);
 \draw  [dashed]  (.618,0) -- (.618,.5);
 \draw (0,0) -- (0,.5);
 \node at (-.382,1) {};
  \node at (-.382,-.07){$-g^2$};
  \node at (0,-.07){$0$};
 \node at (.618,-.07){$g$};
   \node at (-.42,.5){$\tfrac12$};
      \end{tikzpicture}
 $$
 \caption[$\Omega_g$] {\label{fig: omegag}
$\Omega_g.$}
\endminipage
\end{figure}

\section{The case $\alpha\in (\frac12,g]$}
\label{The case alpha in (frac12,g]}
An important implication of our construction of $\Omega_{g}$ in the previous section is that in continued fraction expansions associated to $\Omega_\alpha$, $\alpha \leq g$, the partial quotient $1$ is non-existent. In the current and the following section we will similarly remove all partial quotients $2$ (with plus sign). We determine the interval $(\alpha', \alpha]$, $\alpha' <\alpha \leq g$, such that $t_n \in (\alpha', \alpha]$ implies $t_n=[0;2,\dots]$, in a way similar to the  one at the end of Section (\ref{intervals}): if $x \in [0,g)$ and $\alpha \in (x,g]$, then
$$
2 = \left \lfloor \frac1g+1-g\right \rfloor \leq  \left \lfloor \frac1\alpha +1-\alpha \right \rfloor \leq a_1(\alpha'_\alpha)=\left \lfloor \frac1{\alpha'} +1-\alpha \right \rfloor \leq \left \lfloor \frac1{\alpha'}+1-\alpha' \right \rfloor, 
$$
from which we derive $a_1(\alpha'_\alpha)=2$ if and only if $1/\alpha'+1-\alpha' <3$, which (given that $\alpha' \in [0,g)$) is the case if and only if $\sqrt{2}-1<\alpha'<g$; indeed $\sqrt{2}-1_{\sqrt{2}-1}=[0;3,\overline{-2,-4}]$. We conclude that $\Delta_\alpha(2)=(\sqrt{2}-1,g]$, and therefore we will now investigate the case $\alpha \in (\sqrt{2}-1,g]$.\smallskip

In the current section we will confine ourselves to the case $\alpha \in (1/2,g]$. To make a cut at $\alpha=1/2$ is because obviously $t \in (1/2,g]$ if and only if $t=[0;2,-a,\dots]$, with $a \in \N_{\geq 2}$, while $t \in (\sqrt{2}-1,1/2]$ if and only if $t=[0;2,b,\dots]$, with $b \in \N_{\geq 2}$. We will see how the distinction shows in different ways of transforming $\Omega_\alpha$.

Let $x \in [g-1,g)$ and let $n \geq 0$ be the smallest integer for which $t_n \in [\alpha,g]$, i.e., for which $(t_n,v_n) \in R_\alpha:=[\alpha,g]\times[0,1/2)$. We already know that in this case $x=[0;\varepsilon_1a_1,\dots,\varepsilon_na_n,2,-a_{n+2},\dots]$. Since the continued fraction map of the related dynamical system is $T_g$, for $\alpha\in (1/2,g]$ we have
$$
a_2(\alpha)=\left \lfloor \left| \tfrac1{T_g(\alpha)}\right| +1-g \right \rfloor = \left \lfloor \left|\tfrac1{\tfrac1\alpha-2}\right|+1-g\right \rfloor \geq \left \lfloor \left|\tfrac1{\tfrac1g-2}\right|+1-g\right \rfloor =  G+1+1-g=3,
$$
so in this case 
$$
x=[0;\varepsilon_1a_1,\dots,\varepsilon_na_n,2,-a_{n+2},\dots], \quad a_{n+2} \geq 3.
$$
To remove the partial quotient $2$, we first insert $-1$, so as to write 
$$
x=[0;\varepsilon_1a_1,\dots,\varepsilon_n(a_n+1),-1,1,-a_{n+2},\dots].
$$
Now we singularise $1$ in this expansion and get 
$$
x=[0;\varepsilon_1a_1,\dots,\varepsilon_n(a_n+1),-2,(a_{n+2}-1),\dots].
$$

Having removed the first partial quotient $2$ in the continued fraction expansion of $x$ and so obtained the alternative expansion, we repeat the procedure ad infinitum. Note that the next $+2$ could be $a_{n+2}-1$ in the rewritten expansion, since $a_{n+2} \geq 3$. In the current case we have $R_\alpha=[\alpha,g]\times[0,1/2)$ and, considering that $(t_{n+1},v_{n+1})=((1-2t_n)/t_n,1/(2+v_n))$, $\mathcal T_g(R_\alpha)=[-g^2,(1-2\alpha)/\alpha]\times(2/5,1/2]$, knowing that $(1-2g)/g=-g^2$; see Figure \ref{fig: omegag2}, where we have taken $\alpha=0.52$ for illustrative purposes. Having determined what to remove (in the sense of the previous section), we use similar calculations as before to find what to add. Using the same notation as in the previous case, we find 
\begin{align*}
(t_n^*,v_n^*)&=(t_n-1,\tfrac{v_n}{1+v_n}) \quad {\text{(as in the previous case),}}\\
(t_{n+1}^*, v_{n+1}^*)&=(-\tfrac{t_{n+1}}{1+t_{n+1}},1-v_{n+1})=(\tfrac{1-2t_n}{t_n-1},\tfrac{v_n+1}{v_n+2}) \,\, {\text{and}}\\
(t_{n+2}^*, v_{n+2}^*)&=(t_{n+2}, v_{n+2}).
\end{align*}
So there are two regions to be added: $A_\alpha=[\alpha-1,-g^2]\times[0,1/3)$ and $\mathcal T_g(A_\alpha)=[(1-2\alpha)/(\alpha-1),(1-2g)/(g-1)]\times[1/2,3/5)$. With regard to the last region, we note that $(1-2g)/(g-1)=g$. Because of the latter, we have not yet established to construct $\Omega_\alpha$, since $\mathcal T_g(A_\alpha)$ contains points $(t,v)$ that do not exist in $\Omega_\alpha$, i.e. those for which $t \in (\alpha,g]$; see the greyed region in Figure \ref{fig: omega.521}. We will call the newly constructed region $\Omega_{\alpha,1}^+$, defined by
$$
\Omega_{\alpha,1}^+:=\Omega_g\cup A_\alpha\cup \mathcal T_g(A_\alpha) \setminus (R_\alpha \cup \mathcal T_g(R_\alpha))
$$
Putting $R_{\alpha,1}:=\mathcal T_g(A_\alpha) \cap [\alpha,g]\times[0,3/5)$ as the part yet to be removed, we define
$$
\Omega_{\alpha,1}:=\Omega_{\alpha,1}^+ \setminus R_{\alpha,1}.
$$
\begin{figure}[!htb]
\minipage{0.48\textwidth}
$$
\begin{tikzpicture}[scale =6] 
  \draw[draw=white,fill=gray!20!white] 
 plot[smooth,samples=100,domain=.52:.618] (\x,0) --
 plot[smooth,samples=100,domain=.618:.52] (\x,.5);
 \draw[draw=white,fill=gray!20!white] 
 plot[smooth,samples=100,domain=-.382:-.077] (\x,.4) --
 plot[smooth,samples=100,domain=-.077:-.382] (\x,.5);
 \node at (-.23,.45){$\mathcal T_g(R_\alpha)$};
   \draw [dashed] (0,.5) -- (.25,.5);
    \draw [dashed] (0,.5) -- (.618,.5);
    \draw [dashed] (-.382,.4) -- (-.077,.4);
     \draw (-.382,.4) -- (-.382,.5);
    \draw [dashed] (-.077,0) -- (-.077,.4);
      \draw (-.077,.4) -- (-.077,.5);
   \draw (-.382,.5) -- (0,.5);
    \draw  [dashed]  (-.382,0) -- (-.382,.4);
       \draw (.52,0) -- (.52,.5);
 \draw (-.382,0) -- (.618,0);
 \draw (0,0) -- (0,.5);
  \draw  (.618,0) -- (.618,.5);
  \node at (-.382,-.05){$-g^2$};
 \node at (.57,.25){$R_\alpha$};
 \node at (.618,-.05){$g$};
 \node at (-.12,-.05){$\tfrac{1-2\alpha}{\alpha}$};
 \node at (.01,-.05) {0};
 \node at (0,.63) {};
 \node at (-.43,.4) {$\tfrac25$};
  \node at (-.43,.5) {$\tfrac12$};
   \node at (.52,-.05) {$\alpha$};

 \end{tikzpicture}
 $$
 \caption[$\Omega_g$ with $R_\alpha$ and $\mathcal T_g(R_\alpha)$, $\alpha=0.52$.] {\label{fig: omegag2}
$\Omega_g$ with $R_\alpha$ and $\mathcal T_g(R_\alpha)$, $\alpha=0.52$.}
\endminipage\hfill
\minipage{0.48\textwidth}
$$
\begin{tikzpicture}[scale =6] 
\draw[draw=white,fill=gray!20!white] 
 plot[smooth,samples=100,domain=.52:.618] (\x,.5) --
 plot[smooth,samples=100,domain=.618:.52] (\x,.6);
 \draw (-.48,0) -- (.52,0);
 \draw [dashed] (.52,0) -- (.618,0);
 \draw (-.382,0) -- (-.382,.4);
  \draw (.083,.5) -- (.618,.5);
 \draw [dashed] (.618,0) -- (.618,.5);
 \draw (-.48,0) -- (-.48,.333);
  \draw [dashed] (-.48,.333) -- (-.382,.333);
  \draw [dashed] (.083,0) -- (.083,.6);
   \draw (.083,.5) -- (.083,.6);
  \draw [dashed] (.083,.6) -- (.618,.6);
   \draw (.52,0) -- (.52,.6);
  \draw (0,.4) -- (0,.5);
  \draw [dashed] (0,.5) -- (.083,.5); 
    \draw [dashed] (-.077,0) -- (-.077,.5);
   \draw (.618,.5) -- (.618,.6);
  \draw (-.382,.4) -- (-.077,.4);
  \draw (-.077,.5) -- (0,.5);
 \draw (0,0) -- (0,.4);
  \node at (-.36,-.05){$-g^2$};
  \node at (-.52,-.05) {$\alpha-1$};
 \node at (.618,-.05){$g$};
   \node at (-.41,.4){$\tfrac25$};
     \node at (-.35,.333){$\tfrac13$};
        \node at (-.11,.5){$\tfrac12$};
     \node at (.06,.6){$\tfrac35$};
     \node at (.15,-.05) {$\tfrac{1-2\alpha}{\alpha-1}$};
       \node at (.52,-.05) {$\alpha$};
        \node at (.3,.55){$\mathcal T_g(A_\alpha)$};
         \node at (-.43,.2){$A_\alpha$};
        \node at (.565,.55){$R_{\alpha,1}$};
   \node at (-.13,-.05){$\tfrac{1-2\alpha}{\alpha}$};
   \node at (.01,-.05) {0};
     \node at (0,-.14) {};
  \end{tikzpicture}
 $$
 \caption[$\Omega_{\alpha,1}^+$, $\alpha=0.52$.] {\label{fig: omega.521}
$\Omega_{\alpha,1}^+$, $\alpha=0.52$.}
\endminipage
\end{figure}
We will follow the same procedure of removal and addition regarding $R_{\alpha,1}$ as in the case of $R_\alpha$. Note that this renders the same $t$-coordinates as in the creation of $\Omega_{\alpha,1}$. To efficiently calculate the $v$-coordinates, it is convenient to use M\"obius transformations. We saw that in the current construction of $\Omega_\alpha$ the equation $v_n^*=v_n/(v_n+1)$ holds, that we now write as $v_n^*=\left (\begin{array}{cc}1 & 0 \\
1 & 1 \end{array} \right)(v_n)$; similarly we write $v_{n+1}^*=\left (\begin{array}{cc}1 & 1 \\
1 & 2 \end{array} \right)(v_n)$. It follows that the $v$-coordinate of points in $A_{\alpha,1}$ are given by 
\begin{equation}\label{v_n,1*}
v_{n,1}^*=\left (\begin{array}{cc}1 & 0 \\
1 & 1 \end{array} \right)\left (\begin{array}{cc}1 & 1 \\
1 & 2 \end{array} \right)(v_n), 
\end{equation}
with $v_n$ the $v$-coordinate of points in $R_\alpha$; for the $v$-coordinate of points in $\mathcal T_g (A_{\alpha,1})$ we have 
\begin{equation}\label{v_n+1,1*}
v_{n+1,1}^*=\left (\begin{array}{cc}1 & 1 \\
1 & 2 \end{array} \right)\left (\begin{array}{cc}1 & 1 \\
1 & 2 \end{array} \right)(v_n).
\end{equation}
For the construction of $\Omega_{\alpha,2}$ we also need to calculate $\mathcal T_g(R_{\alpha,1})$. Since
$v_{n+1}=\tfrac1{v_n+2}=\left (\begin{array}{cc}0 & 1\\
1 & 2 \end{array} \right)(v_n)$, we find that the $v$-coordinates of $\mathcal T_g(R_{\alpha,1})$ are given by 
\begin{equation}\label{v_n+1,1}
v_{n+1,1}=\left (\begin{array}{cc} 0 & 1 \\
1 & 2 \end{array} \right)\left (\begin{array}{cc}1 & 1 \\
1 & 2 \end{array} \right)(v_n), 
\end{equation}
$v_n$ still being the $v$-coordinate of points in $R_\alpha$. Since $v \in [0,1/2)$ for points in $R_\alpha$, calculating the boundaries of $\Omega_{\alpha,2}$ is now straightforward. We have written the relevant values of $t$ and $v$ in Figure \ref{fig: omega.522} and write 
$$
\Omega_{\alpha,2}^+:=\Omega_{\alpha,1}^+\cup A_{\alpha,1}\cup \mathcal T_g(A_{\alpha,1}) \setminus (R_{\alpha,1} \cup \mathcal T_g(R_{\alpha,1}))
$$
Defining $R_{\alpha,2}:=\mathcal T_g(A_{\alpha,1}) \cap [\alpha,g]\times[3/5,8/13)$ as the part yet to be removed, we define
$$
\Omega_{\alpha,2}:=\Omega_{\alpha,2}^+ \setminus R_{\alpha,2}.
$$
\begin{figure}[!htb]
$$
\begin{tikzpicture}[scale =11] 
\draw[draw=white,fill=gray!20!white] 
 plot[smooth,samples=100,domain=.52:.618] (\x,.6) --
 plot[smooth,samples=100,domain=.618:.52] (\x,.62);
 \draw (-.48,0) -- (.52,0);
 \draw [dashed] (.52,0) -- (.618,0);
 \draw (-.382,0) -- (-.382,.385);
  \draw [dashed] (.618,0) -- (.618,.6);
 \draw (-.48,0) -- (-.48,.375);
  \draw [dashed] (-.48,.375) -- (-.382,.375);
    \draw [dashed] (.083,0) -- (.083,.5);
  \draw (.083,.5) -- (.083,.62);
  \draw (.083,.6) -- (.618,.6);
   \draw [dashed] (.083,.62) -- (.618,.62);
  \draw [dashed] (.52,0) -- (.52,.6);
   \draw (.52,.6) -- (.52,.62);
     \draw (.618,.6) -- (.618,.62);
  \draw (0,.4) -- (0,.5);
  \draw [dashed] (0,.5) -- (.083,.5); 
    \draw [dashed] (-.077,0) -- (-.077,.5);
     5    \draw (-.382,.385) -- (-.077,.385);
  \draw (-.077,.5) -- (0,.5);
 \draw (0,0) -- (0,.4);
  \node at (-.38,-.02){$-g^2$};
  \node at (-.49,-.02) {$\alpha-1$};
 \node at (.618,-.02){$g$};
   \node at (-.4,.4){$\tfrac5{13}$};
     \node at (-.5,.375){$\tfrac38$};
        \node at (-.09,.5){$\tfrac12$};
     \node at (.065,.59){$\tfrac35$};
        \node at (.065,.63){$\tfrac8{13}$};
     \node at (.09,-.03) {$\tfrac{1-2\alpha}{\alpha-1}$};
       \node at (.52,-.02) {$\alpha$};
        \node at (.57,.612){$R_{\alpha,2}$};
   \node at (-.077,-.03){$\tfrac{1-2\alpha}{\alpha}$};
   \node at (0,-.02) {0};
  \end{tikzpicture}
 $$
 \caption[$\Omega_{\alpha,2}$, $\alpha=0.52$.] {\label{fig: omega.522}
$\Omega_{\alpha,2}$, $\alpha=0.52$.}
\end{figure}
Proceeding this way, we construct the sequence $\Omega_{\alpha,k}$, $k=1,2,\dots$, and bring forth $\Omega_\alpha$ as being $\lim_{k \to \infty} \Omega_{\alpha,k}$. To do so, we need to calculate $\lim_{k \to \infty} v_{n,k}^*$, $\lim_{k \to \infty} v_{n+1,k}^*$ and $\lim_{k \to \infty} v_{n+1,k}$, using formulas associated with (\ref{v_n,1*}) through  (\ref{v_n+1,1}):
$$
\left \{
\begin{aligned}
v_{n,k}^*&=\left (\begin{array}{cc}1 & 0 \\
1 & 1 \end{array} \right)\left (\begin{array}{cc}1 & 1 \\
1 & 2 \end{array} \right)^k(v_n), \,\, k\geq1, \\
v_{n+1,k}^*&=\left (\begin{array}{cc}1 & 1 \\
1 & 2 \end{array} \right)\left (\begin{array}{cc}1 & 1 \\
1 & 2 \end{array} \right)^k(v_n), \,\, k\geq1, \quad \text{and}\\
v_{n+1,k}&=\left (\begin{array}{cc} 0 & 1 \\
1 & 2 \end{array} \right)\left (\begin{array}{cc}1 & 1 \\
1 & 2 \end{array} \right)^k(v_n),\,\, k\geq1.
\end{aligned}
\right.
$$
We note that $\left (\begin{array}{cc}1 & 1 \\
1 & 2 \end{array} \right)$ is the square of $\left (\begin{array}{cc}0 & 1 \\
1 & 1 \end{array} \right)$, the matrix that maps a pair of consecutive Fibonacci numbers $(F_{i-1},F_i)$ to $(F_i,F_{i+1})$; we put $F_0=0$ and $F_1=1$. It follows that 
 $$
 \lim_{k \to \infty} \left (\begin{array}{cc}1 & 1 \\
1 & 2 \end{array} \right)^k(v_n)=\lim_{k \to \infty}  \frac {F_{2k-1}v_n+F_{2k}}{F_{2k}v_n+F_{2k+1}}.
$$
 Now suppose $\lim_{k \to \infty}  \frac {F_{2k-1}v_n+F_{2k}}{F_{2k}v_n+F_{2k+1}}=L$. Then 
$$ 
L = \lim_{k \to \infty}  \frac 1{\tfrac{F_{2k}v_n+F_{2k+1}}{F_{2k-1}v_n+F_{2k}}}=\lim_{k \to \infty}  \frac 1{1+\tfrac{F_{2k-2}v_n+F_{2k-1}}{F_{2k-1}v_n+F_{2k}}}=\frac1{1+L},
$$ from which we derive 
$$
\lim_{k \to \infty} \left (\begin{array}{cc}1 & 1 \\
1 & 2 \end{array} \right)^k(v_n)=g, 
$$
whence 
$$
\left \{
\begin{aligned}
\lim_{k \to \infty} v_{n,k}^*&=\tfrac g{g+1}=g^2, \\
\lim_{k \to \infty} v_{n+1,k}^*&=\tfrac{g+1}{g+2}=g \quad \text{and} \\
\lim_{k \to \infty} v_{n+1,k}&=\tfrac1{g+2}=g^2.
\end{aligned} 
\right.
$$
We conclude that 
\begin{align*}
\Omega_\alpha&=\Omega_g \setminus \left ([\alpha,g]\times[0,\tfrac12) \bigcup_{k=1}^{\infty} [-g^2,\tfrac{1-2\alpha}{\alpha}]\times \left (\frac{F_{2k+1}}{F_{2k+3}},\frac{F_{2k-1}}{F_{2{k+1}}} \right ) \right )\\
&\bigcup_{k=1}^{\infty} \left ([\alpha-1,-g^2]\times \left [\frac{F_{2(k-1)}}{F_{2k}},\frac{F_{2k}}{F_{2(k+1)}} \right )\cup [\tfrac{1-2\alpha}{\alpha-1},\alpha)\times \left [\frac{F_{2k}}{F_{2k+1}},\frac{F_{2k+2}}{F_{2k+3}}\right )\right ),  \\
\end{align*}
from which we derive the following results, again obtained by Nakada in \cite{Na} (cf.~\cite{K}), illustrated by Figures \ref{fig: omega.52infty} and \ref{fig: omega.1/2}.\\
\begin{Theorem}
Let $\alpha \in (\frac12,g]$. Then 
\begin{align*}
\Omega_\alpha=&[\alpha-1,-g^2)\times[0,g^2)\,\cup\,[-g^2,\tfrac{1-2\alpha}{\alpha}]\times[0,g^2]\\
&\cup\,(\tfrac{1-2\alpha}{\alpha},0)\times[0,\tfrac12]\,\cup\,[0,\tfrac{1-2\alpha}{\alpha-1}]\times[0,\tfrac12)\,\cup\,(\tfrac{1-2\alpha}{\alpha-1},\alpha)\times[0,g).
\end{align*}
\end{Theorem}
Moreover, we have
$$
\Omega_{1/2}=[-\tfrac12,0]\times[0,g^2]\,\cup\,(0,\tfrac12)\times[0,g).
$$

\begin{figure}[!htb]
\minipage{0.48\textwidth}
$$
\begin{tikzpicture}[scale =6] 
 \draw (-.48,0) -- (.52,0);
 \draw (-.48,0) -- (-.48,.382);
 \draw [dashed] (-.48,.382) -- (-.382,.382);
 \draw [dashed] (-.382,0) -- (-.382,.382);
    \draw [dashed] (.083,0) -- (.083,.5);
      \draw (.083,.5) -- (.083,.62);
     \draw [dashed] (.083,.62) -- (.52,.62);
   \draw [dashed] (.52,0) -- (.52,.62);
  \draw (0,.4) -- (0,.5);
  \draw [dashed] (0,.5) -- (.083,.5); 
    \draw [dashed] (-.077,0) -- (-.077,.5);
  \draw (-.382,.382) -- (-.077,.382);
  \draw (-.077,.5) -- (0,.5);
 \draw (0,0) -- (0,.4);
  \node at (-.5,-.05) {$\alpha-1$};
   \node at (-.53,.382){$g^2$};
          \node at (-.11,.5){$\tfrac12$};
          \node at (-.35,-.05) {$-g^2$};
     \node at (.06,.618){$g$};
     \node at (.1,-.05) {$\tfrac{1-2\alpha}{\alpha-1}$};
       \node at (.52,-.05) {$\alpha$};
   \node at (-.11,-.05){$\tfrac{1-2\alpha}{\alpha}$};
   \node at (0,-.05) {0};
 \end{tikzpicture}
 $$
 \caption[$\Omega_{\alpha}$, $\alpha=0.52$.] {\label{fig: omega.52infty}
$\Omega_\alpha$, $\alpha=0.52$.}
\endminipage\hfill
\minipage{0.48\textwidth}
$$
\begin{tikzpicture}[scale =6] 
 \draw (-.5,0) -- (.5,0);
 \draw  [dashed] (.5,0) -- (.5,.618);
 \draw(-.5,0) -- (-.5,.382);
   \draw [dashed] (-.382,0) -- (-.382,.382);
     \draw [dashed] (0,.62) -- (.5,.62);
   \draw (0,0) -- (0,.618);
    \draw [dashed] (-.382,.382) -- (0,.382);
    \draw (-.5,.382) -- (-.382,.382);
 \draw (0,0) -- (0,.4);
  \node at (-.36,-.05){$-g^2$};
  \node at (-.5,-.05) {$-\tfrac12$};
   \node at (-.53,.382){$g^2$};
     \node at (-.03,.62){$g$};
       \node at (.5,-.05) {$\tfrac12$};
      \node at (0,-.05) {0};
  \end{tikzpicture}
 $$
 \caption[$\Omega_{\frac12}$] {\label{fig: omega.1/2}
$\Omega_{\frac12}$.}
\endminipage
\end{figure}
In order to find the effect of the current transformation on the sequence of convergents, we recall that we rewrote 
\begin{align*}
&x=[0;\varepsilon_1a_1,\dots,\varepsilon_na_n,2,-a_{n+2},\dots] \quad {\text{as}}\\
&x=[0;\varepsilon_1a_1,\dots,\varepsilon_n(a_n+1),-2,(a_{n+2}-1),\dots].
\end{align*}
Since $p_k^*=p_k$ and $q_k^*=q_k$, $k<n$, we find, applying (\ref{pn en qn}) and omitting the straightforward calculations,
$$
\left \{
\begin{aligned}
p_n^*&=p_n+p_{n-1};\\
p_{n+1}^*&=p_{n+1};\\
p_{n+2}^*&=p_{n+2}.\\
\end{aligned}
\right.
$$
Similarly, we find $q_n^*=q_n+q_{n-1}$, $q_{n+1}^*=q_{n+1}$ and $q_{n+2}^*=q_{n+2}$. It would seem that in the transformation from $\Omega_g$ to $\Omega_\alpha$, $\alpha \in [1/2,g)$, the convergent $p_n/q_n$ is replaced by the so-called {\it{mediant}} $(p_n+p_{n-1})/(q_n+q_{n-1})$. We emphasise, however, that this is only a mediant with respect to a small interval with upper bound $g$; actually, this `$\alpha${\it{-mediant}}' is a convergent of the regular continued fraction. This comes as no surprise: from \cite{K} we know that regular convergents only disappear and are not replaced by mediants when moving from the regular expansion to the $\alpha$-expansion for $1/2\leq \alpha\leq g$. The regular convergent that replaces another one in the transformation from $\Omega_g$ to $\Omega_\alpha$, $\alpha \in [1/2,g)$ had previously been 'singularised away' in the transformation from $\Omega_1$ to $\Omega_g$. We will illustrate this with an example.\smallskip

Let $x=(\sqrt{17}-3)/4=1/(3+(\sqrt{17}-3)/2)=0.28077\dots$; in the second, more intricate expression for $x$, the number $(\sqrt{17}-3)/2$ is explicitly displayed for reasons that will show instantly. We distinguish 
$$
\begin{cases}
x_\alpha=[0;\overline{3,1,1}], \quad &\alpha \in (g,1];\\
x_\alpha=[0;{3,\overline{2,-4}}], \quad &\alpha \in ((\sqrt{17}-3)/2,g];\\
x_\alpha=[0;\overline{4,-2}], \quad &\alpha \in (1/2,(\sqrt{17}-3)/2].
\end{cases}
$$

If we focus at the first partial quotients, we have 
$$
\begin{cases}
 {\text{if}} \,\, x_\alpha=[0;3,1,1,3,1,1\dots], &{\text{then}} \,\, (\tfrac {p_0}{q_0},\dots,\tfrac{p_6}{q_6})=(\tfrac 01,\tfrac 13,\tfrac 14, \tfrac 27, \tfrac 7{25}, \tfrac 9{32}, \tfrac {16}{57});\\
 {\text{if}} \,\, x_\alpha=[0;3,2,-4,2,\dots], &{\text{then}} \,\, (\tfrac {p_0}{q_0},\dots,\tfrac{p_4}{q_4})=(\tfrac 01,\tfrac 13,\tfrac 27, \tfrac 7{25}, \tfrac {16}{57});\\
 {\text{if}} \,\, x_\alpha=[0;4,-2,4,-2,\dots], &{\text{then}}\,\,  (\tfrac {p_0}{q_0},\dots,\tfrac{p_4}{q_4})=(\tfrac 01, \tfrac 14, \tfrac 27, \tfrac 9{32}, \tfrac {16}{57}).
\end{cases}
$$

What happens is that we create $[0;3,2,-4,2,\dots]$ from $[0;3,1,1,3,1,1\dots]$ by singularising the second $1$ from each pair $(1,1)$, in which process some of the (regular) convergents get lost, in accordance with what we saw in Section \ref{The case alpha in (g,1]}. When creating $[0;4,-2,4,-2,\dots]$ from $[0;3,2,-4,2,\dots]$, we first insert a $1$ between each $2$ followed by $-4$, so as to retrieve $[0;3,1,1,3,1,1]$ (and its convergents), and subsequently singularise the first $1$ from each pair $(1,1)$, so as to lose other convergents of the regular case. Generally, we conclude that in the transformation from $\Omega_g$ to $\Omega_\alpha$, $\alpha \in [1/2,g)$, some regular convergents $p_n/q_n$ are replaced by others. \smallskip

Although direct singularisation, as used in for instance \cite{K}, would be more straightforward, it is limited to $\alpha \geq 1/2$, while our approach enables a decrease beyond $\alpha=1/2$, as we will show in the next sections. Moreover, it shows how to transform a domain $\Omega_\alpha$ into another domain $\Omega_{\alpha'}$ in a uniform manner, where $\alpha'$ is smaller than but close to $\alpha$, for any $\alpha \in (g^2,1]$. On top of that, in Section \ref{ergodic systems} we will go into the isomorphism between $\Omega_\alpha$ and $\Omega_g$ for $g^2 \leq \alpha < g$.

\section{The case $\alpha\in (\sqrt{2}-1,\tfrac12]$}
\label{The case alpha in (sqrt{2}-1,tfrac12)}
In this section we will show how to derive $\Omega_{\alpha}$ from $\Omega_{1/2}$, with $\sqrt{2}-1\leq \alpha <1/2$. Let $x \in [-1/2,1/2)$ and let $n \geq 0$ be the smallest integer for which $t_n \in [\alpha,1/2]$, i.e., for which $(t_n,v_n) \in R_\alpha:=[\alpha,1/2]\times[0,g)$. We already know that in this case $x=[0;\varepsilon_1a_1,\dots,\varepsilon_na_n,2,a_{n+2},\dots]$. Since the continued fraction map of the related dynamical system is $T_{1/2}$, for $\alpha\in (\sqrt{2}-1,1/2]$ we have 
$$
a_2(\alpha)=\left \lfloor \tfrac1{T_{1/2}(\alpha)} +1-\tfrac12 \right \rfloor = \left \lfloor \tfrac1{\tfrac1\alpha-2}+\frac12\right \rfloor \geq \left \lfloor \tfrac1{\tfrac1{\sqrt{2}-1}-2}+\frac12\right \rfloor =\left \lfloor \sqrt{2}+\frac32\right \rfloor =2,
$$
so in this case 
$$
x=[0;\varepsilon_1a_1,\dots,\varepsilon_na_n,2,a_{n+2},\dots], \, a_{n+2} \geq2 .
$$
The difference with the case in the previous section is merely the sign of $a_{n+2}$, which shows in the great resemblance of equations that we need to construct $\Omega_\alpha$. Applying insertion and singularisation in a similar way, we find
$$
x=[0;\varepsilon_1a_1,\dots,\varepsilon_n(a_n+1),-2,-(a_{n+2}+1),\dots]
$$ 
this time. The other equations are exactly the same as in the previous section: $(t_{n+1},v_{n+1})=((1-2t_n)/t_n,1/(2+v_n))$, $
(t_n^*,v_n^*)=(t_n-1,v_n/(1+v_n))$, $(t_{n+1}^*, v_{n+1}^*)=((1-2t_n)/(t_n-1),(v_n+1)/(v_n+2))$ and $(t_{n+2}^*, v_{n+2}^*)=(t_{n+2}, v_{n+2})$.\\
Using similar calculations as in the previous section, we find that $R_\alpha=[\alpha,1/2)\times[0,g)$ implies $\mathcal T_{1/2}(R_\alpha) = (0,(1-2\alpha)/\alpha) \times (g^2,1/2]$; see Figure \ref{fig: omega43}, where we have taken $\alpha=0.43$ for illustrative purposes. The first region to be added is: $A_\alpha=[\alpha-1,-\tfrac12]\times[0,g^2)$; the second is $\mathcal T_{1/2}(A_\alpha)=[(1-2\alpha)/(\alpha-1),0]\times[1/2,g)$. An important difference with the previous case is that $(R_\alpha \cup \mathcal T_{1/2}(R_\alpha))\cap((A_\alpha \cup \mathcal T_{1/2}(A_\alpha))=\emptyset$, so all we have to do is remove $R_\alpha$ and $\mathcal T_{1/2}(R_\alpha)$ from $\Omega_{1/2}$ and add  $A_\alpha$ and $\mathcal T_{1/2}(A_\alpha)$ to it. Doing so, we find
\begin{Theorem}
Let $\alpha \in (\sqrt{2}-1,1/2]$. Then 
$$
\Omega_\alpha=[\alpha-1,\alpha)\times[0,g^2)\,\cup\,[\tfrac{1-2\alpha}{\alpha},\alpha)\times[g^2,\tfrac12)\, \cup\,[\tfrac{1-2\alpha}{\alpha-1},\alpha)\times[\tfrac12,g).
$$
\end{Theorem}
Since $(1-2\alpha)/\alpha$ equals $\alpha$ for $\alpha=\sqrt{2}-1$, at $\alpha=\sqrt{2}-2$ the domain $\Omega_\alpha$ is split in two parts:
$$
\Omega_{\sqrt{2}-1}=[\sqrt{2}-2,\sqrt{2}-1)\times[0,g^2)\,\cup\,[\tfrac12\sqrt{2}-1,\sqrt{2}-1)\times[\tfrac12,g);
$$
see Figure \ref{fig: omegasqrt2-1}.
\begin{figure}[!htb]
\minipage{0.48\textwidth}
$$
\begin{tikzpicture}[scale =6] 
  \draw[draw=white,fill=gray!20!white] 
 plot[smooth,samples=100,domain=.43:.5] (\x,0) --
 plot[smooth,samples=100,domain=.5:.43] (\x,0.618);
 \draw[draw=white,fill=gray!20!white] 
 plot[smooth,samples=100,domain=0:.326] (\x,.382) --
 plot[smooth,samples=100,domain=.326:0] (\x,.5);
 \draw (-.57,0) -- (.5,0);
 \draw [dashed] (.5,0) -- (.5,.618);
 \draw (-.5,0) -- (-.5,.382);
     \draw [dashed] (0,.618) -- (.5,.618);
     \draw (-.57,0) -- (-.57,.382);
     \node at (-.53,.21) {$A_\alpha$};   
    \draw [dashed] (-.57,.382) -- (-.5,.382);
   \draw (0,0) -- (0,.618);
    \draw [dashed] (-.382,.382) -- (0,.382);
    \draw (.-.246,.5) -- (0,.5);
    \draw [dashed] (.-.246,.618) -- (0,.618);
    \draw (-.246,.5) -- (-.246,.618);
\draw [dashed] (0,.382) -- (.326,.382);
\draw [dashed] (.43,0) --(.43,.618);
\draw (0,.5) -- (.326,.5);
\draw (.326,.382) -- (.326,.5);
      \draw [dashed]  (-.246,0) -- (-.246,.5);
    \draw [dashed] (-.5,.382) -- (-.382,.382);
 \draw [dashed] (.326,0) --(.326,.382);
  \node at (-.47,-.05) {$-\tfrac12$};
   \node at (.43,-.05) {$\alpha$};
    \node at (.33,-.05) {$\tfrac{1-2\alpha}{\alpha}$};
     \node at (-.28,.5) {$\tfrac12$};
   \node at (-.6,.382){$g^2$};
   \node at (-.62,-.05){$\alpha-1$};
   \node at (-.246,-.05) {$\tfrac{1-2\alpha}{\alpha-1}$};
   \node at (0,.66) {};
     \node at (-.28,.62){$g$};
       \node at (.5,-.05) {$\tfrac12$};
       \node at (.47,.31) {$R_\alpha$};
          \node at (.16,.44) {$\mathcal T_{1/2}(R_\alpha)$};
           \node at (-.125,.55) {$\mathcal T_{1/2}(A_\alpha)$};
      \node at (0,-.06) {0};
 \end{tikzpicture}
 $$
 \caption[$\Omega_{\tfrac12}, \alpha=0.43.$] {\label{fig: omega43}
$\Omega_{\tfrac12}, \alpha=0.43.$}
\endminipage\hfill
\minipage{0.48\textwidth}
$$
\begin{tikzpicture}[scale =6] 
 \draw (-.586,0) -- (.414,0);
 \draw [dashed] (.414,0) -- (.414,.382);
  \draw [dashed] (.414,.5) -- (.414,.618);
  \draw [dashed] (-.293,.62) -- (.414,.62);
  \draw [dashed] (-.293,0) -- (-.293,.5);
     \draw (-.586,0) -- (-.586,.382);
   \draw [dashed] (0,0) -- (0,.618);
    \draw [dashed] (-.586,.382) -- (.414,.382);
    \draw (.-.293,.5) -- (.414,.5);
     \draw (-.293,.5) -- (-.293,.618);
         \node at (-.32,.5) {$\tfrac12$};
      \node at (-.293,-.05) {$\tfrac12\sqrt{2}-1$};
   \node at (-.62,.382){$g^2$};
   \node at (-.6,-.05){$\sqrt{2}-2$};
     \node at (-.32,.62){$g$};
       \node at (.414,-.05) {$\sqrt{2}-1$};
                 \node at (0,-.05) {0};
   \end{tikzpicture}
 $$
 \caption[$\Omega_{\sqrt{2}-1}.$] {\label{fig: omegasqrt2-1}
$\Omega_{\sqrt{2}-1}.$}
\endminipage
\end{figure}

With a great similarity of calculations as in the previous section, we find that in this case also
$$
\left \{
\begin{aligned}
p_n^*&=p_n+p_{n-1};\\
p_{n+1}^*&=p_{n+1};\\
p_{n+2}^*&=p_{n+2}.\\
\end{aligned}
\right.
$$
The difference with the case $\alpha \in (1/2,g]$ is that now the convergent $p_n/q_n$ is really replaced by the mediant $(p_n+p_{n-1})/(q_n+q_{n-1})$ in the sense of the regular case. 

\section{The case $\alpha\in (\tfrac{\sqrt{10}-2}3,\sqrt{2}-1]$}
\label{The case alpha in (tfrac{sqrt{10}-2}3,sqrt{2}-1]}
 In Section \ref{The case alpha in (frac12,g]} we noticed that $\sqrt{2}-1_{\sqrt{2}-1}=[0;3,\overline{-2,-4}]$; actually, $\sqrt{2}-1$ is the largest $\alpha \in [0,1]$ such that no partial quotients $1$ or $2$ (with plus sign) occur in $\alpha_\alpha$. In this section we want to remove points $(t_n,v_n)$ with $t_n=[0;3,\varepsilon_{n+2}a_{n+2},\dots]$, that is, points associated with numbers $x=[0;\varepsilon_1a_1,\dots,\varepsilon_na_n,3,\varepsilon_{n+2}a_{n+2},\dots]$. Similar to the cases in Sections \ref{The case alpha in (frac12,g]} and \ref{The case alpha in (sqrt{2}-1,tfrac12)}, the removal of a $3$ requires compensated insertion, for which we need a partial quotient $2$, as explained in Section \ref{Introduction}. Since we removed all partial quotients $2$ with a plus sign, we are dealing with the following numbers:
$$
x=[0;\varepsilon_1a_1,\dots,\varepsilon_na_n,3,-2,\varepsilon_{n+3}a_{n+3},\dots].
$$
Note that uncompensated insertion would leave us with partial quotients $1$ or $+2$, involving points outside $\Omega_\alpha$. Because of the necessity of compensated insertion, the next interval to explore is not simply $\Delta_\alpha(3)=((\sqrt13-3)/2,\sqrt{2}-1]$ but the much smaller $\Delta_\alpha(3,-2)$. To find the boundaries of this interval, we determine $(\alpha', \alpha]$, $\alpha' <\alpha \leq \sqrt{2}-1$, such that $t_n \in (\alpha', \alpha]$ implies $t_n=[0;3,-2,\dots]$. If $\alpha' \in [0,\sqrt{2}-1)$ and $\alpha \in (\alpha',\sqrt{2}-1]$, then
\begin{align*}
2 &= \left \lfloor \frac1{3-\frac1{\sqrt{2}-1}}+2-\sqrt{2}\right \rfloor \leq  \left \lfloor\frac1{3-\frac1{\alpha}}+1-\alpha \right \rfloor \leq a_2(\alpha'_\alpha) \\
&=\left \lfloor \frac1{3-\frac1{\alpha'}}+1-\alpha \right \rfloor \leq \left \lfloor \frac1{3-\frac1{\alpha'}}+1-\alpha' \right \rfloor, 
\end{align*}
from which we derive $a_2(\alpha'_\alpha)=2$ if and only if $1/(3-1/\alpha')+1-\alpha' <3$. Given that $\alpha' \in [0,\sqrt{2}-1)$), this is the case if and only if $(\sqrt{10}-2)/3<\alpha'<\sqrt{2}-1$; indeed, $(\sqrt{10}-2)/3_{(\sqrt{10}-2)/3}=[0;3,\overline{-3,-2,-3,-4}]$. So $\Delta_\alpha(3,-2)=((\sqrt{10}-2)/3,\sqrt{2}-1]$. Since $\sqrt{2}-1=[0;3,\overline{-2,-4}]$ and $\lim_{a \to \infty}[0;3,-2,-a]=\tfrac25$, we will first show how to construct $\Omega_{\alpha}$ from $\Omega_{\sqrt{2}-1}$, with $2/5 \leq \alpha < \sqrt{2}-1$, which will prove to be more complicated than the work we have done so far. After that, the step from $\Omega_{2/5}$ to $\Omega_{(\sqrt{10}-2)/3}$ will prove to be fairly small.\smallskip

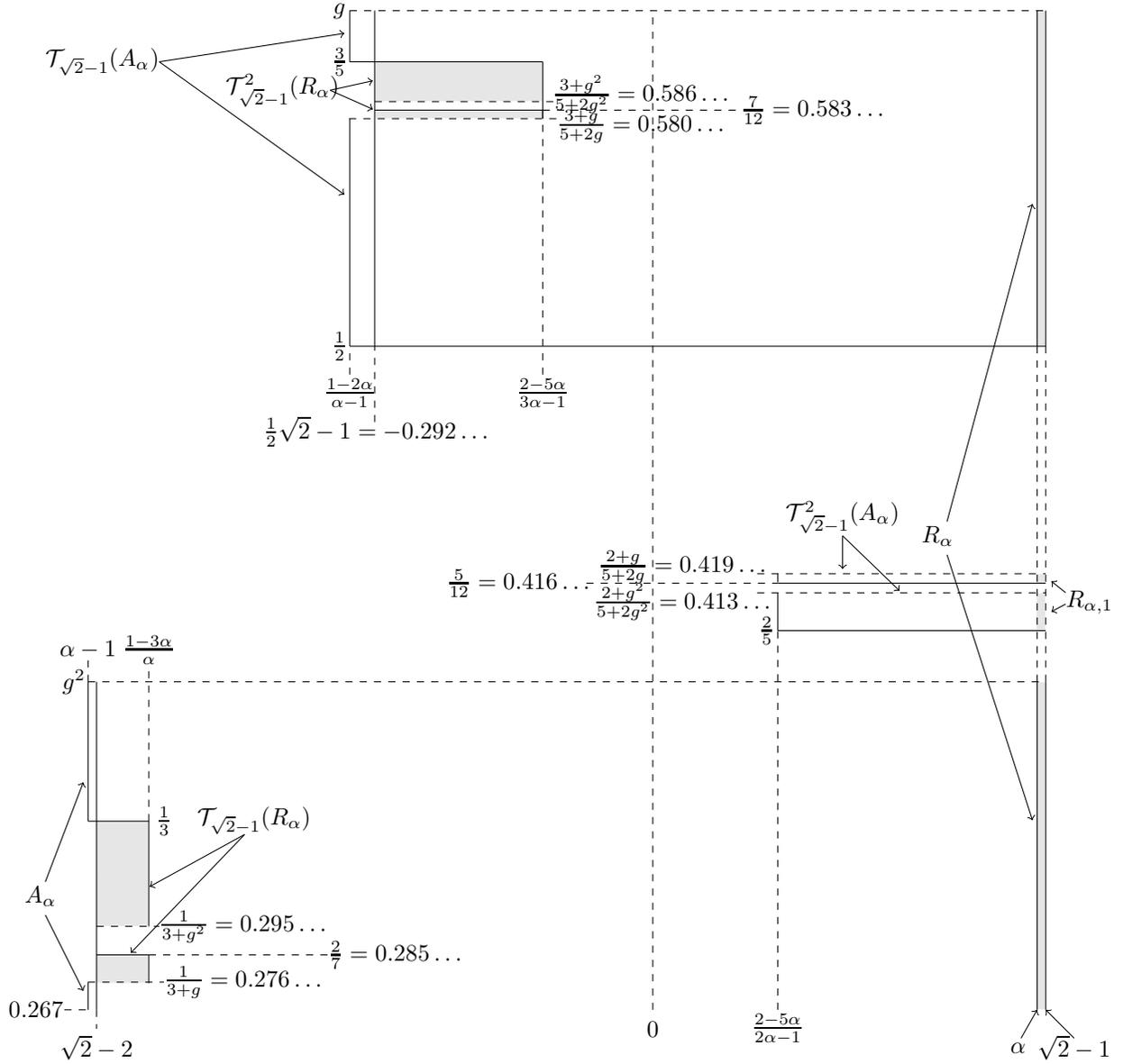
\begin{figure}[!htb]
$$
\begin{tikzpicture}[scale =14] 
  \draw[draw=white,fill=gray!20!white] 
 plot[smooth,samples=100,domain=.405:.414] (\x,0) --
 plot[smooth,samples=100,domain=.414:.405] (\x,0.346);
  \draw[draw=white,fill=gray!20!white] 
 plot[smooth,samples=100,domain=.405:.414] (\x,.7) --
 plot[smooth,samples=100,domain=.414:.405] (\x,1.054); 
    \draw[draw=white,fill=gray!20!white] 
 plot[smooth,samples=100,domain=-.586:-.531] (\x,.028) --
 plot[smooth,samples=100,domain=-.531:-.586] (\x,.058);
\draw[draw=white,fill=gray!20!white] 
 plot[smooth,samples=100,domain=-.586:-.531] (\x,.088) --
 plot[smooth,samples=100,domain=-.531:-.586] (\x,.199);
     \draw[draw=white,fill=gray!20!white] 
 plot[smooth,samples=100,domain=-.293:-.116] (\x,.94) --
 plot[smooth,samples=100,domain=-.116:-.293] (\x,.949);
   \draw[draw=white,fill=gray!20!white] 
 plot[smooth,samples=100,domain=-.293:-.116] (\x,.958) --
 plot[smooth,samples=100,domain=-.116:-.293] (\x,1);
  \draw[draw=white,fill=gray!20!white] 
 plot[smooth,samples=100,domain=.405:.414] (\x,.4) --
 plot[smooth,samples=100,domain=.414:.405] (\x,0.4399);
  \draw[draw=white,fill=gray!20!white] 
 plot[smooth,samples=100,domain=.405:.414] (\x,.45) --
 plot[smooth,samples=100,domain=.414:.405] (\x,.46);
 \draw (.414,0) -- (.414,.346);
  \draw [dashed] (.414,.346) -- (.414,.7);
 \draw (.405,0) -- (.405,.346);
  \draw [->]  (.388,-.03)-- (.405,0);
   \draw [dashed] (.405,.346) -- (.405,.7);
 \draw (.405,.7) -- (.405,1.054);
  \draw [dashed] (-.586,.088) -- (-.52,.088);
   \draw [dashed] (-.293,.62) -- (-.293,.7);
    \draw (-.586,.199) -- (-.531,.199);
    \draw [dashed] (-.531,.199) -- (-.531,.363);
         \draw (-.531,.088) -- (-.531,.199);
    \draw (-.595,0) -- (-.595,.0292);
     \draw (-.595,.199) -- (-.595,.346);
     \draw (-.595,.199) -- (-.586,.199);
          \draw [dashed] (-.595,.346) -- (-.595,.37);
 \draw [dashed] (-.595,.0292) -- (-.515,.0292);
 \draw (.1316,.4) -- (.414,.4);
  \draw (.1316,.4) -- (.1316,.4399);
  \draw [dashed] (.11,.4399) -- (.414,.4399); 
  \draw [->] (.31,.48) -- (.4,.2);
   \draw [->] (.31,.515) -- (.4,.85);
\draw [->] (.435,.44) -- (.42,.452);
   \draw [->] (.435,.428) -- (.42,.42);
\draw [->] (-.64,.135) -- (-.6,.24);
  \draw [->] (-.64,.1) -- (-.6,.02);
\draw [->] (-.34,.97) -- (-.295,.98);
 \draw [->] (-.34,.97) -- (-.295,.95);
 \draw [->] (-.43,.185) -- (-.53,.13);
 \draw [->] (-.43,.185) -- (-.55,.06);
 \draw [->] (-.52,1) -- (-.325,1.03);
 \draw [->] (-.52,1) -- (-.325,.86);
\draw [->] (.2,.5) -- (.2,.465);
 \draw [->] (.2,.5) -- (.26,.442);
  \draw [dashed] (-.07,.45) -- (.1316,.45);
    \draw (.1316,.45) -- (.414,.45);
  \draw (.1316,.451) -- (.1316,.46);
 \draw [dashed] (.11,.46) -- (.414,.46);
    \draw (-.586,.058) -- (-.531,.058);
    \draw [dashed] (-.531,.058) -- (-.34, .058);
    \draw (-.531,.028) -- (-.531,.058);
 \draw [dashed] (-.293,.958) -- (-.1,.958);
    \draw (-.3193,1) -- (-.116,1);
        \draw (-.3193,1) -- (-.3193,1.054);
    \draw (-.116,.94) -- (-.116,1);
        \draw [dashed] (-.116,.67) -- (-.116,.94);
 \draw [dashed] (-.32,.94) -- (-.1,.94);
    \draw (-.293,.949) -- (-.116,.949);
    \draw [dashed] (-.116,.949) -- (.09,.949);
    \draw [dashed] (-.62,0) -- (-.595,0);
    \draw (-.116,.94) -- (-.116,.949);
  \draw (.414,.7) -- (.414,1.054);
  \draw (.444,-.03) [->] -- (.414,0);
  \draw [dashed] (-.3193,1.054) -- (.414,1.054);
  \draw [dashed] (.1316,0) -- (.1316,.4);
  \draw [dashed] (-.3193,.67) -- (-.3193,.7);
     \draw (-.586,0) -- (-.586,.346);
     \draw [dashed] (-.586,0) -- (-.586,.-.02);
   \draw [dashed] (0,0) -- (0,1.054);
    \draw [dashed] (-.595,.346) -- (.414,.346);
    \draw (.-.3193,.7) -- (.414,.7);
    \draw (-.3193,.7) -- (-.3193,.94);
     \draw (-.293,.7) -- (-.293,1.054);
     \node at (-.43,.088) {$\tfrac1{3+g^2}=0.295\dots$};
       \node at (-.011,.966) {$\tfrac{3+g^2}{5+2g^2}=0.586\dots$};
   \node at (.17,.949) {$\tfrac7{12}=0.583\dots$};
     \node at (-.011,.932) {$\tfrac{3+g}{5+2g}=0.580\dots$};
\node at (-.116,.65) {$\tfrac{2-5\alpha}{3\alpha-1}$};
      \node at (-.27,.058) {$\tfrac27=0.285\dots$};
 \node at (-.43,.0292) {$\tfrac1{3+g}=0.276\dots$};
 \node at (-.33,1) {$\tfrac35$};
      \node at (-.33,.7) {$\tfrac12$};
          \node at (-.29,.61) {$\tfrac12\sqrt{2}-1=-0.292\dots$};
          \node at (-.515,.199) {$\tfrac13$};
          \node at (.12,.4) {$\tfrac25$};
          \node at (.033,.431) {$\tfrac{2+g^2}{5+2g^2}=0.413\dots$};
          \node at (-.14,.45) {$\tfrac5{12}=0.416\dots$};
           \node at (.033,.467) {$\tfrac{2+g}{5+2g}=0.419\dots$};
   \node at (-.61,.346){$g^2$};
   \node at (-.586,-.04){$\sqrt{2}-2$};
   \node at (-.531,.38) {$\tfrac{1-3\alpha}{\alpha}$};
   \node at (.385,-.04) {$\alpha$};
   \node at (-.595,.38) {$\alpha-1$};
       \node at (-.33,1.054){$g$};
       \node at (.444,-.04) {$\sqrt{2}-1$};
       \node at (-.65,0) {0.267};
       \node at (0,-.02) {0};
        \node at (.3,.5) {$R_\alpha$};
       \node at (.46,.43) {$R_{\alpha,1}$};
      \node at (-.645,.12) {$A_\alpha$};
        \node at (-.39,.97) {$\mathcal T_{\sqrt{2}-1}^2(R_\alpha)$};
           \node at (-.42,.2) {$\mathcal T_{\sqrt{2}-1}(R_\alpha)$};
             \node at (-.58,1) {$\mathcal T_{\sqrt{2}-1}(A_\alpha)$};
                 \node at (.2,.52) {$\mathcal T_{\sqrt{2}-1}^2(A_\alpha)$};
                 \node at (.1316,-.02) {$\tfrac{2-5\alpha}{2\alpha-1}$};
                 \node at (-.3193,.65) {$\tfrac{1-2\alpha}{\alpha-1}$};
   \end{tikzpicture}
 $$
 \caption[$\Omega_{\alpha,1}$, $\alpha=0.405.$] {\label{fig: omega.4051}
$\Omega_{\alpha,1}$, $\alpha=0.405$. The lower rectangle $[\alpha-1,\sqrt{2}-1] \times [0,.267)$ has been omitted.}
\end{figure}

Let $\alpha \in (2/5, \sqrt{2}-1]$, $x \in [\sqrt{2}-2,\sqrt{2}-1)$ and $n \geq 0$ be the smallest integer for which $t_n \in [\alpha,\sqrt{2}-1)$, i.e., for which 
\begin{equation}\label{tvcase4}
(t_n,v_n) \in R_\alpha:=[\alpha,\sqrt{2}-1]\times[0,g^2)\cup [\alpha,\sqrt{2}-1]\times[\tfrac12,g).
\end{equation}
Considering the fact that $\sqrt{2}-1_{\sqrt{2}-1}=[0;3,\overline{-2,-4}]$, it is not hard to see that in this case 
$$
x=[0;\varepsilon_1a_1,\dots,\varepsilon_na_n,3,-2,-a_{n+3},\dots], \,a_{n+3} \geq4.
$$
First we insert $1$ between $3$ and $-2$, yielding
$$
x=[0;\varepsilon_1a_1,\dots,\varepsilon_na_n,2,1, 1,-a_{n+3},\dots], \,a_{n+3} \geq4.
$$
Next we insert $-1$ between $a_n$ and $2$, yielding
$$
x=[0;\varepsilon_1a_1,\dots,\varepsilon_n(a_n+1),-1,1,1, {\bf{1}},-a_{n+3},\dots], \,a_{n+3} \geq4.
$$
Then we singularise the last $+1$, yielding
$$
x=[0;\varepsilon_1a_1,\dots,\varepsilon_n(a_n+1),-1,{\bf{1}},2, -2,a_{n+3}-1,\dots], \,a_{n+3} \geq4.
$$
Finally we singularise the remaining $+1$, so as to get
$$
x=[0;\varepsilon_1a_1,\dots,\varepsilon_n(a_n+1),-2,-3,a_{n+3}-1,\dots], \,a_{n+3} \geq4.
$$
For the first time we are confronted with the singularised expansion differing at four places from the original expansion, necessitating one more stage of calculating successive $(t,v)$-pairs (in the sense of the previous sections). Omitting the straightforward calculations, we find  
\begin{align*}
\qquad \qquad \qquad \qquad (t_{n+1},v_{n+1})&=(\tfrac{1-3t_n}{t_n},\tfrac1{3+v_n}),  &(t_{n+2},v_{n+2})&=(\tfrac{2-5t_n}{3t_n-1},\tfrac{v_n+3}{2v_n+5}),\qquad \qquad \qquad \qquad\\
(t_n^*,v_n^*)&=(t_n-1,\tfrac{v_n}{1+v_n}), &(t_{n+1}^*, v_{n+1}^*)&=(\tfrac{1-2t_n}{t_n-1},\tfrac{v_n+1}{v_n+2}), \\
(t_{n+2}^*, v_{n+2}^*)&=(\tfrac{2-5t_n}{2t_n-1},\tfrac{v_n+2}{2v_n+5}), &(t_{n+3}^*, v_{n+3}^*)&=(t_{n+3}, v_{n+3}).
\end{align*}
We see that in the current situation (i.e. for the current values of $\alpha$), constructing $\Omega_\alpha$ is much more complicated than in the previous sections. There is one more step of removing and adding regions, and these regions are all split in two disjoint parts. In particular, $R_{\alpha,1}:=\mathcal T_{\sqrt{2}-1}^2(A_\alpha) \cap R_\alpha \neq \emptyset$, similar to Section~\ref{The case alpha in (frac12,g]}. Since (\ref{tvcase4}) holds, we can draw $\Omega_{\alpha,1}$, defined similarly as in Section \ref{The case alpha in (frac12,g]}. We will first sketch roughly what is removed from $\Omega_{\sqrt{2}-1}$ (see Figure \ref{fig: omega.4051, removed, rough}; the removed parts are in grey), consisting of 
\begin{align*}
R_\alpha&=[\alpha,\sqrt{2}-1]\times[0,g^2)\cup [\alpha,\sqrt{2}-1]\times[\tfrac12,g),\\
\mathcal T_{\sqrt{2}-1}(R_\alpha)&=[\sqrt{2}-2,\tfrac{1-3\alpha}{\alpha}]\times(\tfrac1{3+g^2},\tfrac13]\cup[\sqrt{2}-2,\tfrac{1-3\alpha}{\alpha}]\times(\tfrac1{3+g},\tfrac27] \quad {\text{and}}\\
\mathcal T_{\sqrt{2}-1}^2(R_\alpha)&=[\tfrac12\sqrt{2}-1,\tfrac{2-5\alpha}{3\alpha-1}]\times(\tfrac{3+g^2}{5+2g^2},\tfrac35]\cup[\tfrac12\sqrt{2}-1,\tfrac{1-3\alpha}{\alpha}]\times(\tfrac{3+g}{5+2g},\tfrac7{12}].
\end{align*}
The removal of $R_\alpha$ does not show in the resulting figure (it is `cut off' $\Omega_{\sqrt{2}-1}$), but the four parts of which $\mathcal T_{\sqrt{2}-1}(R_\alpha)$ and $\mathcal T_{\sqrt{2}-1}^2(R_\alpha)$ consist, show as coves. In Figure \ref{fig: omega.4051, added, rough} we see Figure \ref{fig: omega.4051, removed, rough} with the added regions; $R_{\alpha,1}$ is in grey:
\begin{align*}
A_\alpha&=[\alpha-1,\sqrt{2}-2]\times[0,\tfrac{g^2}{1+g^2})\cup [\alpha-1,\sqrt{2}-2]\times[\tfrac13,g^2),\\
\mathcal T_{\sqrt{2}-1}(A_\alpha)&=[\tfrac{1-2\alpha}{\alpha-1},\tfrac12\sqrt{2}-1]\times[\tfrac12,\tfrac{1+g^2}{2+g^2})\cup[\tfrac{1-2\alpha}{\alpha-1},\tfrac12\sqrt{2}-1]\times[\tfrac35,g)\,\, {\text{and}} \\
\mathcal T^2_{\sqrt{2}-1}(A_\alpha)&=[\tfrac{2-5\alpha}{2\alpha-1},\sqrt{2}-1]\times[\tfrac25,\tfrac{2+g^2}{5+2g^2})\cup[\tfrac{2-5\alpha}{2\alpha-1},\sqrt{2}-1]\times[\tfrac5{12},\tfrac{2+g}{5+2g}). \\
\end{align*}
Similar to the construction in Section \ref{The case alpha in (frac12,g]}, we are now left with $R_{\alpha,1}$ to remove, consisting of the protuberant parts of two small rectangles between the parts that we already removed.
Figure \ref{fig: omega.4051} shows a detailed example of $\Omega_{\alpha,1}$, $2/5<\alpha<\sqrt{2}-1$. This first stage shows how $\Omega_\alpha$ is not only characterised by more and more protuberant parts and coves as $\alpha$ decreases, but is also splitting in more and more disjoint regions. Nevertheless, the construction of $\Omega_\alpha$ in the current case is very similar to the construction of $\Omega_\alpha$ in Section \ref{The case alpha in (frac12,g]}: likewise, we construct the sequence $\Omega_{\alpha,k}$, $k=1,2,\dots$, and bring forth $\Omega_\alpha$ as being $\lim_{k \to \infty} \Omega_{\alpha,k}$. \\
\begin{figure}[!htb]
\minipage{0.48\textwidth}
$$
\begin{tikzpicture}[scale =5.5] 
 \draw[draw=white,fill=gray!20!white] 
 plot[smooth,samples=100,domain=.3:.414] (\x,0) --
 plot[smooth,samples=100,domain=.414:.3] (\x,0.382);
 \draw[draw=white,fill=gray!20!white] 
 plot[smooth,samples=100,domain=.3:.414] (\x,.5) --
 plot[smooth,samples=100,domain=.414:.3] (\x,0.618);
    \draw[draw=white,fill=gray!20!white] 
 plot[smooth,samples=100,domain=-.586:-.5] (\x,.3) --
 plot[smooth,samples=100,domain=-.5:-.586] (\x,.32);
\draw[draw=white,fill=gray!20!white] 
 plot[smooth,samples=100,domain=-.586:-.5] (\x,.33) --
 plot[smooth,samples=100,domain=-.5:-.586] (\x,.36);
     \draw[draw=white,fill=gray!20!white] 
 plot[smooth,samples=100,domain=-.293:-.27] (\x,.53) --
 plot[smooth,samples=100,domain=-.27:-.293] (\x,.54);
   \draw[draw=white,fill=gray!20!white] 
 plot[smooth,samples=100,domain=-.293:-.27] (\x,.55) --
 plot[smooth,samples=100,domain=-.27:-.293] (\x,.58);
  \draw (-.586,0) -- (.414,0);
 \draw [dashed] (.414,0) -- (.414,.382);
 \draw [dashed] (.3,0) -- (.3,.618);
  \draw [dashed] (.414,.5) -- (.414,.618);
  \draw [dashed] (-.293,.62) -- (.414,.62);
  \draw [dashed] (-.293,0) -- (-.293,.5);
     \draw (-.586,0) -- (-.586,.3);
       \draw (-.586,.3) -- (-.5,.3);
  \draw (-.5,.3) -- (-.5,.32);
  \draw (-.5,.32) -- (-.586,.32);
   \draw (-.586,.32) -- (-.586,.33);
       \draw (-.586,.33) -- (-.5,.33);
  \draw (-.5,.33) -- (-.5,.36);
  \draw (-.5,.36) -- (-.586,.36);
  \draw (-.586,.36) -- (-.586,.382);
   \draw [dashed] (0,0) -- (0,.618);
    \draw [dashed] (-.586,.382) -- (.414,.382);
    \draw (.-.293,.5) -- (.414,.5);
     \draw (-.293,.5) -- (-.293,.53);
      \draw (-.293,.53) -- (-.27,.53);
 \draw (-.27,.53) -- (-.27,.54);
 \draw (-.27,.54) -- (-.293,.54);
 \draw (-.293,.54) -- (-.293,.55);
 \draw (-.293,.55) -- (-.27,.55);
 \draw (-.27,.55) -- (-.27,.58);
  \draw (-.27,.58) -- (-.293,.58);
 \draw (-.293,.58) -- (-.293,.618);
\node at (.3,-.05) {$\alpha$};
         \node at (-.32,.5) {$\tfrac12$};
      \node at (-.293,-.05) {$\tfrac12\sqrt{2}-1$};
   \node at (-.62,.382){$g^2$};
   \node at (-.6,-.05){$\sqrt{2}-2$};
     \node at (-.32,.62){$g$};
       \node at (.43,-.05) {$\sqrt{2}-1$};
                 \node at (0,-.05) {0};
   \end{tikzpicture}
 $$
\caption[$\Omega_{\alpha,1}$, removed parts, rough sketch.$] {\label{fig: omega.4051, removed, rough}
$\Omega_{\alpha,1}$, removed parts, rough sketch.}
\endminipage\hfill
\minipage{0.48\textwidth}
$$
\begin{tikzpicture}[scale =5.5] 
 \draw[draw=white,fill=gray!20!white] 
 plot[smooth,samples=100,domain=.3:.414] (\x,.4) --
 plot[smooth,samples=100,domain=.414:.3] (\x,0.43);
 \draw[draw=white,fill=gray!20!white] 
 plot[smooth,samples=100,domain=.3:.414] (\x,.44) --
 plot[smooth,samples=100,domain=.414:.3] (\x,0.46);
    \draw[draw=white,fill=gray!20!white];
 \draw (-.7,0) -- (.3,0);
 \draw [dashed] (.414,0) -- (.414,.618);
 \draw [dashed] (.3,0) -- (.3,.618);
   \draw [dashed] (-.31,.62) -- (.3,.618);
  \draw [dashed] (-.293,0) -- (-.293,.618);
  \draw [dashed] (-.586,0) -- (-.586,.382);
  
  \draw (.25,.4) -- (.414,.4);
    \draw (.25,.4) -- (.25,.43);
\draw [dashed] (.25,.43) -- (.414,.43);
  \draw (.25,.44) -- (.414,.44);
    \draw (.25,.44) -- (.25,.46);
\draw [dashed] (.25,.46) -- (.414,.46);

     \draw (-.7,0) -- (-.7,.3);
       \draw (-.7,.3) -- (-.5,.3);
  \draw (-.5,.3) -- (-.5,.32);
  \draw (-.5,.32) -- (-.7,.32);
   \draw (-.7,.32) -- (-.7,.33);
       \draw (-.7,.33) -- (-.5,.33);
  \draw (-.5,.33) -- (-.5,.36);
  \draw (-.5,.36) -- (-.7,.36);
  \draw (-.7,.36) -- (-.7,.382);
   \draw [dashed] (0,0) -- (0,.618);
    \draw [dashed] (-.7,.382) -- (.3,.382);
    \draw (.-.31,.5) -- (.3,.5);
     \draw (-.31,.5) -- (-.31,.53);
      \draw (-.31,.53) -- (-.27,.53);
 \draw (-.27,.53) -- (-.27,.54);
 \draw (-.27,.54) -- (-.293,.54);
 \draw (-.293,.54) -- (-.293,.55);
 \draw (-.293,.55) -- (-.27,.55);
 \draw (-.27,.55) -- (-.27,.58);
  \draw (-.27,.58) -- (-.31,.58);
 \draw (-.31,.58) -- (-.31,.618);
\node at (.3,-.05) {$\alpha$};
         \node at (-.34,.5) {$\tfrac12$};
      \node at (-.293,-.05) {$\tfrac12\sqrt{2}-1$};
   \node at (-.74,.382){$g^2$};
   \node at (-.53,-.05){$\sqrt{2}-2$};
     \node at (-.34,.62){$g$};
       \node at (.43,-.05) {$\sqrt{2}-1$};
  \node at (-.72,-.05) {$\alpha-1$};     
                 \node at (0,-.05) {0};
   \end{tikzpicture}
 $$
\caption[$\Omega_{\alpha,1}$, removed and added parts, rough sketch.$] {\label{fig: omega.4051, added, rough}
$\Omega_{\alpha,1}$, removed and added parts, rough sketch.}
\endminipage
\end{figure}
We remark that fractions involving $g$ or $G$ can be represented in many ways. In this particular case, substituting $g$ for $v_n$ in $1/(3+v_n)$ (associated with the $v$-bounds of $\mathcal T_{\sqrt{2}-1}(R_\alpha)$) yields $1/(3+g)$, which is equal to $g^2/(1+g^2)$; this is the fraction we find when we substitute $g^2$ for $v_n$ in $v_n/(1+v_n)$, which is associated with the $v$-bounds of $A_\alpha$; see the lower corner on the left side of Figure \ref{fig: omega.4051}. 

To obtain $\Omega_\alpha$ as a `limit' of the sequence $\Omega_{\alpha,k}$, we need to infinitely remove and add regions. Similar to the case $1/2 < \alpha \leq g$, only new $v$-values need to be calculated, to which end we use the M\"obius transformations associated with $v_n^*$, $v_{n+1}^*$, $v_{n+2}^*$, $v_{n+1}$ and $v_{n+2}$. The matrices involved are immediately derived from the way $v_n^*$, $v_{n+1}^*$ et cetera are expressed in $v_n$. In the recurrent process of removing and adding regions, the matrix $\left (\begin{array}{cc}1 & 2 \\ 2 & 5 \end{array}\right)$ - the one belonging to $v_{n+2}^*$ - plays a central role, because it determines the domain of $v$-values in each new stage of constructing $\Omega_{\alpha,k}, k \geq 1$, similar to $\left (\begin{array}{cc}1 & 1 \\
1 & 2 \end{array} \right)$ in Section~\ref{The case alpha in (frac12,g]}. Using similar notations for new $v$-values also, we find 
$$
\left \{
\begin{aligned}
v_{n,k}^*&=\left (\begin{array}{cc}1 & 0 \\
1 & 1 \end{array} \right)\left (\begin{array}{cc}1 & 2 \\
2 & 5 \end{array} \right)^k(v_n), \,\, k\geq1, \\
v_{n+1,k}^*&=\left (\begin{array}{cc}1 & 1 \\
1 & 2 \end{array} \right)\left (\begin{array}{cc}1 & 2 \\
2 & 5 \end{array} \right)^k(v_n), \,\, k\geq1, \\
v_{n+2,k}^*&=\left (\begin{array}{cc} 1 & 2 \\
2 & 5 \end{array} \right)\left (\begin{array}{cc}1 & 2 \\
2 & 5 \end{array} \right)^k(v_n),\,\, k\geq1
\end{aligned}
\right.
$$
and
$$
\left \{
\begin{aligned}
v_{n+1,k}&=\left (\begin{array}{cc}0 & 1 \\
1 & 3 \end{array} \right)\left (\begin{array}{cc}1 & 2 \\
2 & 5 \end{array} \right)^k(v_n), \,\, k\geq1, \quad \text{and}\\
v_{n+2,k}&=\left (\begin{array}{cc} 1 & 3 \\
2 & 5 \end{array} \right)\left (\begin{array}{cc}1 & 2 \\
2 & 5 \end{array} \right)^k(v_n),\,\, k\geq1.
\end{aligned}
\right.
$$
It be remarked that graphic representations of $\Omega_{\alpha,k}, k \geq 2 $ would show merely very small changes compared to $\Omega_{\alpha,1}$. The six regions that are removed as well as the six regions that are added in each new stage are very small. \smallskip
We have, for example:
$$
A_{\alpha,1}=[\alpha-1,\sqrt{2}-2]\times[\tfrac{2}{7},\tfrac{2+g^2}{7+3g^2})\cup [\alpha-1,\sqrt{2}-2]\times[\tfrac{5}{17},\tfrac{2+g}{7+3g}), 
$$
while 
$\frac{2}{7}=0.285714\dots;\frac{2+g^2}{7+3g^2}=0.292412\dots;\frac{5}{17}=0.294117\dots$ and $\frac{2+g}{7+3g}=0.295686\dots$;
$$
\mathcal T_{\sqrt{2}-1}(A_{\alpha,1})=[\tfrac{1-2\alpha}{\alpha-1},\tfrac12\sqrt{2}-1]\times[\tfrac{7}{12},\tfrac{7+3g^2}{12+5g^2})\cup [\tfrac{1-2\alpha}{\alpha-1},\tfrac12\sqrt{2}-1]\times[\tfrac{17}{29},\tfrac{7+3g}{12+5g}), 
$$
while 
$\frac{7}{12}=0.583333\dots;\frac{7+3g^2}{12+5g^2}=0.585621\dots;\frac{17}{29}=0.586206\dots$ and $\frac{7+3g}{12+5g}=0.586746\dots$;
$$
\mathcal T_{\sqrt{2}-1}^2(A_{\alpha,1})=[\tfrac{2-5\alpha}{2\alpha-1},\sqrt{2}-1]\times[\tfrac{12}{29},\tfrac{12+5g^2}{29+12g^2})\cup[\tfrac{2-5\alpha}{2\alpha-1},\sqrt{2}-1]\times[\tfrac{29}{70},\tfrac{12+5g}{29+12g}), 
$$
while 
$\frac{12}{29}=0.413793\dots;\frac{12+5g^2}{29+12g^2}=0.414185\dots;\frac{29}{70}=0.414285\dots$ and $\frac{12+5g}{29+12g}=0.414378\dots$;
$$
\mathcal T_{\sqrt{2}-1}(R_{\alpha,1})=[\sqrt{2}-2,\tfrac{1-3\alpha}{\alpha}]\times(\tfrac{5+2g^2}{17+7g^2},\tfrac5{17}]\cup[\sqrt{2}-2,\tfrac{1-3\alpha}{\alpha}]\times(\tfrac{5+2g}{17+7g},\tfrac{12}{41}], 
$$
while 
$\frac5{17}=0.294117\dots;\frac{5+2g^2}{17+7g^2}=0.292975\dots;\frac{12}{41}=0.292682\dots$ and $\frac{5+2g}{17+7g}=0.292412\dots$ and
$$
\mathcal T_{\sqrt{2}-1}^2(R_{\alpha,1})=[\tfrac12\sqrt{2}-1,\tfrac{2-5\alpha}{3\alpha-1}]\times(\tfrac{17+7g^2}{29+12g^2},\tfrac{17}{29}]\cup[\tfrac12\sqrt{2}-1,\tfrac{1-3\alpha}{\alpha}]\times(\tfrac{17+7g}{29+12g},\tfrac{41}{70}], 
$$
while
$\frac{17}{29}=0.586206\dots;\frac{17+7g^2}{29+12g^2}=0.585814\dots;\frac{41}{70}=0.585714\dots$ and $\frac{17+7g}{29+12g}=0.585621\dots$.\\
In Figure \ref{fig: omega405,2zoom} the leftmost corner of $\Omega_{\alpha,2}$ is shown. Similar to Figure \ref{fig: omega.4051}, we have omitted the lower rectangle (where $v < \tfrac27$), in order not to waste too much blank space. Even more than Figure \ref{fig: omega405,2zoom}, Figure \ref{fig: omega405,2new} illustrates the rapid decrease of the surface of the rectangles, the union of which finally builds $\Omega_{\alpha}$. In particular,
$$
R_{\alpha,2}=[\alpha,\sqrt{2}-1]\times[\tfrac{12}{29},\tfrac{12+5g^2}{29+12g^2})\cup[\alpha,\sqrt{2}-1]\times[\tfrac{29}{70},\tfrac{12+5g}{29+12g}),
$$
the union of two very small rectangles between the two rectangles of $R_{\alpha,1}$. In the next step, we would find $R_{\alpha,3}$ consisting of two very small rectangles between the two rectangles of $R_{\alpha,2}$ - and so on.\\
Since $\left (\begin{array}{cc}1 & 2 \\
2 & 5 \end{array}\right)$= $\left (\begin{array}{cc}0 & 1 \\
1 & 2 \end{array}\right)^2$, we see that at the base of our current construction lies a matrix that is almost equal to the one in Section~\ref{The case alpha in (frac12,g]} that renders consecutive Fibonacci numbers. In this case, we have a matrix that renders the recurrent sequence $E_n$ given by $E_{n+1}=2E_{n}+E_{n-1}, n\geq 0$, with $E_{-1}=1$ and $E_0=0$. Now suppose $\lim_{k \to \infty}  \frac {E_{2k-1}v_n+E_{2k}}{E_{2k}v_n+E_{2k+1}}=L$. Then 
$$
L = \lim_{k \to \infty}  \frac 1{\tfrac{E_{2k}v_n+E_{2k+1}}{E_{2k-1}v_n+E_{2k}}}=\lim_{k \to \infty}  \frac 1{2+\tfrac{E_{2k-2}v_n+E_{2k-1}}{E_{2k-1}v_n+E_{2k}}}=\frac1{2+L},
$$
from which we derive 
$$
\lim_{k \to \infty} \left (\begin{array}{cc}1 & 1 \\
1 & 2 \end{array} \right)^k(v_n)=\sqrt{2}-1, \quad {\text{whence}}
$$
$$
\left \{
\begin{aligned}
\lim_{k \to \infty} v_{n,k}^*&=\tfrac {\sqrt{2}-1}{\sqrt{2}}=1-\tfrac12\sqrt{2}=0.292893\dots, \\
\lim_{k \to \infty} v_{n+1,k}^*&=\tfrac{\sqrt{2}}{\sqrt{2}+1}=2-\sqrt{2}=0.585786\dots, \\
\lim_{k \to \infty} v_{n+2,k}^*&=\tfrac{\sqrt{2}+1}{2\sqrt{2}+3}=\sqrt{2}-1=0.414213\dots, \\
\lim_{k \to \infty} v_{n+1,k}&=\tfrac1{\sqrt{2}+2}=1-\tfrac12\sqrt{2} \quad {\text{and}}\\
\lim_{k \to \infty} v_{n+2,k}&=\tfrac{\sqrt{2}+2}{2\sqrt{2}+3}=2-\sqrt{2}.
\end{aligned}
\right.
$$
\begin{figure}[!htb]
\minipage{0.5\textwidth}
$$
\begin{tikzpicture}[scale =60] 
\draw[draw=white,fill=gray!20!white] 
 plot[smooth,samples=100,domain=-.586:-.531] (\x,.067) --
 plot[smooth,samples=100,domain=-.531:-.586] (\x,0.07);
  \draw[draw=white,fill=gray!20!white] 
 plot[smooth,samples=100,domain=-.586:-.531] (\x,.073) --
 plot[smooth,samples=100,domain=-.531:-.586] (\x,0.084);
  \draw (-.595,0) -- (-.595,.067);
  \draw [dashed] (-.595,.067) -- (-.531,.067);
  \draw [->] (-.591,.077) -- (-.591,.083);
   \draw [->] (-.591,.074) -- (-.591,.068);
    \draw [->] (-.56,.062) -- (-.565,.066);
   \draw [->] (-.56,.062) -- (-.555,.072);
  \draw (-.531,.067) -- (-.531,.07);
  \draw (-.531,.07) -- (-.586,.07);
  \draw (-.586,.07) -- (-.586,.073);
  \draw [dashed]  (-.586,.073) -- (-.531,.073);
  \draw [dashed] (-.531,0) -- (-.524,0);
  \draw [dashed] (-.6,.07) -- (-.586,.07);
    \draw [dashed] (-.6,.084) -- (-.595,.084);
  \draw (-.531,.073) -- (-.531,.084);
  \draw (-.531,.084) -- (-.595,.084);
  \draw (-.595,.084) -- (-.595,.1);
  \draw [dashed] (-.595,.1) -- (-.531,.1);
  \draw [dashed] (-.586,.1) -- (-.586,0);
 
   \draw [dashed] (-.531,.07) -- (-.531,0);
      \draw [dashed] (-.531,.1) -- (-.531,.084);
    \node at (-.597,-.003) {\scriptsize{$\alpha-1$}};
      \node at (-.584,-.003) {\scriptsize{$\sqrt{2}-2$}};
          \node at (-.531,-.003) {\scriptsize{$\tfrac{1-3\alpha}{\alpha}$}};
             \node at (-.51,.1) {\scriptsize{$\tfrac{2+g}{7+3g}=0.2956\dots$}};
      \node at (-.613,.084) {\scriptsize{$\tfrac{5}{17}=0.2941\dots$}};
            \node at (-.51,.0745) {\scriptsize{$\tfrac{5+2g^2}{17+7g^2}=0.2929\dots$}};
          \node at (-.613,.07) {\scriptsize{$\tfrac{12}{41}=0.2926\dots$}};
        \node at (-.52,.066) {\scriptsize{$\tfrac{2+g^2}{7+3g^2}$}};
        \node at (-.522,.058) {\scriptsize{$=\tfrac{5+2g}{17+7g}$}};
\node at (-.518,.052) {\scriptsize{$=0.2924\dots$}};
     \node at (-.509,.000) {\scriptsize{$\tfrac27=0.2857\dots$}};
     \node at (-.591,.076) {\scriptsize{$A_{\alpha,1}$}};
   \node at (-.56,.06) {\scriptsize{$\mathcal T_\alpha(R_{\alpha,1})$}};
           \end{tikzpicture}
  $$
 \caption[The leftmost corner of$\Omega_{\alpha,2}, \alpha=0.405$.] {\label{fig: omega405,2zoom}
The leftmost corner of $\Omega_{\alpha,2}, \alpha=0.405$.}
\endminipage\hfill
\minipage{0.5\textwidth}
$$
\begin{tikzpicture}[scale =30] 
\draw[draw=white,fill=gray!20!white] 
 plot[smooth,samples=100,domain=.2025:.2071] (\x,.1379) --
 plot[smooth,samples=100,domain=.2071:.2025] (\x,.1418);
  \draw[draw=white,fill=gray!20!white] 
 plot[smooth,samples=100,domain=.2025:.2071] (\x,.1429) --
 plot[smooth,samples=100,domain=.2071:.2025] (\x,.1438);
\draw(.0658,0) -- (.2025,0);
\draw (.0658,0) -- (.0658,.1325);
\draw (.2025,0) -- (.2025,.1325);
\draw [dashed] (.0658,0) -- (.055,0);
\draw [dashed] (.0658,.1325) -- (.2025,.1325);
\draw [dashed] (.2071,-.01) -- (.2071,.1379);
\draw [dashed] (.2025,.1379) -- (.2025,.1438);
\draw (.0658,.1379) -- (.2071,.1379);
\draw [->] (.193,.151) -- (.202,.14);
\draw [->] (.193,.151) -- (.204,.144);
\draw (.0658,.1379) -- (.0658,.1418);
\draw (.2071,.1379) -- (.2071,.1418);
\draw [dashed] (.0658,.1418) -- (.202,.1418);
\draw(.0658,.1667) -- (.2025,.1667);
\draw (.0658,.1667) -- (.0658,.1982);
\draw (.2025,.1667) -- (.2025,.1982);
\draw [dashed] (.0658,.1982) -- (.2025,.1982);
\draw (.0658,.1429) -- (.2071,.1429);
\draw (.0658,.1429) -- (.0658,.1438);
\draw (.2071,.1429) -- (.2071,.1438);
\draw [dashed] (.0658,.1438) -- (.2025,.1438);
\node at (.05,-.03) {};
\node at (.2025,-.005) {\scriptsize{$\alpha$}};
\node at (.2071,-.015) {\scriptsize{$\sqrt{2}-1$}};
\node at (.07,-.008) {\scriptsize{$\tfrac{2-5\alpha}{2\alpha-1}$}};
\node at (.058,0) {\scriptsize{$\tfrac25$}};
\node at (.058,.134) {\scriptsize{$\tfrac{12}{29}$}};
\node at (.058,.145) {\scriptsize{$\tfrac{29}{70}$}};
\node at (.058,.1667) {\scriptsize{$\tfrac5{12}$}};
\node at (.225,.124) {\scriptsize{$\tfrac{2+g^2}{5+2g^2}$}};
\node at (.225,.138) {\scriptsize{$\tfrac{12+5g^2}{29+12g^2}$}};
\node at (.225,.152) {\scriptsize{$\tfrac{12+5g}{29+12g}$}};
\node at (.215,.1982) {\scriptsize{$\tfrac{2+g}{5+2g}$}};
\node at (.19,.155) {\scriptsize{$R_{\alpha,2}$}};
    \node at (.1325,.11) {\scriptsize{\it{The fractions on the left}}};
    \node at (.1325,.098) {\scriptsize{\it{are associated with}}};
    \node at (.1325,.085) {\scriptsize{\it{the lower boundaries}}};
    \node at (.1325,.072) {\scriptsize{\it{of the rectangles;}}};
    \node at (.1325,.058) {\scriptsize{\it{the fractions on the right}}};
    \node at (.1325,.046) {\scriptsize{\it{are associated with}}};
    \node at (.1325,.033) {\scriptsize{\it{the dashed upper boundaries.}}};
\end{tikzpicture}
  $$
\caption[The rectangles between the upper and lower block of $\Omega_{\alpha,2}, \alpha=0.405$.] {\label{fig: omega405,2new}
The rectangles between the upper and lower block of $\Omega_{\alpha,2}, \alpha=0.405$.}
\endminipage
\end{figure}

Omitting straightforward calculations, we conclude that 
\begin{align*}
\Omega_\alpha=&\,\,\Omega_{\sqrt{2}-1} \setminus \left (\bigcup_{k=0}^{\infty} R_{\alpha,k} \cup \mathcal T_{\sqrt{2}-1}(R_{\alpha,k})  \cup \mathcal T^2_{\sqrt{2}-1}(R_{\alpha,k}) \right ) \\
&\bigcup_{k=0}^{\infty} \left (A_{\alpha,k} \cup \mathcal T_{\sqrt{2}-1}(A_{\alpha,k}) \cup \mathcal T^2_{\sqrt{2}-1}(A_{\alpha,k})^*\right ),
\end{align*}
with
\begin{align*}
R_{\alpha,k} =&[\alpha,\sqrt{2}-1]\times[\tfrac{E_{2k}}{E_{2k+1}},\tfrac{E_{2k}+g^2E_{2k-1}}{E_{2k+1}+g^2E_{2k}}) \cup [\alpha,\sqrt{2}-1]\times[\tfrac{E_{2k+1}}{E_{2k+2}},\tfrac{E_{2k}+gE_{2k-1}}{E_{2k+1}+gE_{2k}});\\
\mathcal T_{\sqrt{2}-1}(R_{\alpha,k}) = &[\sqrt{2}-2,\tfrac{1-3\alpha}{\alpha}]\times(\tfrac{E_{2k+1}+g^2E_{2k}}{E_{2k+2}+E_{2k+1}+g^2(E_{2k+1}+E_{2k})},\tfrac{E_{2k+1}}{E_{2k+2}+E_{2k+1}}] \\
& \cup [\sqrt{2}-2,\tfrac{1-3\alpha}{\alpha}]\times(\tfrac{E_{2k+1}+gE_{2k}}{E_{2k+2}+E_{2k+1}+g(E_{2k+1}+E_{2k})},\tfrac{E_{2k+2}}{E_{2k+3}+E_{2k+2}}];\\ 
\mathcal T^2_{\sqrt{2}-1}(R_{\alpha,k}) = &[\tfrac12\sqrt{2}-1,\tfrac{2-5\alpha}{3\alpha-1}]\times(\tfrac{E_{2k+2}+E_{2k+1}+g^2(E_{2k+1}+E_{2k})}{E_{2k+3}+g^2E_{2k+2}},\tfrac{E_{2k+2}+E_{2k+1}}{E_{2k+3}}] \\
&\cup [\tfrac12\sqrt{2}-1,\tfrac{2-5\alpha}{3\alpha-1}]\times(\tfrac{E_{2k+2}+E_{2k+1}+g(E_{2k+1}+E_{2k})}{E_{2k+3}+gE_{2k+2}},\tfrac{E_{2k+3}+E_{2k+2}}{E_{2k+4}}];\\ 
A_{\alpha,k} =& [\alpha-1,\sqrt{2}-2]\times[\tfrac{E_{2k}}{E_{2k+1}+E_{2k}},\tfrac{E_{2k}+g^2E_{2k-1}}{E_{2k+1}+E_{2k}+g^2(E_{2k}+E_{2k-1})}) \\
&\cup [\alpha-1,\sqrt{2}-2]\times[\tfrac{E_{2k+1}}{E_{2k+2}+E_{2k+1}},\tfrac{E_{2k}+gE_{2k-1}}{E_{2k+1}+E_{2k}+g(E_{2k}+E_{2k-1})}) ;\\ 
\end{align*}
\begin{align*}
\mathcal T_{\sqrt{2}-1}(A_{\alpha,k}) =& [\tfrac{1-2\alpha}{\alpha-1},\tfrac12\sqrt{2}-1]\times[\tfrac{E_{2k+1}+E_{2k}}{E_{2k+2}},\tfrac{E_{2k+1}+E_{2k}+g^2(E_{2k}+E_{2k-1})}{E_{2k+2}+g^2E_{2k+1}}) \\
&\cup [\tfrac{1-2\alpha}{\alpha-1},\tfrac12\sqrt{2}-1]\times[\tfrac{E_{2k+2}+E_{2k+1}}{E_{2k+3}},\tfrac{E_{2k+1}+E_{2k}+g(E_{2k}+E_{2k-1})}{E_{2k+2}+gE_{2k+1}}) ;\\ 
\mathcal T^2_{\sqrt{2}-1}(R_{\alpha,k})^* = &[\tfrac{2-5\alpha}{2\alpha-1},\alpha)\times[\tfrac{E_{2k+2}}{E_{2k+3}},\tfrac{E_{2k+2}+g^2E_{2k+1}}{E_{2k+3}+g^2E_{2k+2}}) \cup [\tfrac{2-5\alpha}{2\alpha-1},\alpha)\times[\tfrac{E_{2k+3}}{E_{2k+4}},\tfrac{E_{2k+2}+gE_{2k+1}}{E_{2k+3}+gE_{2k+2}}),
\end{align*}
where $\mathcal T^2_{\sqrt{2}-1}(R_{\alpha,k})^*:= \mathcal T^2_{\sqrt{2}-1}(R_{\alpha,k})\setminus R_{\alpha,k+1}$, $R_{\alpha,0}:=R_\alpha$ and $A_{\alpha,0}:=A_\alpha$.\\
In order to find the effect of the current transformation on the sequence of convergents, we recall that we rewrote 
\begin{align*}
x&=[0;\varepsilon_1a_1,\dots,\varepsilon_na_n,3,-2,-a_{n+3},\dots] \quad {\text{as}}\\
x&=[0;\varepsilon_1a_1,\dots,\varepsilon_n(a_n+1),-2,-3,a_{n+3}-1,\dots]
\end{align*}
Since $p_k^*=p_k$ and $q_k^*=q_k$, $k<n$, we find, applying (\ref{pn en qn}) and this time omitting the straightforward calculations,
$$
\left \{
\begin{aligned}p_n^*&=p_n+p_{n-1};\\
p_{n+1}^*&=p_{n+1}-p_n;\\
p_{n+2}^*&=p_{n+2};\\
p_{n+3}^*&=p_{n+3}.\\
\end{aligned}
\right.
$$

Of course, we find similar relations for $q_n^*$ through $q_{n+3}^*$. We conclude that in the transformation from $\Omega_{\sqrt{2}-1}$ to $\Omega_\alpha$, $\alpha \in [2/5,\sqrt{2}-1)$, the convergent $p_n/q_n$ is replaced by the mediant $(p_n+p_{n-1})/(q_n+q_{n-1})$ and that $p_{n+1}/q_{n+1}$ is replaced by the mediant $(p_{n+1}-p_n)/(q_{n+1}-q_n)$.\\

\begin{figure}[!htb]
\minipage{0.48\textwidth}
$$
\begin{tikzpicture}[scale =7] 
\draw[draw=white,fill=gray!20!white] 
 plot[smooth,samples=100,domain=-.6:-.498] (\x,0) --
 plot[smooth,samples=100,domain=-.498:-.6] (\x,.0292);
 \draw[draw=white,fill=gray!20!white] 
 plot[smooth,samples=100,domain=-.6:-.498] (\x,.057) --
 plot[smooth,samples=100,domain=-.498:-.6] (\x,.077);
 \draw[draw=white,fill=gray!20!white] 
 plot[smooth,samples=100,domain=-.6:-.498] (\x,.082) --
 plot[smooth,samples=100,domain=-.498:-.6] (\x,.087);
  \draw[draw=white,fill=gray!20!white] 
 plot[smooth,samples=100,domain=-.6:-.498] (\x,.2) --
 plot[smooth,samples=100,domain=-.498:-.6] (\x,.346);
   \draw[draw=white,fill=gray!20!white] 
 plot[smooth,samples=100,domain=-.502:.4] (\x,0) --
 plot[smooth,samples=100,domain=.4:-.502] (\x,.346);
   \draw[draw=white,fill=gray!20!white] 
 plot[smooth,samples=100,domain=0:.4] (\x,.4) --
 plot[smooth,samples=100,domain=.4:0] (\x,.44); 
 \draw[draw=white,fill=gray!20!white] 
 plot[smooth,samples=100,domain=0:.4] (\x,.45) --
 plot[smooth,samples=100,domain=.4:0] (\x,.46);
 \draw[draw=white,fill=gray!20!white] 
 plot[smooth,samples=100,domain=-.333:.001] (\x,.7) --
 plot[smooth,samples=100,domain=0.001:-.333] (\x,.94);
  \draw[draw=white,fill=gray!20!white] 
 plot[smooth,samples=100,domain=-.333:.001] (\x,.95) --
 plot[smooth,samples=100,domain=0.001:-.333] (\x,.956);
 \draw[draw=white,fill=gray!20!white] 
 plot[smooth,samples=100,domain=-.333:.001] (\x,.957) --
 plot[smooth,samples=100,domain=0.001:-.333] (\x,.96);
  \draw[draw=white,fill=gray!20!white] 
 plot[smooth,samples=100,domain=-.333:.001] (\x,1) --
 plot[smooth,samples=100,domain=0.001:-.333] (\x,1.054);
  \draw[draw=white,fill=gray!20!white] 
 plot[smooth,samples=100,domain=-.001:.4] (\x,.7) --
 plot[smooth,samples=100,domain=.4:-.001] (\x,1.054);

  \draw (.4,0) -- (.4,.346);
   \draw [dashed] (.4,.346) -- (.4,.7);
 \draw (.4,.7) -- (.4,1.054);
  \draw [dashed] (-.586,.087) -- (-.5,.087);
    \draw (-.586,.199) -- (-.5,.199);
         \draw (-.5,.087) -- (-.5,.199);
    \draw (-.6,0) -- (-.6,.0292);
       \draw (-.6,.199) -- (-.6,.346);
     \draw (-.6,.199) -- (-.586,.199);
        \draw [dashed] (-.6,.0292) -- (-.5,.0292);
  \draw (0,.4) -- (.4,.4);
  \draw (0,.4) -- (0,.4399);
  \draw [dashed] (0,.4399) -- (.4,.4399); 
      \draw (0,.45) -- (.4,.45);
    \draw (0,.45) -- (0,.46);
    \draw [dashed] (-.333,.7) -- (-.333,.65);
 \draw [dashed] (0,.46) -- (.4,.46);
    \draw (-.6,.057) -- (-.5,.057);
    \draw (-.6,.057) -- (-.6,.077);
    \draw [dashed] (-.6,.077) -- (-.5,.077);
    \draw (-.5,.077) -- (-.5,.078);
       \draw (-.5,.079) -- (-.5,.082);
    \draw (-.6,.082) -- (-.5,.082);
    \draw (-.6,.082) -- (-.6,.087);
    \draw [dashed] (-.6,.087) -- (-.586,.087);
   \draw (-.5,.028) -- (-.5,.057);
    \draw (-.333,1) -- (0,1);
        \draw (-.333,1) -- (-.333,1.054);
    \draw (0,.96) -- (0,1);
        \draw [dashed] (-.333,.94) -- (0,.94);
  \draw [dashed] (-.333,1.054) -- (.4,1.054);
   \draw [dashed] (0,0) -- (0,.4);
   \draw [dashed] (0,.457) -- (0,.94);
      \draw [dashed] (-.6,.346) -- (.405,.346);
    \draw (.-.333,.7) -- (.4,.7);
    \draw (-.333,.7) -- (-.333,.94);
      \draw (0,.94) --(0,.95);
     \draw (0,.95) -- (-.333,.95);
     \draw (-.333,.95) -- (-.333,.956);
     \draw [dashed] (-.333,.956) -- (0,.956);
     \draw (0,.956) -- (0,.957);
     \draw (0,.957) -- (-.333,.957);
     \draw (-.333,.957) -- (-.333,.96);
     \draw [dashed] (-.333,.96) -- (0,.96);
     \node at (-.35,1) {\scriptsize{$\tfrac35$}};
      \node at (-.35,.7) {\scriptsize{$\tfrac12$}};
            \node at (-.48,.199) {\scriptsize{$\tfrac13$}};
          \node at (-.02,.4) {\scriptsize{$\tfrac25$}};
        \node at (-.615,.346) {\scriptsize{$g^2$}};
     \node at (-.5,-.03) {\scriptsize{$-\tfrac12$}};
   \node at (.4,-.03) {\scriptsize{$\tfrac25$}};
   \node at (-.6,-.03) {\scriptsize{$-\tfrac35$}};
       \node at (-.35,1.054) {\scriptsize{$g$}};
            \node at (0,-.03) {\scriptsize{0}};
            \node at (.05,.94) {\scriptsize{$\tfrac{3+g}{5+2g}$}};
      \node at (-.45,.087) {\scriptsize{$\tfrac1{3+g^2}$}};
 \node at (-.45,.0292) {\scriptsize{$\tfrac1{3+g}$}};
                 \node at (-.05,.469) {\scriptsize{$\tfrac{2+g}{5+2g}$}};
           \node at (-.33,.64) {\scriptsize{$-\tfrac13$}};
              \end{tikzpicture}
  $$
 \caption[$\Omega_{\tfrac25}$, 'rough sketch'.] {\label{fig: omega2/5}
$\Omega_{\tfrac25}$, 'rough sketch'.}
\endminipage\hfill
\minipage{0.48\textwidth}
$$
\begin{tikzpicture}[scale =7] 
\draw[draw=white,fill=gray!20!white] 
 plot[smooth,samples=100,domain=-.613:-.418] (\x,0) --
 plot[smooth,samples=100,domain=-.418:-.613] (\x,.0292);
 \draw[draw=white,fill=gray!20!white] 
 plot[smooth,samples=100,domain=-.613:-.418] (\x,.057) --
 plot[smooth,samples=100,domain=-.418:-.613] (\x,.077);
 \draw[draw=white,fill=gray!20!white] 
 plot[smooth,samples=100,domain=-.613:-.418] (\x,.082) --
 plot[smooth,samples=100,domain=-.418:-.613] (\x,.087);
  \draw[draw=white,fill=gray!20!white] 
 plot[smooth,samples=100,domain=-.613:-.418] (\x,.2) --
 plot[smooth,samples=100,domain=-.418:-.613] (\x,.346);
   \draw[draw=white,fill=gray!20!white] 
 plot[smooth,samples=100,domain=-.42:.387] (\x,0) --
 plot[smooth,samples=100,domain=.387:-.42] (\x,.346);
   \draw[draw=white,fill=gray!20!white] 
 plot[smooth,samples=100,domain=-.279:.387] (\x,.4) --
 plot[smooth,samples=100,domain=.387:-.279] (\x,.44); 
 \draw[draw=white,fill=gray!20!white] 
 plot[smooth,samples=100,domain=-.279:.387] (\x,.45) --
 plot[smooth,samples=100,domain=.387:-.279] (\x,.46);
 \draw[draw=white,fill=gray!20!white] 
 plot[smooth,samples=100,domain=-.368:.387] (\x,.7) --
 plot[smooth,samples=100,domain=.387:-.368] (\x,.94);
  \draw[draw=white,fill=gray!20!white] 
 plot[smooth,samples=100,domain=-.368:.387] (\x,.95) --
 plot[smooth,samples=100,domain=.387:-.368] (\x,.956);
 \draw[draw=white,fill=gray!20!white] 
 plot[smooth,samples=100,domain=-.368:.387] (\x,.957) --
 plot[smooth,samples=100,domain=.387:-.368] (\x,.96);
  \draw[draw=white,fill=gray!20!white] 
 plot[smooth,samples=100,domain=-.368:.387] (\x,1) --
 plot[smooth,samples=100,domain=.387:-.368] (\x,1.054);
   \draw [dashed] (.387,0) -- (.387,1.054);
             \draw (-.419,.087) -- (-.419,.2);
    \draw (-.613,0) -- (-.613,.0292);
       \draw (-.613,.2) -- (-.613,.346);
         \draw [dashed] (-.613,.0292) -- (-.419,.0292);
    \draw (-.279,.4) -- (.387,.4);
  \draw (-.279,.4) -- (-.279,.4399);
  \draw [dashed] (-.279,.4399) -- (.387,.4399); 
  \draw (-.419,.2) -- (-.613,.2);
      \draw (-.279,.45) -- (.387,.45);
      \draw (-.279,.45) -- (-.279,.46);
 \draw [dashed] (-.279,.46) -- (.387,.46);
     \draw (-.279,.441) -- (.387,.441);
      \draw (-.279,.441) -- (-.279,.4426);
 \draw [dashed] (-.279,.4426) -- (.387,.4426);
     \draw (-.279,.4429) -- (.387,.4429);
      \draw (-.279,.4429) -- (-.279,.4431);
 \draw [dashed] (-.279,.4431) -- (.387,.4431);
    \draw (-.613,.057) -- (-.419,.057);
    \draw (-.613,.057) -- (-.613,.077);
    \draw [dashed] (-.613,.077) -- (-.419,.077);
    \draw (-.419,.077) -- (-.419,.078);
         \draw (-.419,.079) -- (-.419,.082);
    \draw (-.613,.082) -- (-.419,.082);
  \draw (-.613,.082) -- (-.613,.087);
  \draw [dashed] (-.613,.087) -- (-.419,.087);
     \draw (-.419,.0292) -- (-.419,.057);
\draw [dashed] (-.419,0) -- (-.419,.0292);
    \draw (-.368,1) -- (.387,1);
        \draw (-.368,1) -- (-.368,1.054);
   \draw [dashed] (-.368,.66) -- (-.368,.7);
   \draw [dashed] (-.279,0) -- (-.279,.4);
   \draw [dashed] (-.368,1.054) -- (.387,1.054);
      \draw [dashed] (-.613,.346) -- (.387,.346);
    \draw (.-.368,.7) -- (.387,.7);
   \draw (-.368,.7) -- (-.368,.94);
   \draw [dashed] (-.368,.94) -- (.387,.94);
      \draw (.-.368,.95) -- (.387,.95);
   \draw (-.368,.95) -- (-.368,.957);
   \draw (-.368,.957) -- (.387,.957);
   \draw [dashed]  (.-.368,.956) -- (.387,.956);
   \draw (-.368,.956) -- (-.368,.96);
   \draw [dashed] (-.368,.96) -- (.387,.96);
   \node at (-.613,-.04) {\scriptsize{$\tfrac{\sqrt{10}-5}{3}$}};
   \node at (-.425,-.04) {\scriptsize{$\tfrac{\sqrt{10}-4}{2}$}};
    \node at (-.368,.63) {\scriptsize{$\tfrac{\sqrt{10}-5}{5}$}};
    \node at (-.27,-.04) {\scriptsize{$\tfrac{\sqrt{10}-4}{3}$}};
   \node at (.37,-.04) {\scriptsize{$\tfrac{\sqrt{10}-2}{3}$}};
 \node at (-.65,.0292) {\scriptsize{$\tfrac1{3+g}$}};
 \node at (-.635,.2) {\scriptsize{$\tfrac13$}};
 \node at (-.635,.346) {\scriptsize{$g^2$}};
    \node at (-.3,.4) {\scriptsize{$\tfrac25$}};
      \node at (-.33,.469) {\scriptsize{$\tfrac{2+g}{5+2g}$}};
        \node at (-.385,.7) {\scriptsize{$\tfrac12$}};
       \node at (-.385,1.054) {\scriptsize{$g$}};
        \end{tikzpicture}
$$
 \caption[$\Omega_{\tfrac{\sqrt{10}-2}3}$, 'rough sketch'.] {\label{fig: omega(sqrt10-2)/3}
$\Omega_{\tfrac{\sqrt{10}-2}3}$, 'rough sketch'.}
\endminipage
\end{figure}
We will now show how to construct $\Omega_{\alpha}$ from $\Omega_{2/5}$, with $\sqrt{10}-2)/3<\alpha<2/5$. Let $\alpha \in ((\sqrt{10}-2)/3,2/5)$, $x \in [-3/5,2/5)$ and $n \geq 0$ be the smallest integer for which $t_n \in [\alpha,2/5]$, i.e., for which 
$$
(t_n,v_n) \in \bigcup_{k=0}^{\infty} \left ([\alpha,\tfrac25]\times\left ([\tfrac{E_{2k}}{E_{2k+1}},\tfrac{E_{2k}+g^2E_{2k-1}}{E_{2k+1}+g^2E_{2k}}) \cup [\tfrac{E_{2k+1}}{E_{2k+2}},\tfrac{E_{2k}+gE_{2k-1}}{E_{2k+1}+gE_{2k}})\right )\right ).
$$
We already know that in this case 
$$
x=[0;\varepsilon_1a_1,\dots,\varepsilon_na_n,3,-2,a_{n+3},\dots],
$$
that we rewrite as
$$
x=[0;\varepsilon_1a_1,\dots,\varepsilon_n(a_n+1),-2,-3,-(a_{n+3}+1),\dots];
$$
it is not hard to find that $a_{n+3} \geq 3$. Remember that in Section \ref{The case alpha in (sqrt{2}-1,tfrac12)} we saw how much the case $\alpha \in (\sqrt{2}-1, 1/2]$ resembled the case $\alpha \in (1/2,g]$. Here we find a similar resemblance between $\alpha\in ((\sqrt{10}-2)/3, 2/5)$ and $\alpha\in [2/5,\sqrt{2}-1]$. Of course, in the current case the continued fraction map is $T_{2/5}$ instead of $T_{\sqrt{2}-1}$, but straightforward calculations show that exactly the same equations for calculating the $(t,v)$-pairs for the regions $R_\alpha$ through $\mathcal T^2(A_\alpha)$ hold, which proves to be very convenient for the construction of $\Omega_\alpha$, $\alpha \in ((\sqrt{10}-2)/3, 2/5)$. It appears that the transformation consists solely of expanding the regions already added and removed, in other words: both the protuberant parts as the coves are extended. Of special interest are the coves in the upper block: as $\alpha$ decreases, they scoop out the upper block until it splits in infinitely many rectangles at $t = (\sqrt{10}-2)/3$; see Figure \ref{fig: omega(sqrt10-2)/3}. Finally, extending our research from $[2/5,\sqrt{2}-1]$ to $((\sqrt{10}-2)/3,\sqrt{2}-1]$ makes no difference for the way convergents are replaced by mediants as described earlier in this section.\\
We have now proved the following theorem:
\begin{Theorem}
Let $\alpha \in ((\sqrt{10}-2)/3,\sqrt{2}-1]$ and let the sequence $E_n, \, n \geq -1$, be defined by $E_{n+1}:=2E_{n}+E_{n-1}, n\geq 0$, with $E_{-1}:=1$ and $E_0:=0$. Define
\begin{align*}
&V_{1,k}:=[\tfrac{E_{2k+2}}{E_{2k+3}},\tfrac{E_{2k+2}+g^2E_{2k+1}}{E_{2k+3}+g^2E_{2k+2}}) \cup [\tfrac{E_{2k+3}}{E_{2k+4}},\tfrac{E_{2k+2}+gE_{2k+1}}{E_{2k+3}+gE_{2k+2}});\\
&V_{2,k}:=[\tfrac{E_{2k+1}+E_{2k}}{E_{2k+2}},\tfrac{E_{2k+1}+E_{2k}+g^2(E_{2k}+E_{2k-1})}{E_{2k+2}+g^2E_{2k+1}})  \cup [\tfrac{E_{2k+2}+E_{2k+1}}{E_{2k+3}},\tfrac{E_{2k+1}+E_{2k}+g(E_{2k}+E_{2k-1})}{E_{2k+2}+gE_{2k+1}});\\
&V_{3,k}:=[\tfrac{E_{2k}}{E_{2k+1}+E_{2k}},\tfrac{E_{2k}+g^2E_{2k-1}}{E_{2k+1}+E_{2k}+g^2(E_{2k}+E_{2k-1})})\cup[\tfrac{E_{2k+1}}{E_{2k+2}+E_{2k+1}},\tfrac{E_{2k}+gE_{2k-1}}{E_{2k+1}+E_{2k}+g(E_{2k}+E_{2k-1})}).
\end{align*}
 Then 
\begin{align*}
\Omega_\alpha =&[\tfrac{1-3\alpha}{\alpha},\alpha) \times [0,g^2) \cup [\tfrac{2-5\alpha}{3\alpha-1},\alpha)\times [\tfrac12,g)\\
& \bigcup_{k=0}^{\infty} \left ( [\tfrac{2-5\alpha}{2\alpha-1},\alpha)\times V_{1,k} \cup [\tfrac{1-2\alpha}{\alpha-1},\tfrac{2-5\alpha}{3\alpha-1}) \times V_{2,k} \cup [\alpha-1,\tfrac{1-3\alpha}\alpha)\times V_{3,k} \right ).
\end{align*}
\end{Theorem}

\section{The long way down to $g^2$}
\label{The long way down to g^2}

In the previous section we rewrote 
$$
x=[0;\varepsilon_1a_1,\dots,\varepsilon_na_n,3,-2,\varepsilon_{n+3}a_{n+3},\dots]
$$
as
$$
x=[0;\varepsilon_1a_1,\dots,\varepsilon_n(a_n+1),-2,-3,-\varepsilon_{n+3}(a_{n+3}+\varepsilon_{n+3}),\dots].
$$
In terms of M\"obius transformations this is
\begin{equation}\label{moebius to g^2}
\begin{split}
&\left (\begin{array}{cc}0 & \varepsilon_n \\
1 & a_n \end{array} \right)\left (\begin{array}{cc}0 & 1 \\
1 & 3 \end{array} \right)\left (\begin{array}{cc}0 & -1 \\
1 & 2 \end{array} \right)\left (\begin{array}{cc}0 & \varepsilon_{n+3} \\
1 & a_{n+3} \end{array} \right)\\
=&\left (\begin{array}{cc}0 & \varepsilon_n \\
1 & a_n+1 \end{array} \right)\left (\begin{array}{cc}0 & -1 \\
1 & 2 \end{array} \right)\left (\begin{array}{cc}0 & -1 \\
1 & 3 \end{array} \right)\left (\begin{array}{cc}0 & -\varepsilon_{n+3} \\
1 & a_{n+3}+\varepsilon_{n+3} \end{array} \right)\\
=&\left (\begin{array}{cc}5\varepsilon_n & 5\varepsilon_na_{n+3}+3\varepsilon_n\varepsilon_{n+3} \\
5a_n+2 & 5a_na_{n+3}+3\varepsilon_{n+3}a_n\\
{}&+2a_{n+3}+\varepsilon_{n+3}  \end{array} \right);
\end{split}
\end{equation}
We can extend and then generalise this: just as 
\begin{equation}\label{matrix product 1}
\begin{split}
&\left (\begin{array}{cc}0 & \varepsilon_n \\
1 & a_n \end{array} \right)\left (\begin{array}{cc}0 & 1 \\
1 & 3 \end{array} \right)\underbrace{\left (\begin{array}{cc}0 & -1 \\
1 & 3 \end{array} \right)\dots\left (\begin{array}{cc}0 & -1 \\
1 & 3 \end{array} \right)}_\text{$k-1$ times} \left (\begin{array}{cc}0 & -1 \\
1 & 2 \end{array} \right)\left (\begin{array}{cc}0 & \varepsilon_{n+k+2} \\
1 & a_{n+k+2} \end{array} \right)\\
=&\left (\begin{array}{cc}0 & \varepsilon_n \\
1 & a_n+1 \end{array} \right)\left (\begin{array}{cc}0 & -1 \\
1 & 2 \end{array} \right)\underbrace{\left (\begin{array}{cc}0 & -1 \\
1 & 3 \end{array} \right)\dots\left (\begin{array}{cc}0 & -1 \\
1 & 3 \end{array} \right)}_\text{$k$ times} \left (\begin{array}{cc}0 & -\varepsilon_{n+k+2} \\
1 & a_{n+k+2}+\varepsilon_{n+k+2} \end{array} \right)\\
=&\left (\begin{array}{cc}F_{2k+5}\varepsilon_n & F_{2k+5}\varepsilon_na_{n+k+2}+F_{2k+4}\varepsilon_n\varepsilon_{n+k+2} \\
F_{2k+5}a_n+F_{2k+4} & F_{2k+5}a_na_{n+k+2}+F_{2k+4}\varepsilon_{n+k+2}a_n\\{} & +F_{2k+3}a_{n+k+2}+F_{2k+2}\varepsilon_{n+k+2}  \end{array} \right),
\end{split}
\end{equation}
which can be derived from (\ref{moebius to g^2}) by inserting $\left (\begin{array}{cc}0 & -1 \\
1 & 3 \end{array} \right)$-matrices at the right place, we have 
\begin{align*}
&[0;\varepsilon_1a_1,\dots,\varepsilon_na_n,3,\underbrace{-3,\dots,-3}_\text{$k-1$ times},-2,\varepsilon_{n+k+2}a_{n+k+2},\dots]\\
=& [0;\varepsilon_1a_1,\dots,\varepsilon_n(a_n+1),-2,\underbrace{-3,\dots,-3}_\text{$k$ times},-\varepsilon_{n+k+2}(a_{n+k+2}+\varepsilon_{n+k+2}),\dots].
\end{align*}
We observe that 
$$
\lim_{k \to \infty} [0;3,\underbrace{-3,\dots,-3}_\text{$k-1$ times}]=g^2.
$$
The good part of this is that it appears to be possible to apply compensated insertion beyond $(\sqrt{10}-2)/3$, by means of extending the $(3,-2)$-insertion by adding partial quotients $-3$ between $3$ and $-2$. The bad part is that this possibility vanishes once we reach $g^2$. Taking the approach we took so far, the thing to do is letting $\alpha$ decrease from $(\sqrt{10}-2)/3$ to $g^2$ via the sequence $[0;3,\underbrace{-3,\dots,-3}_\text{$k-1$ times},-2,\varepsilon_{k+2}a_{k+2},\dots]$, starting with the case $k=2$, which is $[0;3,-3,-2,\varepsilon_4a_4,\dots]$.\smallskip

In the previous sections we have seen how our approach of singularisation and insertion gradually involved more details to process. We will show how things rapidly become even more complex, making it impossible to continue this approach much longer. So far, we have been able to construct $\Omega_\alpha$ for all $\alpha \in [(\sqrt{10}-2)/3,1)$. The next set of numbers to tackle would seem to consist of numbers $x=[0;\varepsilon_1a_1,\dots,\varepsilon_na_n,3,-3,-2,\varepsilon_{n+4}a_{n+4},\dots]$ and $t_n=[0;3,-3,-2,\varepsilon_{n+4}a_{n+4},\dots]$, the largest of which is $(\sqrt{10}-2)/3=[0;3,\overline{-3,-2,-3,-4}]$. Let us first determine the interval $(\alpha', \alpha)$, $\alpha' <\alpha \leq (\sqrt{10}-2)/3$, such that $t_n \in (\alpha', \alpha)$ implies $t_n=[0;3,-3,-2,\dots]$. We remark that 
$$
f_2(x):=\cfrac1{3-\cfrac1{3 - \frac1{x}}},
$$
where the index $2$ indicates the number of successive $3$s, is a decreasing function of $x$ on $\R \setminus \{0\}$, so if $\alpha' \in [0,(\sqrt{10}-2)/3)$ and $\alpha \in (\alpha',(\sqrt{10}-2)/3]$, then
\begin{equation}\label{boundary between 3,-3,-2 and 3,-3,-3}
2 = \left \lfloor f((\sqrt{10}-2)/3)+1-(\sqrt{10}-2)/3\right \rfloor \leq a_3(\alpha'_\alpha) =\left \lfloor f(\alpha')+1-\alpha \right \rfloor
\leq \left \lfloor f(\alpha')+1-\alpha' \right \rfloor, 
\end{equation}
from which we derive $a_3(\alpha'_\alpha)=2$ if and only if $f(\alpha')+1-\alpha'<3$. Given that $\alpha' \in [0,(\sqrt{10}-2)/3)$), this is the case if and only if $(\sqrt{65}-5)/8 <\alpha'<(\sqrt{10}-2)/3$; indeed 
$$
(\sqrt{65}-5)/8_{(\sqrt{65}-5)/8}=[0;3,\overline{-3,-3,-2,-3,-3,-4}].
$$
So the next step of our investigation would be the case $\alpha \in \Delta_\alpha(3,-3,-2)=((\sqrt{65}-5)/8,(\sqrt{10}-2)/3]$. Proceeding in a similar way, a sequence of successive intervals to be investigated would come into view, the boundaries of which are the positive roots $R_k$ of 
$$
f_k(x):=1/\underbrace{(3-1/(3-\dots1/(3}_{k \,\, {\text{times a}} \,\, 3}-1/x)\dots))-x-2.
$$
Straightforward calculation shows that 
$$
R_k=\frac{\sqrt{F_{2k+1}F_{2k+3}}-F_{2k+1}}{F_{2k+2}}
$$
and
$$
\lim_{k \to \infty} R_k=\lim_{k \to \infty} \left (\sqrt{\frac{F_{2k+1}F_{2k+3}}{F_{2k+2}^2}}-\frac {F_{2k+1}}{F_{2k+2}}\right )=\sqrt{g\cdot G}-g=1-g=g^2,
$$
as was to be expected.
For our investigations we would also need that
$$
{R_k}_{R_k}=[0;3,\overline{\underbrace{-3,\dots,-3}_{k \,\,{\text{times}}},-2,\underbrace{-3,\dots,-3}_{k \,\,{\text{times}}},-4}].
$$
Although this perspective would already promise a lot of intricacies, when we turn back to the case $\alpha \in ((\sqrt{65}-5)/8,(\sqrt{10}-2)/3]$, we are confronted with a problem we had not encountered yet: if we apply compensated insertion to numbers equal to or only slightly smaller than $(\sqrt{10}-2)/3$, we get
$$
\alpha_\alpha=[0;3,-3,-2,-3,-a,\dots] = [1;-2,-3,-3,2,-a,\dots],\,\,\,a\geq4.
$$
Of course, the occurrence of the partial quotient $2$ is not allowed, which calls for one more compensated insertion, so as to get
$$
\alpha_\alpha=[0;3,-3,-2,-3,-a,\dots] = [1;-2,-3,-3,2,-a,\dots]= [1;-2,-3,-4,-2,a-1,\dots], \,\,\,a\geq4.
$$
To determine the interval where these numbers occur, we have to solve 
$$
\cfrac1{\displaystyle 2 - \cfrac 1{\displaystyle 3 - \cfrac 1{\displaystyle 3 - \cfrac 1{\displaystyle \alpha }}}}+1-\alpha=4.
$$
The positive solution of this equation is 
$$
(5\sqrt{13}-13)/13_{(5\sqrt{13}-13)/13}=[0;3,\overline{-3,-2,-4,-2,-3,-4,-2,-4}]=0.386750\ldots.
$$
So the next case to investigate is actually $\alpha \in \Delta_\alpha(3,-3,-2,-3)=((5\sqrt{13}-13)/13,(\sqrt{10}-2)/3]$. More precisely, let $x \in [(\sqrt{10}-5)/3,(\sqrt{10}-2)/3)$ and $n \geq 0$ be the smallest integer for which $t_n \in [\alpha,(\sqrt{10}-2)/3]$, i.e., for which 
$$
(t_n,v_n) \in \bigcup_{k=0}^{\infty} \left \{[\alpha,(\sqrt{10}-2)/3]\times\left ([\tfrac{E_{2k}}{E_{2k+1}},\tfrac{E_{2k}+g^2E_{2k-1}}{E_{2k+1}+g^2E_{2k}}) \cup [\tfrac{E_{2k+1}}{E_{2k+2}},\tfrac{E_{2k}+gE_{2k-1}}{E_{2k+1}+gE_{2k}})\right )\right \}.
$$
Above we saw that in this case 
\begin{equation}\label{last rewriting}
\begin{aligned}
x&=[0;\varepsilon_1a_1,\dots,\varepsilon_na_n,3,-3,-2,-3,-a_{n+5},\dots] \,\,\, ,a_{n+5}\geq4, \,\,\, {\text{that we rewrite as}}\\
x&=[0;\varepsilon_1a_1,\dots,\varepsilon_n(a_n+1),-2,-3,-4,-2,a_{n+5}-1,\dots] \,\,\, ,a_{n+5}\geq4.
\end{aligned}
\end{equation}
We now have $(t_n^*,v_n^*)=(t_n-1,\tfrac{v_n}{v_n+1})$ and
\begin{align*}
&(t_{n+1},v_{n+1})=(\tfrac{1-3t_n}{t_n},\tfrac1{v_n+3}) & (t_{n+1}^*,v_{n+1}^*)&=(\tfrac{1-2t_n}{t_n-1},\tfrac{v_n+1}{v_n+2})\\
&(t_{n+2},v_{n+2})=(\tfrac{3-8t_n}{3t_n-1},\tfrac{v_n+3}{3v_n+8}) & (t_{n+2}^*,v_{n+2}^*)&=(\tfrac{2-5t_n}{2t_n-1},\tfrac{v_n+2}{2v_n+5})\\
&(t_{n+3},v_{n+3})=(\tfrac{5-13t_n}{8t_n-3},\tfrac{3v_n+8}{5v_n+13}) & (t_{n+3}^*,v_{n+3}^*)&=(\tfrac{7-18t_n}{5t_n-2},\tfrac{2v_n+5}{7v_n+18})\\
&(t_{n+4},v_{n+4})=(\tfrac{12-31t_n}{13t_n-5},\tfrac{5v_n+13}{12v_n+31}) & (t_{n+4}^*,v_{n+4}^*)&=(\tfrac{12-31t_n}{18t_n-7},\tfrac{7v_n+18}{12v_n+31}),\\
& (t_{n+5},v_{n+5})=(t_{n+5}^*,v_{n+5}^*) &{}
\end{align*}
Although it would still be possible to explicitly give $\Omega_\alpha$ for $\alpha\in ((5\sqrt{13}-13)/13,(\sqrt{10}-2)/3]$, it may be clear that it becomes quite unmanageable. Determining the effect of the current transformation on the sequence of convergents is still easy, though. Since $p_k^*=p_k$ and $q_k^*=q_k$, $k<n$, we find, applying (\ref{pn en qn}) and referring to (\ref{last rewriting}) (while omitting the straightforward calculations),
$$
\left \{
\begin{aligned}
p_n^*&=p_n+p_{n-1};\\
p_{n+1}^*&=p_{n+1}-p_n;\\
p_{n+2}^*&=p_{n+2}-p_{n+1};\\
p_{n+3}^*&=p_{n+4}-p_{n+3};\\
p_{n+4}^*&=p_{n+4};\\
p_{n+5}^*&=p_{n+5}.\\
\end{aligned}
\right.
$$
Of course, we find similar relations for $q_n^*$ through $q_{n+5}^*$. We conclude that in the transformation from $\Omega_{(\sqrt{10}-2)/3}$ to $\Omega_\alpha$, $\alpha \in [(5\sqrt{13}-13)/13,(\sqrt{10}-2)/3)$, the convergent $p_n/q_n$ is replaced by the mediant $(p_n+p_{n-1})/(q_n+q_{n-1})$, that $p_{n+1}/q_{n+1}$ is replaced by the mediant $(p_{n+1}-p_n)/(q_{n+1}-q_n)$, $p_{n+2}/q_{n+2}$ is replaced by the mediant $(p_{n+2}-p_{n+1})/(q_{n+2}-q_{n+1})$ and that $p_{n+3}/q_{n+3}$ is replaced by the mediant $(p_{n+4}-p_{n+3})/(q_{n+4}-q_{n+3})$.

A substantial part of our construction of $\Omega_\alpha$, $(\sqrt{10}-2)/3 \leq \alpha <1$, consisted in infinitely repeating the removal and addition of rectangles in $[(\sqrt{10}-5)/3,1]\times[0,1]$, in the cases that part of the added regions overlapped part of the regions to be removed. In Section \ref{The case alpha in (frac12,g]}, for instance, we worked at the case $x=[0;\varepsilon_1a_1,\dots,\varepsilon_na_n,2,-a_{n+2},\dots], a_{n+2} \geq 3,
$ and rewrote it as
\begin{equation}\label{singularisation between 1/2 and g re}
x=[0;\varepsilon_1a_1,\dots,\varepsilon_n(a_n+1),-2,(a_{n+2}-1),\dots].
\end{equation}
The overlapping stems from the possibility that (\ref{singularisation between 1/2 and g re}) can be written as 
\begin{equation}\label{singularisation between 1/2 and g special case}
x=[0;\varepsilon_1a_1,\dots,\varepsilon_n(a_n+1),-2,2,\varepsilon_{n+3}a_{n+3},\dots],\,\, a_{n+3}\geq3,
\end{equation}
i.e. when $a_{n+2}=3$. 
One could ask whether the work of infinite removal and addition could have possibly been avoided by immediately performing one more compensated insertion and so rewriting (\ref{singularisation between 1/2 and g special case}) as
$$
x=[0;\varepsilon_1a_1,\dots,\varepsilon_n(a_n+1),-3,-2,-\varepsilon_{n+3}(a_{n+3}+\varepsilon_{n+3}),\dots].
$$
This, however, confronts us with some intricacies, among which $t_{n+2}^*\neq t_{n+2}$, that appear to at least cancel the benefits of this approach. Earlier in this section we showed that for $\alpha \leq (\sqrt{10}-2)/3$, this repeated compensated insertion is unavoidable, as is a further complication of the construction of $\Omega_\alpha$. One obvious problem is the growth of the number of regions to be removed and added in each step of construction, due to the lengthening of the continued fraction expansions of the numbers involved on the way down to $g^2=0.381966\dots$, the first of which are (in decreasing order):
\begin{align*}
(5\sqrt{13}-13)/13_{(5\sqrt{13}-13)/13}&=[0;3,\overline{-3,-2,-4,-2,-3,-4,-2,-4}]\\
&=0.3867504\ldots,\\
5/13_{5/13}&=[0;3,-3,-2]\\
&=0.3846153\ldots,\\
(\sqrt{65}-5)/8_{(\sqrt{65}-5)/8}&=[0;3,\overline{-3,-3,-2,-3,-3,-4}]\\
&=0.3827822\ldots,\\
\end{align*}
\begin{align*}
(\sqrt{14401}-89)/81_{(\sqrt{14401}-89)/81}&=[0;3,-3,-3,-2,-3,-4,-2,-3,-3,-4,-3,-2,\dots]\\
&=0.3827674\ldots,\\
(\sqrt{2210}-34)/34_{(\sqrt{2210}-34)/34}&=[0;3,\overline{-3,-3,-2,-4,-2,-3,-3,-4,-2,-4}]\\
&=0.3826657\ldots,\\
13/34_{13/34}&=[0;3,-3,-3,-2]\\
&=0.3823529\ldots,\\
(\sqrt{442}-13)/21_{(\sqrt{442}-13)/21}&=[0;3,\overline{-3,-3,-3,-2,-3,-3,-3,-4}]\\
&=0.3820855\ldots,\\
(\sqrt{670762}-610)/547_{(\sqrt{670762}-610)/547}&=[0;3,-3,-3,-3,-2,-3,-3,-4,-2,-3,-3,-3\dots]\\
&=0.3820852\dots
\end{align*}
The red thread in this list is the first appearance of the partial quotient 2 with minus sign in the continued fraction expansion $\alpha_\alpha$ of the related decreasing numbers $\alpha$, always precedented by the sequence $3,\underbrace{-3,\dots,-3}_{k \,\,{\text{times}}},\,\, k > 0$. With $k$ increasing, the number of partial quotients $3$ with minus sign following this $2$ before the first appearance of the partial quotient $4$ with minus sign may be from $0$ to $k$, each case calling for an $\alpha$-fundamental interval to be investigated separately. On our way down from $(\sqrt{65}-5)/8 = [0;3,\overline{-3,-3,-2,-3,-3,-4}]$ to $(\sqrt{442}-13)/21_{(\sqrt{442}-13)/21}=[0;3,\overline{-3,-3,-3,-2,-3,-3,-3,-4}]$, for instance, we have to distinguish between the cases associated with the sequences starting with $3,-3,-3,-2,-3,-3,-4$, then $3,-3,-3,-2,-3,-4$, then $3,-3,-3,-2,-4$ and $3,-3,-3,-2$ before arriving at $3,-3,-3,-3,\dots$.\\

In sections \ref{The case alpha in (g,1]}, \ref{The case alpha in (sqrt{2}-1,tfrac12)} and \ref{The case alpha in (tfrac{sqrt{10}-2}3,sqrt{2}-1]} we saw how (compensated) insertion in most cases involves loss of convergents, although for $\alpha \leq 1/2$ each convergent is replaced by a mediant of the form $(p_k \pm p_{k-1})/(q_k \pm q_{k-1})$ for some $k \in \N$. As $\alpha$ decreases beyond $(5\sqrt{13}-13)/13$, this will remain the case, which we can best illustrate by again using M\"obius transformations. If $p_{n-1}$ and $p_n$ are two consecutive denominators of a convergent of some number $x=[0;\varepsilon_1a_1,\dots,\varepsilon_{n-1}a_{n-1},\varepsilon_na_n,\dots]$, then
$$
\left (\begin{array}{c} p_{n+1} \\
p_n \end{array} \right)=\left (\begin{array}{cc} a_{n+1} & \varepsilon_{n+1}  \\
1 & 0 \end{array} \right) \left (\begin{array}{c} p_n \\
p_{n-1} \end{array} \right).
$$
The denominators $p_{n+k}^*,\,\,k \in \{0,\dots,5\}$ that we calculated in the previous section, can easily be derived by comparing the partial products of
\begin{align*}
\left (\begin{array}{c} p_{n+5}^* \\
p_{n+4}^* \end{array} \right)= 
& \left (\begin{array}{cc} a_{n+5}-1 & 1  \\
1 & 0 \end{array} \right) 
 \left (\begin{array}{cc} 2 & -1  \\
1 & 0 \end{array} \right) 
 \left (\begin{array}{cc} 4 & -1  \\
1 & 0 \end{array} \right) 
 \left (\begin{array}{cc} 3 & -1 \\
1 & 0 \end{array} \right) \\
 & \cdot \left (\begin{array}{cc} 2 & -1 \\
1 & 0 \end{array} \right) 
 \left (\begin{array}{cc} a_n+1 & \varepsilon_n \\
1 & 0 \end{array} \right) 
\left (\begin{array}{c} p_{n-1} \\
p_{n-2} \end{array} \right)
\end{align*}
with those of
\begin{align*}
\left (\begin{array}{c} p_{n+5} \\
p_{n+4} \end{array} \right)=  & \left (\begin{array}{cc} a_{n+5} & 1  \\
1 & 0 \end{array} \right) 
 \left (\begin{array}{cc} 3 & -1  \\
1 & 0 \end{array} \right) 
 \left (\begin{array}{cc} 2 & -1  \\
1 & 0 \end{array} \right) 
 \left (\begin{array}{cc} 3 & -1 \\
1 & 0 \end{array} \right) \\
& \cdot  \left (\begin{array}{cc} 3 & 1 \\
1 & 0 \end{array} \right) 
 \left (\begin{array}{cc} a_n & \varepsilon_n \\
1 & 0 \end{array} \right) 
\left (\begin{array}{c} p_{n-1} \\
p_{n-2} \end{array} \right),
\end{align*}
associated with the sequence $\varepsilon_n(a_n+1),-2,-3,-4,-2,a_{n+5}-1$ as derived from the sequence $\varepsilon_na_n,3,-3,-2,-3,-a_{n+5}$ by applying compensated insertion.\\
On the interval $((\sqrt{14401}-89)/81,(\sqrt{65}-5)/8]$ we could similarly compare the partial products of
\begin{align*}
\left (\begin{array}{c} p_{n+7}^* \\
p_{n+6}^* \end{array} \right)= 
& \left (\begin{array}{cc} a_{n+7}-1 & 1  \\
1 & 0 \end{array} \right) 
 \left (\begin{array}{cc} 2 & -1  \\
1 & 0 \end{array} \right) 
 \left (\begin{array}{cc} 3 & -1  \\
1 & 0 \end{array} \right) 
 \left (\begin{array}{cc} 4 & -1 \\
 1 & 0 \end{array} \right) 
 \left (\begin{array}{cc} 3 & -1 \\
1 & 0 \end{array} \right) \\
& \cdot  \left (\begin{array}{cc} 3 & -1 \\
1 & 0 \end{array} \right) 
\left (\begin{array}{cc} 2 & -1 \\
1 & 0 \end{array} \right) 
 \left (\begin{array}{cc} a_n+1 & \varepsilon_n \\
1 & 0 \end{array} \right) 
\left (\begin{array}{c} p_{n-1} \\
p_{n-2} \end{array} \right)
\end{align*}
with those of
\begin{align*}
\left (\begin{array}{c} p_{n+7} \\
p_{n+6} \end{array} \right)=  & \left (\begin{array}{cc} a_{n+7} & 1  \\
1 & 0 \end{array} \right) 
 \left (\begin{array}{cc} 3 & -1  \\
1 & 0 \end{array} \right) 
 \left (\begin{array}{cc} 3 & -1  \\
1 & 0 \end{array} \right) 
 \left (\begin{array}{cc} 2 & -1 \\
 1 & 0 \end{array} \right) 
 \left (\begin{array}{cc} 3 & -1 \\
1 & 0 \end{array} \right) \\
 & \cdot \left (\begin{array}{cc} 3 & -1 \\
1 & 0 \end{array} \right) 
\left (\begin{array}{cc} 3 & 1 \\
1 & 0 \end{array} \right) 
 \left (\begin{array}{cc} a_n & \varepsilon_n \\
1 & 0 \end{array} \right) 
\left (\begin{array}{c} p_{n-1} \\
p_{n-2} \end{array} \right),
\end{align*}
associated with the sequence $\varepsilon_n(a_n+1),-2,-3,-3,-4,-3,-2,a_{n+7}-1$ as derived from the sequence $\varepsilon_na_n,3,-3,-3,-2,-3,-3,-a_{n+7}$ by applying compensated insertion. It is not hard to see that as $\alpha$ decreases to $g^2$, it is just the number of lost convergents that increases, but not the way in which they are replaced by mediants as described above; the difference is merely the number of matrices $\left (\begin{array}{cc} 3 & -1 \\
1 & 0 \end{array} \right)$ in both products, similar to the difference made by the number of matrices $\left (\begin{array}{cc} 0 & -1 \\
1 & 3 \end{array} \right)$ in (\ref{matrix product 1}).

\section{The ergodic systems $(\Omega_\alpha,\mathcal B, \mu_\alpha, \mathcal T_\alpha)$}
\label{ergodic systems}
Starting with $\Omega_1$, we constructed domains $\Omega_{\alpha'}$ by removing sets of points from a given $\Omega_\alpha$, with $\alpha' < \alpha$, and adding other sets to it. Each stage of construction -- corresponding with successive sections of this paper -- had its own characteristic set of singularisations, such as simply singularising a partial quotient $1$ for $\alpha' \in (g,1]$. The first step of constructing consisted of fixing a subset $R_{\alpha'}$ of $\Omega_\alpha$, consisting of all points $(t,v) \in \Omega_\alpha$ for which $\alpha'  \leq t < \alpha$.\\
Now let the collection of subsets of $\Omega_\alpha$ be denoted by $\mathcal B$ and $\mu_\alpha$ be defined as the probability measure with density 
$$
\frac1{N_\alpha}\cdot\frac1{(1+tv)^2}
$$
on $(\Omega_\alpha,\mathcal B)$, where $N_\alpha$ is a normalising constant. In \cite{Na}, Nakada showed that for $\alpha \in [1/2,1]$ the dynamical systems $(\Omega_\alpha,\mathcal B,\mu_\alpha,\mathcal T_\alpha)$ are ergodic (he actually obtained stronger mixing properties); see also \cite{K}. Here we will use this result for $\alpha=1$ only, which was the starting point in \cite{K} as well. Note that for $g^2 <\alpha<g$ our construction actually is a bijection between $\Omega_{\alpha'}$ and $\Omega_\alpha$, $\alpha'<\alpha$. Clearly this is not the case for $g <\alpha'<\alpha \leq 1$, where one region is added to $\Omega_\alpha$ while two are removed from it when constructing $\Omega_{\alpha'}$.\smallskip

 We will now show that $\mu_\alpha$ is an invariant measure on $\Omega_\alpha$. 
Let $x_\alpha=[0;\varepsilon_1a_1,\dots,\varepsilon_na_n,\varepsilon_{n+1}a_{n+1},\dots]$ and $D:= [t_1,t_2] \times [v_1,v_2] \subset R_{\alpha'}$. Then
$$
m(D):= \iint_D \frac1{(1+tv)^2} \,dt\,dv =\log \frac{(1+t_2v_2)(1+t_1v_1)}{(1+t_2v_1)(1+t_1v_2)}.
$$
Let $D_k:=\mathcal T^k_\alpha(D)$ and $M_n$ as the M\"obius transformation associated with the matrix $\left (\begin{array}{cc} 0 & \varepsilon_n\\
1 & a_n \end{array}\right )$. We define $\widetilde{M}_k:=M_{n+1} \cdots M_{n+k}$. From (\ref{relation v_{n+k} and t_{n+k}}) we derive that 
$$
D_k=\left [\widetilde{M}_k^{-1}(t_1),\widetilde{M}_k^{-1}(t_2)\right ] \times \left [\widetilde{M}_k^T(v_1),\widetilde{M}_k^T(v_2)\right ],
$$
from which it follows that
\begin{align*}
m(D_k)&=\iint_{D_k} \frac1{(1+tv)^2} \,dt\,dv\\
&=\log \frac{(1+\widetilde{M}_k^{-1}(t_2)\widetilde{M}_k^T(v_2))(1+\widetilde{M}_k^{-1}(t_1)\widetilde{M}_k^T(v_1))}{(1+\widetilde{M}_k^{-1}(t_2)\widetilde{M}_k^T(v_1))(1+\widetilde{M}_k^{-1}(t_1)\widetilde{M}_k^T(v_2))}\\
&=\log \frac{|\widetilde{M}_k|(t_2v_2+1)|\widetilde{M}_k|(t_1v_1+1)}{|\widetilde{M}_k|(t_2v_1+1)|\widetilde{M}_k|(t_1v_2+1)}\\
&=\log \frac{(1+t_2v_2)(1+t_1v_1)}{(1+t_2v_1)(1+t_1v_2)}=m(D). \qquad (*)
\end{align*}
Rewriting 
\begin{align*}
x_\alpha&=[0;\varepsilon_1a_1,\varepsilon_2a_2,\dots,\varepsilon_na_n,a_{n+1},\varepsilon_{n+2}a_{n+2},\dots] \,\,\,{\text{as}}\\
x_\alpha&=[0;\varepsilon_1a_1,\varepsilon_2a_2,\dots,\varepsilon_n(a_n+1),-a_{n+1}^*,\varepsilon_{n+2}^*a_{n+2}^*,\dots],
\end{align*}
as we have done throughout this paper, we define $D^*:= [t_1-1,t_2-1] \times [v_1/(1+v_1),v_2/(1+v_2)] \in A_{\alpha'}$. With a straightforward calculation one finds that $m(D^*)=m(D)$. Similarly to what we have written above, we define $D^*_k:=\mathcal T^k_\alpha(D^*)$ and let $M^*_n$ be the M\"obius transformation associated with the matrix $\left (\begin{array}{cc} 0 & \varepsilon^*_n\\
1 & a^*_n \end{array}\right )$. We define $\widetilde{M}^*_k:=\left (\begin{array}{cc}1 & 1\\
0 & 1 \end{array} \right)M^*_{n+1} \cdots M^*_{n+k}$, where $\left (\begin{array}{cc}1 & 1\\
0 & 1 \end{array} \right)$ is associated with mapping $R_{\alpha'}$ to $A_{\alpha'}$. Again applying (\ref{relation v_{n+k} and t_{n+k}}), we find that 
$$
D^*_k=\left [(\widetilde{M}^*_k)^{-1}(t_1),(\widetilde{M}^*_k)^{-1}(t_2)\right ] \times \left [(\widetilde{M}^*_k)^T(v_1),(\widetilde{M}^*_k)^T(v_2)\right ],
$$
from which it follows that
\begin{align*}
m(D^*_k)&=\iint_{D_k} \frac1{(1+tv)^2} \,dt\,dv\\
&=\log \frac{(1+(\widetilde{M}^*_k)^{-1}(t_2)(\widetilde{M}^*_k)^T(v_2))(1+(\widetilde{M}^*_k)^{-1}(t_1)(\widetilde{M}^*_k)^T(v_1))}{(1+(\widetilde{M}^*_k)^{-1}(t_2)(\widetilde{M}^*_k)^T(v_1))(1+(\widetilde{M}^*_k)^{-1}(t_1)(\widetilde{M}^*_k)^T(v_2))}\\
&=\log \frac{|\widetilde{M}^*_k|(t_2v_2+1)|\widetilde{M}^*_k|(t_1v_1+1)}{|\widetilde{M}^*_k|(t_2v_1+1)|\widetilde{M}^*_k|(t_1v_2+1)}\\
&=\log \frac{(1+t_2v_2)(1+t_1v_1)}{(1+t_2v_1)(1+t_1v_2)}=m(D).
\end{align*}
From this and (*) it follows that $\mu(D_k^*)=\mu(D)$ for any pair of $\alpha$-measurable sets $D_k$ and $D$. Since the dynamical system $(\Omega_1,\mathcal B, \mu_1, T_1)$ is ergodic, from our construction it easily follows that $(\Omega_\alpha,\mathcal B, \mu_\alpha, \mathcal T_\alpha)$ forms an ergodic system as well. Indeed, from our construction and regarding the above observations on $m$, it follows that for $\alpha \in [g^2,g)$, all $\Omega_\alpha$ are isomorphic to $\Omega_g$, so all dynamical systems $(\Omega_\alpha,\mathcal B, \mu_\alpha, \mathcal T_\alpha)$, $\alpha \in [g^2,1)$, will `inherit' ergodicity from $(\Omega_g,\mathcal B, \mu_g, \mathcal T_g)$. Obviously, for $\alpha \in [g,1)$ the dynamical systems $(\Omega_\alpha,\mathcal B, \mu_\alpha, \mathcal T_\alpha)$ are {\it{not}} isomorphic to $\Omega_1$, but from \cite{K} we know they are isomorphic to an induced transformation of $\Omega_1$, which is also clear from our construction; for $g<\alpha <1$, let $I_\alpha \subseteq \Omega_\alpha$ be invariant under $\mathcal T_\alpha$. We define $I_{\alpha,1}:=I_\alpha \cap A_\alpha$ and $I_{\alpha,2}:=I_\alpha \setminus I_{\alpha,1}$. Then $I_{\alpha,2} \subseteq \Omega_1$. Furthermore, putting $I_{1,1}:=\mathcal M^{-1}(I_{\alpha,1})$ (with $\mathcal M$ the map defined in Section \ref{The case alpha in (g,1]}), $I_{1,2}:=\mathcal T_1(I_{1,1})$, and $I_1:=I_{\alpha,2}\cup I_{1,1}\cup I_{1,2}$, then $I_1\subseteq \Omega_1$ is by construction $\mathcal T_1$-invariant, and therefore $I_1$ has $\mu_1$-measure $0$ or $1$. Since $m(I_{\alpha,1})=m(I_{1,1})=m(I_{1,2})$, it follows that $m(I_\alpha)=m(I_1)-m(I_{1,1})$. If $\mu_1(I_1)=0$, we are done. If $\mu_1(I_1)=1$, we see that $m(I_\alpha)=\log 2(\mu_1(I_1)-\mu_1(I_{1,1}))=\log 2(1-\mu_1(R_\alpha))=\log 2 - m(R_\alpha)=N_\alpha$, and we see that 
$$
\mu_\alpha(I_\alpha)=\frac1{N_\alpha}m(I_\alpha)=1;
$$
so if $I_\alpha \subseteq \Omega_\alpha$ is a measurable $\mathcal T_\alpha$-invariant set, we find that $\mu_\alpha(I_\alpha) \in \{0,1\}$, and we conclude that $\mathcal T_\alpha$ is ergodic. 

Writing $h(\mathcal T_\alpha)$ for the {\it {entropy}} of $\mathcal T_\alpha$, in \cite{Na} Nakada also proved that 
$$
h(\mathcal T_\alpha)=\begin{cases}\frac{\pi^2}{6\log(1+\alpha)},\quad \alpha \in (g,1]; \\
\frac{\pi^2}{6\log(G)},\quad \,\,\,\,\, \alpha \in [\tfrac12,g].
\end{cases}
$$
For $g<\alpha<1$ this result is essentially due to the fact that $\Omega_\alpha$ is isomorphic to an appropriate induced transformation of $(\Omega_1,\mathcal B, \mu_1, T_1)$ and Abramov's formula; see \cite{K} for explicit details. This explains the increase of entropy as $\alpha$ decreases from $1$ to $g$. In the previous sections, we have shown that for $\alpha \in (g^2,g]$, the loss of $R_\alpha$ is completely compensated by the addition of $A_\alpha$, from which it follows that the systems $(\Omega_\alpha,\mathcal B, \mu_\alpha, \mathcal T_\alpha)$ are metrically isomorphic to $(\Omega_g,\mathcal B, \mu_g, \mathcal T_g)$. We conclude that $h(\mathcal T_\alpha)=\pi^2/(6\log G), \alpha \in (g^2,g]$, a result previously obtained also in \cite{KSS}. For $\alpha$ smaller than $g^2$, the entropy map $h(\alpha)$ behaves behaves quite irregularly; see \cite{LM,NaNa} and Giulio Tiozzo's thesis \cite{T} for further details. 

\section{The $\alpha$-Legendre constant}
\label{legendre}

In (\cite{Na2}), Nakada obtains a very strong connection between two constants that play an important role in the study of continued fractions: the {\it{Legendre constant}} and the {\it{Lenstra constant}}.\\
The Legendre constant is associated with the theorem of Legendre, mentioned on page~\pageref{Legendre Theorem}, where we also introduced the $\alpha$-Legendre constant
 $$
L(\alpha):=\sup\{c>0: q^2\left |x-\tfrac pq \right |<c, \gcd(p,q)=1 \Rightarrow \tfrac pq = \tfrac{p_n}{q_n}, n \geq 0\}.
$$

The Lenstra constant is associated with the following conjecture by W.~Doeblin and Hendrik Lenstra on regular continued fractions, proved by Wieb Bosma, Hendrik Jager and Freek Wiedijk in \cite{BJW}:
\begin{Theorem}
For almost all $x$ and all $t \in [0,1]$, the limit 
$$
\lim_{n \to \infty} \frac1n \sum _{j=1}^\infty \{1\leq j \leq n: \theta_j \leq t\}
$$
exists, and is equal to the distribution function
$$
F(t)=\begin{cases}
\frac t{\log 2}, & 0 \leq t \leq 1/2,\\
\frac1{\log 2}(1-t+\log 2t), & 1/2 \leq t \leq 1.
\end{cases}
$$
\end{Theorem}

The number $1/2$, the upper bound of the interval where $F(t)$ is linear, is called the {\it{Lenstra constant of regular continued fractions}}; for $\alpha$-expansions the Lenstra constant can be defined similarly. \smallskip

Let $\tilde{g}:=\dfrac{g-2+\sqrt{g^2+4}}{2g}=0.57549\dots$\,. In 1988, Shunji Ito (\cite{I2}; see also \cite{K}) proved that 
$$
L(\alpha)=\begin{cases} \dfrac{\alpha}{\alpha+1},\quad & g \leq \alpha \leq 1;\\
1-\alpha, \quad & \tilde{g} \leq \alpha \leq g;\\
\dfrac{\alpha}{1+g\alpha}, \quad & \frac12 \leq \alpha \leq \tilde{g}.
\end{cases}
$$
In 2010, Nakada proved in \cite{Na2} that for a continued fraction algorithm the Legendre constant equals the Lenstra constant whenever the Legendre constant exists; one year later Natsui (\cite{Nat}) proved the existence of the $\alpha$-Legendre constant for $0<\alpha \leq 1$. She also extends Ito's result to 
$$
L(\alpha)=\frac{\alpha}{1+g\alpha},\quad \sqrt{2}-1 \leq \alpha \leq \tilde{g}.
$$
In this section we will in turn extend Natsui's result to $g^2\leq\alpha<\tilde{g}$. Moreover, we will show how we can use a result of Laura Luzzi and Stefano Marmi (\cite{LM}) to make some observations on $L(\alpha)$ for $\alpha < g^2$. \smallskip

\begin{figure}[!htb]
\begin{tikzpicture}[scale =10] 
\draw[draw=white,fill=gray!20!white] 
 plot[smooth,samples=100,domain=-.613:-.418] (\x,0) --
 plot[smooth,samples=100,domain=-.418:-.613] (\x,.0292);
 \draw[draw=white,fill=gray!20!white] 
 plot[smooth,samples=100,domain=-.613:-.418] (\x,.057) --
 plot[smooth,samples=100,domain=-.418:-.613] (\x,.077);
 \draw[draw=white,fill=gray!20!white] 
 plot[smooth,samples=100,domain=-.613:-.418] (\x,.082) --
 plot[smooth,samples=100,domain=-.418:-.613] (\x,.087);
  \draw[draw=white,fill=gray!20!white] 
 plot[smooth,samples=100,domain=-.613:-.418] (\x,.2) --
 plot[smooth,samples=100,domain=-.418:-.613] (\x,.346);
   \draw[draw=white,fill=gray!20!white] 
 plot[smooth,samples=100,domain=-.42:.387] (\x,0) --
 plot[smooth,samples=100,domain=.387:-.42] (\x,.346);
   \draw[draw=white,fill=gray!20!white] 
 plot[smooth,samples=100,domain=-.279:.387] (\x,.4) --
 plot[smooth,samples=100,domain=.387:-.279] (\x,.44); 
 \draw[draw=white,fill=gray!20!white] 
 plot[smooth,samples=100,domain=-.279:.387] (\x,.45) --
 plot[smooth,samples=100,domain=.387:-.279] (\x,.46);
 \draw [domain=-.613:.387,smooth,variable=\x] plot ({\x},{.9378/(1-.3126*\x)-.8});
\draw [domain=-.613:.387,smooth,variable=\x] plot ({\x},{0.96/(1-.32*\x)-.8});
\draw[black,fill=black] (-.42,.029) circle (.03ex);
\draw[black,fill=black] (.387,.345) circle (.03ex);
   \draw [dashed] (.387,0) -- (.387,.46);
             \draw (-.419,.087) -- (-.419,.2);
    \draw (-.613,0) -- (-.613,.0292);
       \draw (-.613,.2) -- (-.613,.346);
         \draw [dashed] (-.613,.0292) -- (-.419,.0292);
    \draw (-.279,.4) -- (.387,.4);
  \draw (-.279,.4) -- (-.279,.4399);
  \draw [dashed] (-.279,.4399) -- (.387,.4399); 
  \draw (-.419,.2) -- (-.613,.2);
      \draw (-.279,.45) -- (.387,.45);
      \draw (-.279,.45) -- (-.279,.46);
 \draw [dashed] (-.279,.46) -- (.387,.46);
     \draw (-.279,.441) -- (.387,.441);
      \draw (-.279,.441) -- (-.279,.4426);
 \draw [dashed] (-.279,.4426) -- (.387,.4426);
     \draw (-.279,.4429) -- (.387,.4429);
      \draw (-.279,.4429) -- (-.279,.4431);
 \draw [dashed] (-.279,.4431) -- (.387,.4431);
    \draw (-.613,.057) -- (-.419,.057);
    \draw (-.613,.057) -- (-.613,.077);
    \draw [dashed] (-.613,.077) -- (-.419,.077);
    \draw (-.419,.077) -- (-.419,.078);
         \draw (-.419,.079) -- (-.419,.082);
    \draw (-.613,.082) -- (-.419,.082);
  \draw (-.613,.082) -- (-.613,.087);
  \draw [dashed] (-.613,.087) -- (-.419,.087);
     \draw (-.419,.0292) -- (-.419,.057);
\draw [dashed] (-.419,0) -- (-.419,.0292);
        \draw [dashed] (-.613,.346) -- (.387,.346);
   \node at (-.613,-.025) {$\alpha-1$};
      \node at (-.425,-.025) {$\tfrac{1-3\alpha}{\alpha}$};
      \node at (.387,-.025) {$\alpha$};
 \node at (-.647,.028) {$\tfrac1{3+g}$};
 \node at (-.632,.2) {$\tfrac13$};
 \node at (-.632,.346) {$g^2$};
    \node at (-.3,.4) {$\tfrac25$};
      \node at (-.325,.463) {$\tfrac{2+g}{5+2g}$};
          \node at (.14,.13) {$\theta_{n-1}=\frac{\alpha}{1+g\alpha}$};
         \node at (-.12,.17) {$\theta_{n-1}>\frac{\alpha}{1+g\alpha}$};
         \node at (-.385,.028) {$M_\alpha$};
          \node at (.42,.345) {$N_\alpha$};
        \end{tikzpicture}
        \caption[The lower part of $\Omega_\alpha$, with $\alpha$ near $\tfrac{\sqrt{10}-2}3$.] {\label{fig: alpha with thetas}
The lower part of $\Omega_\alpha$, with $\alpha$ near $\tfrac{\sqrt{10}-2}3$.}
\end{figure}
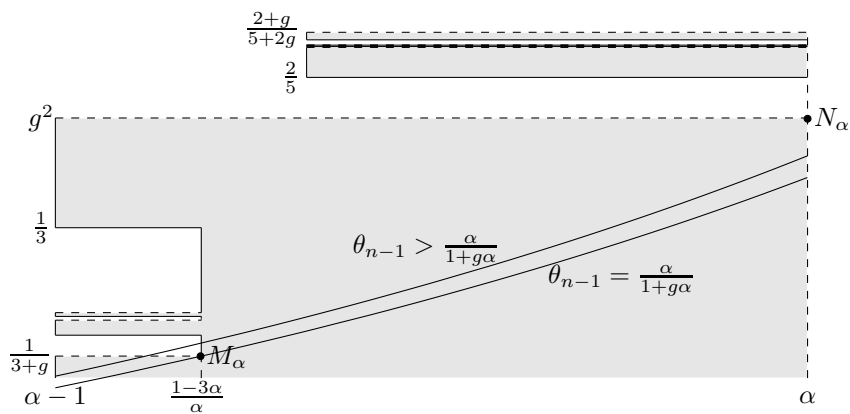

From Section \ref{The long way down to g^2} it is clear that it is both very hard and unmanageable to give $\Omega_\alpha$ explicitly for $\alpha\in [g^2,(\sqrt{10}-2)/3]$. We can, however, give some useful characteristics for this case with regard to $L(\alpha)$. The first thing we note is that $0\leq v < g$ for $(t,v) \in \Omega_\alpha$, where $\alpha\in [g^2,1/2]$. This follows from the fact that in our method of constructing $\Omega_\alpha$ from $\Omega_{\alpha'}$, with $\alpha <\alpha'$, any $-2$ that is introduced by insertion and singularisation is precedented by either a $3$ with minus sign or at least a $4$. Any second coordinate $v$ in $\Omega_\alpha$ will therefore be smaller than $1/(2-g^2)=g$. From this it follows that the bottom side of the lowermost cove in Figure \ref{fig: alpha with thetas} is on the line $v=1/(3+g)$; see our calculations in Section \ref{The case alpha in (tfrac{sqrt{10}-2}3,sqrt{2}-1]}. Figure \ref{fig: alpha with thetas} is an example of $\Omega_\alpha$ with $\alpha$ near $(\sqrt{10}-2)/3$, and is more or less the same as Figures \ref{fig: omega2/5} and \ref{fig: omega(sqrt10-2)/3}. The difference between Figure  \ref{fig: alpha with thetas} and the other two figures is the addition of curves associated with $\theta_{n-1}=\alpha/(1+g\alpha)$ and $\theta_{n-1}>\alpha/(1+g\alpha)$, the choice of which we will explain now.\smallskip

 Recall that in the introduction we defined $\theta_n(x):=q_n^2|x-p_n/q_n|$ and mentioned the equations $\theta_{n-1}=v_n/(1+t_nv_n)$ and $\theta_n=(\varepsilon_nt_n)/(1+t_nv_n)$ (omitting the suffix `($x$)'). In \cite{K}, $L(\alpha)$ is found by determining the largest $C$ for which the curve $\theta_n=C$ contains no points in what is called a singularisation area: indeed, this area is associated with convergents that are replaced (when $1/2 < \alpha \leq g$) or lost (when $g<\alpha<1$ or $\alpha \leq 1/2$), as we have shown in Sections \ref{The case alpha in (g,1]} through \ref{The long way down to g^2}. Note that this $C$ is also the Lenstra constant for the $\alpha$-expansion under consideration. In the current case ($g^2 \leq \alpha <\sqrt{2}-1$), we find it convenient to consider curves $\theta_{n-1}=C$ rather than $\theta_n=C$, if only for the numerator of $\theta_n=\varepsilon_nt_n/(1+t_nv_n)$, yielding two branches that we would have to take into account separately. Since curves associated with the equation $\theta_{n-1}=C$, $C$ small enough, will be in the lower part of $\Omega_\alpha$, a first guess would be that $L(\alpha)$ is determined by either the points $M_\alpha:=((1-3\alpha)/\alpha,1/(3+g))$ or $N_\alpha:=(\alpha, g^2)$. No matter how intricate $\Omega_\alpha$ becomes for $\alpha \leq \sqrt{2}-1$, from our earlier obeservations it follows that these $v$-coordinates $1/(3+g)$ and $g^2$ do not change as $\alpha$ decreases to $g^2$. In Figure \ref{fig: alpha with thetas} we see how curves associated with $\theta_{n-1}>\alpha/(1+g\alpha)$ contain points in the cove left and above $M_\alpha$. 
We find $\theta_{n-1}(M_\alpha)=\alpha/(1+g\alpha)$, while $\theta_{n-1}(N_\alpha)=1/(G^2+\alpha)$. Since
$$
\min \{\alpha/(1+g\alpha),1/(G^2+\alpha)\} =\alpha/(1+g\alpha)
$$
on $[g^2,\sqrt{2}-1)$, we conclude that 
$$
L(\alpha)=\frac{\alpha}{1+g\alpha},\quad g^2 \leq \alpha < \sqrt{2}-1,
$$
hence extending Natsui's result. 

\begin{figure}[!htb]
$$
\begin{tikzpicture}[scale =12] 
\draw[draw=white,fill=gray!20!white] 
 plot[smooth,samples=100,domain=-.75:.25] (\x,0) --
 plot[smooth,samples=100,domain=.25:-.75] (\x,.2087);
 \draw[draw=white,fill=gray!20!white] 
 plot[smooth,samples=100,domain=0:.25] (\x,.208) --
 plot[smooth,samples=100,domain=.25:0] (\x,.264);
\draw [domain=-.75:.25,smooth,variable=\x] plot ({\x},{.2087/(1-.2087*\x)});
\draw[black,fill=black] (0,.2087) circle (.03ex);
 \draw (-.75,0) -- (.25,0);
 \draw (-.75,0) -- (-.75,.2087);
 \draw (-.75,.2087) -- (0,.2087);
 \draw (0,.2087) -- (0,.264);
 \draw (0,.264) -- (.25,.264);
 \draw (.25,0) -- (.25,.264);
 \draw [dashed] (0,0) -- (0,.2087);
 \node at (-.35,.17) {$\theta_{n-1}=L(1/4)$};
  \node at (-.79,.21) {$\tfrac{5-\sqrt{21}}2$};
  \node at (-.04,.264) {$\tfrac{\sqrt{21}-3}6$};
  \node at (0,-.02){$0$};
  \node at (-.75,-.02) {$-3/4$};
       \node at (.25,-.02) {$1/4$};
   \end{tikzpicture}
 $$
 \caption[The bottom part of $\Omega_{1/4}$.] {\label{fig: The bottom part of omega.25}
The bottom part of $\Omega_{1/4}$.}
\end{figure}
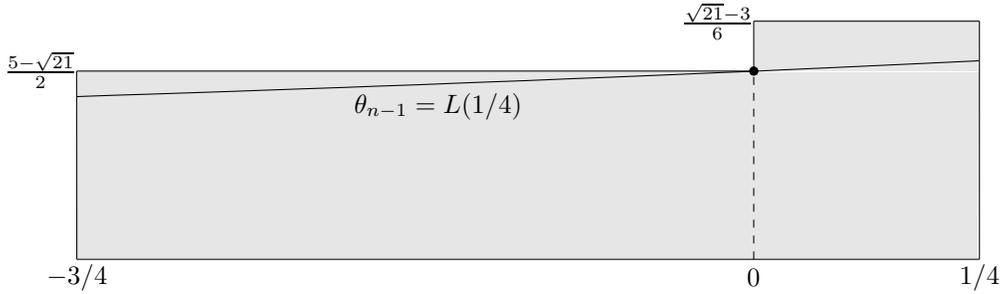

For various reasons, this approach cannot be extended to values of $\alpha<g^2$. For one thing, $g$ will cease to be the upper bound for the values of $v$, since insertions and singularisations in the case $\alpha<g^2$ involve strings of more than one $-2$, yielding values of $v$ larger than $g$. This implies that points have to be removed with $v$-values under $1/(3+g)$. Unfortunately, we cannot say much more about $\Omega_\alpha$ for general $\alpha$ smaller than $g^2$. Still, a result of Luzzi and Marmi (see \cite{LM}, pages 27 and 28) enables us to at least find how $L(\alpha)$ evolves for these rational values of $\alpha$. In \cite{LM}, Luzzi and Marmi use the aforementioned strings of partial quotients $2$ with minus sign to give a description of $\Omega_\alpha$ for $\alpha=1/r$, $r \in \N_{\geq 3}$. Although their description is not very explicit, it is not hard to find that it yields that $\Omega_{1/r}$ consists of the bottom part 
$$
[-1+1/r,0]\times[0,(r+1-\sqrt{(r+1)^2-4})/2]\cup[0,1/r]\times[0,(\sqrt{(r+1)^2-4}-r+1)/(2r-2)]
$$ 
and a myriad of vertically disjunct rectangles above it; see Figure \ref{fig: The bottom part of omega.25} for the bottom part. Taking the approach of considering curves $\theta_{n-1}=C$ as we did above, we see that the Lenstra constant for $\alpha = 1/r$, $L(1/r)$ is determined by the vertex $(0,(r+1-\sqrt{(r+1)^2-4})/2)$ (which is the point $(0,(5-\sqrt{21})/2)$ in Figure \ref{fig: The bottom part of omega.25}). Thanks to Nakada's result in \cite{Na2} it follows that
\begin{equation}\label{L(alpha)}
L\left (\tfrac1r \right )=\frac{r+1-\sqrt{(r+1)^2-4}}2.
\end{equation}
Although this formula is different from the one we found for $g^2 \leq \alpha \leq 1/2$, we remark that $L(1/2)=g^2$ fits both formulas. What's more, we suspect that with (\ref{L(alpha)}) a sharpening is possible of Natsui's result on the boundaries of $L(\alpha)$ in \cite{Nat}, which is: 
\begin{Theorem}
Let $r \in \N$. If $\frac1{r+1}<\alpha<\frac1r$, then $\frac1{r+2}<L(\alpha)<\frac1r$, with $r \in \N_{\geq 3}$.
\end{Theorem}
Since (\ref{L(alpha)}) yields 
\begin{equation}
L(\tfrac1r)=\frac{r+1-\sqrt{(r+1)^2-4}}2=\frac1{r+1-L(\tfrac1r)},
\end{equation}
we suspect that in fact (if $1/(r+1)<\alpha<1/r$)
$$
\frac1{r+2}<L(\alpha)<\frac1{r+1-L(\tfrac1r)}.
$$

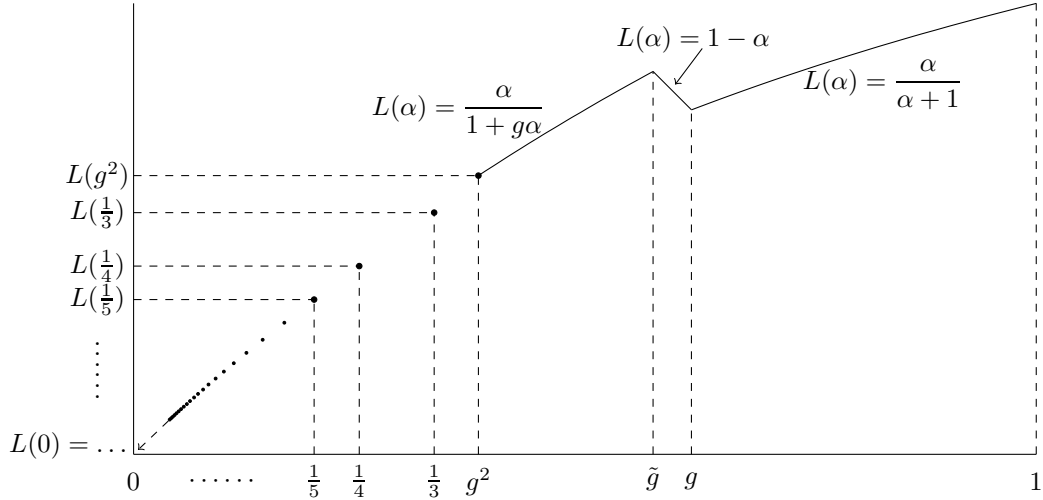
\begin{figure}[!htb]
$$
\begin{tikzpicture}[scale =12] 
\draw [domain=.618:1,smooth,variable=\x] plot ({\x},{\x/(\x+1)});
\draw [domain=.5755:.618,smooth,variable=\x] plot ({\x},{1-\x});
\draw [domain=.382:.5755,smooth,variable=\x] plot ({\x},{\x/(.618*\x+1)});
\draw[black,fill=black] (.382,.309) circle (.02ex);
\draw[black,fill=black] (.333,.268) circle (.02ex);
\draw[black,fill=black] (.25,.2087) circle (.02ex);
\draw[black,fill=black] (.2,.17157) circle (.02ex);
\draw[black,fill=black] (.167,.1459) circle (.01ex);
\draw[black,fill=black] (.143,.127) circle (.01ex);
\draw[black,fill=black] (.125,.1125) circle (.01ex);
\draw[black,fill=black] (.111,.101) circle (.01ex);
\draw[black,fill=black] (.1,.0917) circle (.01ex);
\draw[black,fill=black] (.0909,.084) circle (.01ex);
\draw[black,fill=black] (.083,.0774) circle (.01ex);
\draw[black,fill=black] (.077,.0718) circle (.01ex);
\draw[black,fill=black] (.0714,.067) circle (.01ex);
\draw[black,fill=black] (.067,.063) circle (.01ex);
\draw[black,fill=black] (.0625,.059) circle (.01ex);
\draw[black,fill=black] (.059,.0557) circle (.01ex);
\draw[black,fill=black] (.0556,.0528) circle (.01ex);
\draw[black,fill=black] (.0526,.0501) circle (.01ex);
\draw[black,fill=black] (.05,.0477) circle (.01ex);
\draw[black,fill=black] (.0476,.0456) circle (.01ex);
\draw[black,fill=black] (.0455,.0436) circle (.01ex);
\draw[black,fill=black] (.0435,.0417) circle (.01ex);
\draw[black,fill=black] (.0417,.0401) circle (.01ex);
\draw[black,fill=black] (.04,.0385) circle (.01ex);
\draw [dashed] (.2,0) -- (.2,.17157);
\draw [dashed] (0,.17157) -- (.2,.17157);
\draw [dashed] (.25,0) -- (.25,.2087);
\draw [dashed] (0,.2087) -- (.25,.2087);
\draw [dashed] (.333,0) -- (.333,.26795);
\draw [dashed] (0,.26795) -- (.333,.26795);
\draw [dashed] (.382,0) -- (.382,.309);
\draw [dashed] (0,.309) -- (.382,.309);
\draw [dashed] (.5755,0) -- (.5755,.4245);
\draw [dashed] (.618,0) -- (.618,.382);
\draw [dashed] (1,0) -- (1,.5);
\draw [->] (.63,.45) -- (.6,.41);
\draw [->] [dashed] (.0385,.037) -- (.005,.005);
\node at (1,-.03) {$1$};
\node at (.618,-.03) {$g$};
\node at (.5755,-.027) {$\tilde{g}$};
\node at (.382,-.03) {$g^2$};
\node at (.333,-.03) {$\tfrac13$};
\node at (.25,-.03) {$\tfrac14$};
\node at (.2,-.03) {$\tfrac15$};
\node at (-.04,.309) {$L(g^2)$};
\node at (-.04,.268) {$L(\tfrac13)$};
\node at (-.04,.2087) {$L(\tfrac14)$};
\node at (-.04,.17157) {$L(\tfrac15)$};
\node at (-.07,.01) {$L(0)=\dots$};
\node at (-.04,.12) {$\vdots$};
\node at (-.04,.085) {$\vdots$};
\node at (.1,-.03) {$\cdots\cdots$};
\node at (0,-.03) {$0$};
\node at (.83,.41) {$L(\alpha)=\dfrac{\alpha}{\alpha+1}$};
\node at (.36,.38) {$L(\alpha)=\dfrac{\alpha}{1+g\alpha}$};
\node at (.62,.46) {$L(\alpha)=1-\alpha$};
\draw (0,0) -- (1,0);
\draw (0,0) -- (0,.5);
\end{tikzpicture}
$$
\caption[The graph of $L(\alpha)$ as it is currently known.] {\label{fig: L(alpha)}
The graph of $L(\alpha)$ as it is currently known.}
\end{figure}
Considering this figure, one would expect that $L(0)=0$; this is actually true, and we will return on this in a forthcoming note.

\section*{Acknowledgements}
The first author of this paper is very thankful to the Dutch organisation for scientific research NWO, that funded his research for this paper with a so-called {\it{Promotiebeurs voor Leraren}}, with grant number 023.003.036.\\

\end{document}